\documentclass[envcountsame,envcountchap]{svmono_temp}

\usepackage{makeidx}         % allows index generation
\usepackage{graphicx}        % standard LaTeX graphics tool
\usepackage{multicol}        % used for the two-column index
\usepackage[bottom]{footmisc}% places footnotes at page bottom
\usepackage{amssymb}

\def \det{\mbox{{\rm det$\,$}}}
\def \span{\mbox{{\rm span$\,$}}}
\def \col{\mbox{{\rm col$\,$}}}

\def \rank{\mbox{{\rm rank$\,$}}}

\def \blkdiag{\mbox{\rm blkdiag$\,$}}
\def \ad{\mbox{{\rm ad$\,$}}}

\def \RR{{\mathbb{R}\,\!}}

\newtheorem{theo}{Theorem}[section]

\newtheorem{defi}{Definition}[section]
\newtheorem{rema}{Remark}[section]
\newtheorem{coro}{Corollary}[section]
\newtheorem{assumption}{Assumption}

\newcommand{\ben}{\begin{enumerate}}
\newcommand{\een}{\end{enumerate}}

\newcommand{\mm}[1]{\left[\matrix{#1}\right]}
\newcommand{\be}{\begin{equation}} % \be is same as \eq
\newcommand{\ee}{\end{equation}}   % \ee is same as \en

\def \squ{\hfill \mbox{$\Box$}}

\def \rma{{\rm a}}
\def \rmb{{\rm b}}
\def \rmc{{\rm c}}
\def \rmd{{\rm d}}
\def \rme{{\rm e}}

\def \rmi{{\rm i}}

\def \rmo{{\rm o}}
\def \rmp{{\rm p}}

\def \rms{{\rm s}}

\def \rmI{{\mbox{\tiny {\rm I$\,$}}}}

\def \rmO{{\mbox{\tiny {\rm O$\,$}}}}

\def \rmT{{\mbox{\tiny {\rm T}}}}

\def \bbR{{\mathbb{R}\,\!}}

\def \calA{{\cal A}}
\def \calB{{\cal B}}
\def \calC{{\cal C}}
\def \calD{{\cal D}}
\def \calE{{\cal E}}

\def \calK{{\cal K}}
\def \calL{{\cal L}}

\def \calO{{\cal O}}

\def \calX{{\cal X}}
\def \calY{{\cal Y}}
\def \calZ{{\cal Z}}

\makeindex             % used for the subject index
\usepackage[numbers,sort&compress]{natbib}

\usepackage{color}

\begin{document}

\author{Xinmin Liu\\
%Zongli Lin\\
%**** ***
%Ben M. Chen
}
\title{
%Geometric Structures of Control Systems\\
%Structures of Dynamical Systems\\
%Structures of Multivariable Systems\\
%Structures of Differential Equation Systems\\
%Structures of Control Systems\\
%Nonlinear Multivariable Control\\
%Nonlinear Control Systems\\
On Differential Geometric Approach to Nonlinear Systems Affine in Control
%{\small SPIN Springer's internal project number, if known}
%Dynamical Control Systems
}

\subtitle{
%A Structural Analysis Approach\\
%A Structural Approach\\
%\A Differential Geometric Approach\\
%From Linear to Nonlinear
} \maketitle

\frontmatter%%%%%%%%%%%%%%%%%%%%%%%%%%%%%%%%%%%%%%%%%%%%%%%%%%%%%%

\include{pref}

\tableofcontents

\mainmatter%%%%%%%%%%%%%%%%%%%%%%%%%%%%%%%%%%%%%%%%%%%%%%%%%%%%%%%
%\include{part}
%\include{newideas}

\chapter{Introduction}

\section{Literature Overview}

The note is concerned with nonlinear systems affine in control described by ordinary differential equations of the following form,
\be
\label{nonsys_nfn9}
\left\{\begin{array}{rcl}\dot x &=& f(x)+g(x)u,\cr
                           y &=& h(x),\end{array}\right.
\ee
where $x$, $u$ and $y$ denote the system state, input and output, respectively.

The differential geometric approach to nonlinear control has been proven to be a powerful tool
to deal with fundamental questions in the state space formulation of nonlinear control systems.
Elliott \cite{elliott1997bookreview}, and Nijmeijer \& Schaft \cite{nind90}
had good reviews on the development of differential geometric control theory.
In the 1960's, the popularity of Pontryagin's Maximal Principle led to the need to
understand controllability, and the researchers
realized that some technical assumptions about the nonlinear systems,
such as smoothness and analyticity, could lead to a general mathematical approach.
Hermann \cite{hermann1963accessibility,hermann1968differential,hermann1977nonlinear} studied controllability with methods based on
vector fields and differential forms, which is analogous to Kalman's
criterion for linear systems.
In the early 1970's Brockett, Boothby, Elliott, et al.
were promoting the use of Lie algebra methods to study controllability.
Brockett \cite{brockett1979feedback,brockett1983asymptotic} and Willems
also considered systems invariants equivalent by coordinate change and a class of feedback transformations.
Isidori, Krener, Gori-Giorgi \& Monaco\cite{isidori1981nonlinear}, and Hirschorn \cite{hirschorn1981mathcal}
used the concept of controlled invariant distribution for the solving of the problem of decoupling problems.
Many concepts of differential geometric control on nonlinear systems are indeed
the generalization of concepts of geometric control of linear systems.
Wonham and Morse \cite{mosi73,wolm79,wonham1970decoupling,wonham1972feedback}
and Basile and Marro \cite{basile1969controlled,basile1992controlled}
developed a systematic geometric approach to solving the problems of
pole placement, noninteracting control, disturbance decoupling, and
regulation. This approach depends on global linear space
structure. Isidori \cite{isnc95,isnc99} generalized a local approach of this nature
to nonlinear control problems.
He brought the geometry and
Volterra series methods together and used them appropriately
for stabilization, regulation, disturbance decoupling, noninteracting
control, tracking and regulation \cite{elliott1997bookreview}.

The nonlinear analogues of linear system structural properties,
such as relative degree (or infinite zero structure), zero dynamics
(or finite zero structure) and invertibility properties, have played
critical roles in recent literature on the analysis and control
design for nonlinear systems
The normal forms that are associated with these structural
properties, along with some basic tools,
have enabled many major breakthroughs in nonlinear control theory.

A single input single output system
has a relative degree, if the system can be reduced to
the zero dynamics cascaded with a clean chain of integrators linking
the input to the output.
Here by {\em clean} we mean that no other signal enters the middle of the chain.
This structural feature is extended to nonlinear systems with more
than one input/output pair. For a square invertible nonlinear
system, the notion of vector relative degree was introduced in \cite{byas91,saac89},
and the systems can be transferred into the zero dynamics connecting to
{\em clean} chains of integrators.

The clean chains of integrators are called the prime form in \cite{mosi73} for linear systems,
and the necessary and sufficient geometric
conditions for the existence of prime forms for nonlinear systems is were established \cite{marino1994equivalence}.
The lengths of chains of integrators are the nonlinear extension of infinite zeros.
However, vector relative degree is a rather restrictive structural property that
not even all square invertible linear systems, with the freedom
of choosing coordinates for the state, output and input spaces and state feedback,
could possess.

A major generalization of the normal form representations was made in
\cite{byls88,isnc95,isnc99,scgn99}, where
square invertible systems are considered.
With the assumption that the rank of certain matrices are constant
on a sequence of nested submanifolds, or
with some stronger assumptions \cite{scgn99,isnc99},
the nonlinear systems can be represented
by the zero dynamics cascaded with chains of integrators.
Note that chains of integrators here need not to be {\em clean}.
Interconnection between chains of integrators are allowed.
This greatly enlarges the class of nonlinear systems that normal forms can represent.
But in these normal forms, the lengths of chains of integrators are
no longer the nonlinear extension of infinite zeros.
The applications of these normal forms in solving the problem of
asymptotic stabilization, disturbance decoupling, tracking and regulation
can be found in \cite{isnc95,isnc99} and the references therein.

In the note,
we make an attempt to study structural properties of affine nonlinear systems.
We will develop a constructive algorithm to represent nonlinear systems in normal forms.
In the special
case when the system is square and invertible, our normal forms take forms
similar to those in \cite{isnc95,isnc99,scgn99}, but with an additional
property that allows the normal forms to reveal the nonlinear extension of infinite
zeros of linear systems. In addition, our algorithms require fewer assumptions,
can apply to general nonlinear systems that are not necessarily square, and
can explicitly show invertibility structures of the systems.
We will also study the applications of these new normal forms to solving the problems of
global stabilization, semi-global stabilization and disturbance attenuation.

\section{Note Outline}

The note focuses on
the differential geometric approach to the study of nonlinear systems that are affine in control.
We first develop normal forms for nonlinear system affine in control.
Based on these normal forms, we then address the problems of global stabilization,
semi-global stabilization and disturbance attenuation.
The results presented are based on the works
\cite{
chu2002numerical,
liu2003problem,
liu2005linear,
liu2006symbolic,
liu2007further,
chu2007numerical,
chen2008interconnection,
lisn08,
liu2009assignment,
liu2009semi3,
liu2010disturbance,
liu2011normal,
liu2012backstepping,
liu2012grid,
liu2013dynamic,
belcher2014development,liu2015robotic,belcher2016spatial,
liu2016modeling,
liu2017design,
liu2017use}.

%\cite{belcher2014development,liu2015robotic1,liu2015robotic,belcher2016spatial,liu2016constrained,liu2017use}
%\cite{belcher2014development,liu2015robotic1,liu2015robotic,belcher2016spatial,liu2017use}

The note can naturally be divided into three parts.

The first part is Chapter 2, which presents a brief introduction to the differential geometric concepts for use in the note.
It includes the fundamental concepts of manifolds, submanifolds, tangent vectors, vector fields,
and distributions.

The second part is Chapter 3.
In this chapter, we propose constructive algorithms for decomposing a
nonlinear system that is affine in control but otherwise general.
These algorithms require
modest assumptions on the system and apply to
general multiple input multiple output systems that do not necessarily
have the same number of inputs and outputs. They lead
to various normal form representations and
reveal the structure at infinity, the
zero dynamics and the invertibility properties, all of which
represent nonlinear extensions of relevant linear system structural
properties of the system they represent.

The third part of the note consists of Chapters 4, 5 and 6.
They contain some applications of the structural decomposition developed in Chapter 3.
In Chapter 4, we exploit the properties of such
a decomposition for the purpose of solving the stabilization problem.
In particular, this structural
decomposition simplifies
the conventional backstepping design and motivates new backstepping design
procedures that are able to stabilize some systems on which the conventional
backstepping is not applicable.

In Chapter 5, we exploit the properties of such
a decomposition for the purpose of solving the semi-global stabilization problem
for minimum phase nonlinear systems without vector relative degrees.
By taking advantage of the special structure of the decomposed system, we first apply
the low gain design to the part of system that possesses linear dynamics.
The low gain design results in an augmented zero dynamics that is locally stable
at the origin with a domain of attraction that can be made arbitrarily large by
lowering the gain. With this augmented zero dynamics, the backstepping design procedure is
then applied to achieve semi-global stabilization of the overall system.

Chapter 6 considers the problems of disturbance attenuation and almost disturbance decoupling,
which have played a central role in control theory.
By employing the structural decomposition of multiple
input multiple output nonlinear systems and the backstepping
procedures that we have developed, we show that these two problems can be
solved for a larger class of nonlinear systems.

Finally, Chapter 7 is the conclusions to the note, and some topics for the future research
are also mentioned.

\chapter{Manifolds, Tangent Vectors, Vector Fields, Distributions}

The chapter recalls some basic concepts and facts of differential geometry that
will be used in the following chapters.
The detail can be found in \cite{agrachev2004control,bullo2005geometric,doolin1990introduction,isnc95,nind90,waner2005introduction}.

Differential geometry is a discipline on curves and surfaces.
It studies the functions that define curves and surfaces, and the transformations between the coordinates
that are used to specify curves and surfaces. It also treats the differential relations that put
pieces of curves or surfaces together.

\section{Manifolds}

A manifold is a mathematical space that on a small enough scale resembles the Euclidean space of a specific dimension.
A line and a circle are one-dimensional manifolds, and a plane and the surface of a ball are two-dimensional
manifolds.
Although manifolds resemble Euclidean spaces near each point locally, the global structure of a manifold is more complicated.
A chart of a manifold is an invertible map between a subset of the manifold and
the Euclidean space such that both the map and its inverse preserve the desired structure.
The description of most manifolds requires more than one chart.
A specific collection of charts which covers a manifold is called an atlas.
Charts in an atlas may overlap and a single point of a manifold may be represented in several charts.
Given two overlapping charts, a transition map can be defined which goes from an open ball in Euclidean space
to the manifold and then back to another open ball in Euclidean space.

Topological spaces are structures that define convergence, connectedness, and continuity.
A topological space is a set $X$ together with $\Omega$, a collection of subsets of $X$, satisfying the following axioms:
  \ben
   \item[1)] The empty set and X are in $\Omega$.
   \item[2)] The union of any collection of sets in $\Omega$ is also in $\Omega$.
   \item[3)] The intersection of any finite collection of sets in $\Omega$ is also in $\Omega$.
  \een
The collection $\Omega$ is called a topology on $X$. The elements of $X$ are usually called points.
It is customary to require that the space be Hausdorff and second countable.

A topological manifold is a topological space locally homeomorphic to a Euclidean space,
which means that every point has a neighborhood for which there exists a homeomorphism
(a bijective continuous function whose inverse is also continuous)
mapping that neighborhood to a Euclidean space.

A differentiable manifold is a topological manifold that allow one to do differential calculus.
The primary object of study in differential calculus is the derivative. We now consider the derivative
of a function $f$ with domain an open subset $U$ of $\RR^n$ and with range in $\RR^m$.
The function $f$ is differentiable at $x\in U$ if there is a linear map $A(x)$, a Jacobian matrix,
from $\RR^n$ to $\RR^m$ such that
\[
\lim_{|h|\rightarrow 0} \frac{|f(x+h)-f(x)-A(x)h|}{|h|}=0.
\]
Then $A(x)$ is called the derivative of $f$.
A $\calC^k$ manifold is a differential manifold with an atlas whose transition maps are all $k$-times continuously differentiable.

A smooth manifold ( $\calC^\infty$ manifold ) is a differentiable manifold for which all the transition maps are smooth. That is, derivatives of all orders exist.
An analytic manifold, ( $\calC^\omega$ manifold ) is a smooth manifold with the additional condition that each transition map is analytic: the Taylor expansion is absolutely convergent on some open ball.

Consider a topological space $(X, \Omega)$. Suppose that for any $p\in X$, there exists an open set $U\in \Omega$ with $p\in U$, and a bijection
$\phi$ mapping $U$ onto an open subset of $\RR^n$,
\[
\phi : U \rightarrow \phi(U)\subset \RR^n.
\]
The grid defined on $\phi(U)\subset \RR^n$ is transforms into a grid on $U$. A coordinate chart is the pair $(U,\phi)$.
The map $\phi$ can be represented as a set $(\phi_1,\phi_2,\cdots,\phi_n)$ and $\phi_i: U\rightarrow \RR$ is called
the $i$-th coordinate function. The $n$-tuple of real numbers $(\phi_1(p),\phi_2(p),\cdots,\phi_n(p))$ is called
the set of local coordinates of $p$ in the coordinate chart $(U,\phi)$.

For example, the helix represented by
\begin{eqnarray*}
z_1&=&\cos x_1\\
z_2&=&\sin x_1\\
z_3&=& x_1
\end{eqnarray*}
is a smooth path embedded in Euclidean space $\RR^3$.
It is 1-dimensional smooth manifold.
The parameters $x_1$ is local coordinate, and $z_1$, $z_2$ and $z_3$ are global coordinates
or ambient coordinates.

The sphere $z_1^2+z_2^2+z_3^2=1$ is a smooth surface embedded in Euclidean space $\RR^3$. It is
2-dimensional smooth manifold.
Using spherical polar coordinates, the sphere is represented by
\begin{eqnarray*}
z_1&=&\sin x_1 \cos x_2\\
z_2&=&\sin x_1 \sin x_2\\
z_3&=&\cos x_1
\end{eqnarray*}
For points other than $(0,0,\pm 1)$,
\be
\left.\begin{array}{lll}\label{example_chart}
x_1&=&\arccos z_3\\
x_2&=&
\left\{\begin{array}{ll}\arccos (\frac{z_1}{\sqrt{z_1^2+z_2^2}})& \mbox{if}\;\; z_2\geq 0 \cr
                        2\pi-\arccos (\frac{z_1}{\sqrt{z_1^2+z_2^2}})& \mbox{if}\;\; z_2< 0 .\end{array}\right.
\end{array}
\right.
\ee
The chart of the sphere is the pair of functions in (\ref{example_chart}).
The parameters $x_1$ and $x_2$ are called local coordinates, while $z_1$, $z_2$ and $z_3$ are called global coordinates
or ambient coordinates.
The ambient coordinates are
superfluous data that often have nothing to do with the problem at hand.
It is a tremendous advantage to be able to work with manifolds, without the excess baggage of such an ambient
space.

Let $(U,\phi)$ and $(V,\varphi)$ be two coordinate charts on a manifold $N$ with $U\cap V \neq 0$.
The coordinates transformation on $U\cap V$
\[
\varphi \circ \phi^{-1} : \phi(U\cap V) \rightarrow \varphi(U\cap V)
\]
transfers the set of the local coordinate $(\phi_1(p),\phi_2(p),\cdots,\phi_n(p))$ to
the set of the local coordinate $(\varphi_1(p),\varphi_2(p),\cdots,\varphi_n(p))$. Two coordinate charts
$(U,\phi)$ and $(V,\varphi)$ are $\calC^\infty$-compatible if $\varphi \circ \phi^{-1}$ is smooth ($\calC^\infty$),
{\em i.e.,} $\varphi \circ \phi^{-1}$ is a diffeomorphism.

The set $(\phi_1(p),\phi_2(p),\cdots,\phi_n(p))$ can be represented as an $n$-vector
$x=\col\{x_1,x_2,\cdots,x_n\}$, and the set $(\varphi_1(p),\varphi_2(p),\cdots,\varphi_n(p))$ as \newline $y=\col\{y_1,y_2\cdots,y_n\}$.
Therefore, the coordinate transformation $\varphi\circ \phi^{-1}$ can be represented as
\[
y=\pmatrix{
y_1(x_1,x_2,\cdots,x_n)\cr
y_2(x_1,x_2,\cdots,x_n)\cr
\vdots\cr y_n(x_1,x_2,\cdots,x_n)}=y(x).
\]
and $\phi\circ \varphi^{-1}$ as
\[
x=\pmatrix{
x_1(y_1,y_2,\cdots,y_n)\cr
x_2(y_1,y_2,\cdots,y_n)\cr
\vdots\cr x_n(y_1,y_2,\cdots,y_n)}=x(y).
\]

A $\calC^\infty$ atlas on a manifold $N$ is a collection $\calA=\{ (U^i, \phi^i) : i\in I\}$ of pairwise $\calC^\infty$-compatible coordinate charts
with $\cup_{i\in I} U^i = N$.
An atlas is complete if not properly contained in any other atlas.
A smooth manifold is a manifold equipped with a complete $\calC^\infty$ atlas.

Let $N$ and $M$ be manifolds of dimension $n$ and $m$,
$(U,\phi)$ and $(V,\varphi)$ be coordinate charts on the manifolds $N$ and $M$, respectively.
$F: N \rightarrow M$ is a mapping. The mapping
\[
\check{F}=\varphi \circ F \circ \phi^{-1}
\]
is called an expression of $F$ in local coordinates.

Let $N$ and $M$ be smooth manifolds of dimension $n$.
A mapping $F: N\rightarrow M$ is a smooth mapping if for each $p\in N$ there exist
coordinate charts $(U,\phi)$ of $N$ and $(V,\varphi)$ of $M$, with $p\in U$ and $F(p)\in V$,
such that the expression of $F$ in local coordinates is $\calC^\infty$.

Let $N$ and $M$ be smooth manifolds of dimension $n$.
A mapping $F: N\rightarrow M$ is a diffeomorphism if $F$ is bijective and both $F$ and $F^{-1}$
are smooth mappings. Two manifolds $N$ and $M$ are diffeomorphic if
there exists a diffeomorphism $F: N\rightarrow M$.

\section{Submanifolds}

Let $N$ be a smooth manifold of dimension $n$. A non-empty open set $V\subset N$
is itself a smooth manifold of dimension $m$ with
coordinate charts obtained by restricting the coordinate charts for $N$ to $V$. $V$ is called
an open submanifold of $N$.

Let $N$ be a smooth manifold of dimension $n$. A subset $N'$ of $N$ is an embedded submanifold of dimension
$m<n$ if and only if for each $p\in N'$ there exists a cubic coordinate chart $(U,\phi)$ of $N$,
with $p\in U$, such that $U\cap N'$ coincides with an $n$-dimensional slice of $U$ passing through $p$.

Let $F : N\rightarrow M$ be a smooth mapping of manifolds.
$F$ is an immersion if $\rank(F)=\dim(N)$ for all $p\in N$.
$F$ is an univalent immersion if $F$ is an immersion and is injective.
$F$ is an embedding if $F$ is an univalent immersion and the topology induced on $F(N)$
by the one of $N$ coincides with the topology of $F(N)$ as a subset of $M$.

The image $F(N)$ of a univalent immersion is called an immersed submanifold of $M$.
The image $F(N)$ of an embedding is called an embedded submanifold of $M$.

Let $F: N\rightarrow M$ be an immersion. For each $p\in N$ there exists a neighborhood $U$ of $p$
such that the restriction of $F$ to $U$ is an embedding.

For $F: N\rightarrow M$, let $M'=F(N)$ and $F': N\rightarrow M'$. If the topology of $M'$ is the one induced by
one of $N$, $F'$ is a homeomorphism.  Any coordinate chart $(U, \phi)$ of $N$ induces a coordinate chart $(V,\varphi)$
of $M'$, i.e.,
\[
V=F'(U),\quad \varphi=\phi\circ (F')^{-1}.
\]
The smooth manifold $M'$ is diffeomorphic  to the smooth manifold $N$.

\section{Tangent Vectors}

Let $N$ be a smooth manifold of dimension $n$, and $x$ be a point in $N$.
A tangent space is a real vector space that tangentially pass through the point $x$.
The elements of the tangent space are called tangent vectors at $x$.

All the tangent spaces can be ``glued together" to form a new differentiable manifold of twice the dimension, the tangent bundle of the manifold.

Let $N$ be a smooth manifold.
A real-valued function $\lambda$ is said to be smooth in a neighborhood of $p$,
if the domain of $\lambda$ includes an open set $U$ of $N$ containing $p$
and the restriction of $\lambda$ to $U$ is a smooth function.
The set of all smooth functions in a neighborhood of $p$ is denoted
$\calC^\infty (p)$.
Consider $\lambda\in \calC^\infty (p)$, $\gamma \in \calC^\infty (p)$, and $a\in\RR$, $b\in\RR$.
Define the functions $a\lambda+b\gamma$ and $\lambda\gamma$ as
\[
(a\lambda+b\gamma)(q)=a\lambda(q)+b\gamma(q),
\]
\[
(\lambda\gamma)(q)=\lambda(q) \gamma(q),
\]
for all $q$ in the neighborhood of $p$.
It is obvious that $a\lambda+b\gamma\in \calC^\infty (p)$ and  $\lambda \gamma\in \calC^\infty (p)$.
So $\calC^\infty (p)$ forms a vector space over the field $\RR$.

A tangent vector $v$ at $p$ is a map $v : \calC^\infty (p) \rightarrow \RR$ with
\[
v(a\lambda+b\gamma)=a v(\lambda) +bv(\gamma),
\]
\[
v(\lambda\gamma)=\gamma(p)v(\lambda)+\lambda(p)v(\gamma),
\]
for all $\lambda,\gamma\in\calC^\infty(p)$ and $a,b\in\RR$.

Let $N$ be a smooth manifold. The tangent space to $N$ at $p$, denoted by $T_p N$, is the set of all tangent vectors at $p$.
The set $T_pN$ forms a vector space over the field $\RR$ under the normal rules of scalar multiplication and addition.

Let $N$ be smooth manifold of dimension $n$. Let $p$ be any point
of $N$, and $(U,\phi)$ be a coordinate chart around $p$. In this coordinate, the tangent
vectors $(\frac{\partial}{\partial \phi_1})_p$, $\cdots$, $(\frac{\partial}{\partial \phi_n})_p$ form a basis of $T_pN$,
which is called the natural basis of $T_pN$ induced by the coordinate chart $(U,\phi)$.
Let $v$ be a tangent vector at $p$, we have
\[
v=\sum_{i=1}^n v_i \Bigl( \frac{\partial}{\partial \phi_i} \Bigl)_p,
\]
where $v_1, \cdots, v_n$ are real numbers.

Let $(U,\phi)$ and $(V,\varphi)$ be coordinate charts around $p$. If $v$ is a tangent
vector, then
\[
v=\sum_{i=1}^n v_i\Bigl(\frac{\partial}{\partial \phi_i}\Bigl)_p
=\sum_{i=1}^n w_i\Bigl(\frac{\partial}{\partial \varphi_i}\Bigl)_p,
\]
where
\[
\pmatrix{v_1\cr v_2\cr \vdots\cr v_n}=
\mm{
\displaystyle\frac{\partial x_1}{\partial y_1} & \displaystyle\frac{\partial x_1}{\partial y_2} & \displaystyle\cdots & \displaystyle\frac{\partial x_1}{\partial y_n}\cr
\displaystyle\frac{\partial x_2}{\partial y_1} & \displaystyle\frac{\partial x_2}{\partial y_2} & \displaystyle\cdots & \displaystyle\frac{\partial x_1}{\partial y_n}\cr
\vdots & \vdots & \ddots &\vdots\cr
\displaystyle\frac{\partial x_n}{\partial y_1} & \displaystyle\frac{\partial x_n}{\partial y_2}  & \cdots & \displaystyle\frac{\partial x_n}{\partial y_n}
}
\pmatrix{w_1\cr  w_2\cr \vdots\cr w_n},
\]
and $x=x(y)$ represents the coordinate transformation $\phi\circ\varphi^{-1}$.

Let $N$ and $M$ be smooth manifolds. Let $F:N\rightarrow M$ be a smooth
mapping. The differential of $F$ at $p\in N$ is the map
\[
F_\star : T_p N\rightarrow T_p \rightarrow T_{F(p)} M
\]
defined as
\[
(F_\star(v))(\lambda)=v(\lambda\circ F),
\]
where $v\in T_p N$ and $\lambda\in\calC^\infty(F(p))$.

Let $(U,\phi)$ be a coordinate chart around $p$, $(V,\varphi)$ a coordinate
chart around $q=F(p)$.
The natural basis of $T_pN$ and $T_qM$ are $\Bigl\{\bigl(\frac{\partial}{\partial \phi_1}\bigl)_p,\bigl(\frac{\partial}{\partial \phi_2}\bigl)_p,\cdots,
\bigl(\frac{\partial}{\partial \phi_n}\bigl)_p\Bigl\}$
and
$\Bigl\{\bigl(\frac{\partial}{\partial \varphi_1}\bigl)_q,\bigl(\frac{\partial}{\partial \varphi_2}\bigl)_q,$ $\cdots,
\bigl(\frac{\partial}{\partial \varphi_n}\bigl)_q\Bigl\}$, respectively.
Denote the mapping $\varphi\circ F\circ \phi^{-1}$ as
\[
F(x)=F(x_1,x_2,\cdots,x_n)=\pmatrix{F_1(x_1,x_2,\cdots,x_n)\cr F_2(x_1,x_2,\cdots,x_n)\cr \vdots\cr F_m(x_1,x_2,\cdots,x_n)}.
\]
Suppose $v\in T_pN$ and $w=F_\star(v)\in T_{F(p)}M$ are expressed as
\[
v=\sum_{i=1}^n v_i \Bigl(\frac{\partial}{\partial \phi_i}\Bigl)_p, \quad
w=\sum_{i=1}^m w_i \Bigl(\frac{\partial}{\partial \varphi_i}\Bigl)_q,
\]
then
\[
\pmatrix{w_1\cr w_2\cr \vdots\cr w_m}=
\mm{
\displaystyle\frac{\partial F_1}{\partial x_1} & \displaystyle\frac{\partial F_1}{\partial x_2} & \cdots & \displaystyle\frac{\partial F_1}{\partial x_n}\cr
\displaystyle\frac{\partial F_2}{\partial x_1} & \displaystyle\frac{\partial F_2}{\partial x_2} & \cdots & \displaystyle\frac{\partial F_1}{\partial x_n}\cr
\vdots&\vdots &\ddots &\vdots\cr
\displaystyle\frac{\partial F_m}{\partial x_1} & \displaystyle\frac{\partial F_m}{\partial x_2} & \cdots & \displaystyle\frac{\partial F_m}{\partial x_n}
}
\pmatrix{v_1\cr v_2\cr \vdots\cr v_n}.
\]

\section{Vector Fields}

Consider a smooth manifold $N$ of dimension $n$. A vector field $f$
on $N$ is a mapping assigning to each point $p\in N$ a tangent vector $f(p)$ in $T_pN$.
A vector field $f$ is smooth if for each $p\in N$ there exists a coordinate chart $(U,\phi)$ about
$p$ and $n$ real-valued smooth function $f_1$, $f_2, \cdots$, $f_n$ defined on $U$ such that
for all $q\in U$
\[
f(q)=\sum_{i=1}^n f_i(q) \Bigl(\frac{\partial}{\partial \phi_i}\Bigl)_q.
\]

In local coordinates, $f_i$ can be expressed as
\[
\check{f}_k=f_i\circ \phi^{-1}.
\]
If $p$ is a point of coordinates $(x_1,x_2,\cdots,x_n)$ in the chart $(U,\phi)$, $f(p)$ is a
tangent vector of coefficients $(\check{f}_1(x_1,x_2,\cdots,x_n), \check{f}_2(x_1,x_2, \cdots,x_n), \cdots, \newline \check{f}_n(x_1,x_2,\cdots,x_n))$
in the basis
$\{(\frac{\partial}{\partial \phi_1})_p, (\frac{\partial}{\partial \phi_2})_p, \cdots, (\frac{\partial}{\partial \phi_n})_p\}$ of
$T_pN$.
Usually, $f_i$ is used to replace $f_i\circ \phi^{-1}$, therefore, $f$ in the local coordinates
is given by $f=\col(f_1,f_2,\cdots,f_n)$.

A smooth curve $\sigma: (t_1,t_2)\rightarrow N$ is an integral curve of $f$ if
\[
\sigma_\star\Bigl(\frac{d}{dt}\Bigl)_t=f(\sigma(t))
\]
for all $t\in(t_1,t_2)$.
By
\[
f(\sigma(t))=\sum_{i=1}^n f_i(\sigma_1(t),\sigma_2(t),\cdots,\sigma_n(tt)) \Bigl(\frac{\partial}{\partial \phi_i}\Bigl)_{\sigma(t)}
\]
\[
\sigma_\star \Bigl( \frac{d}{dt}\Bigl)_t=\sum_{i=1}^n \frac{d\sigma_i}{dt}\Bigl(\frac{\partial}{\partial \phi_i}\Bigl)_{\sigma(t)}.
\]
One obtains
\[
\frac{d\sigma_i}{dt}=f_i(\sigma_1(t),\sigma_2(t),\cdots,\sigma_n(t)).
\]

Let $f$ be a smooth vector field on $N$ and $\lambda$ a smooth real valued function
on $N$. The derivative of $\lambda$ along $f$ is a function $N\rightarrow \RR$, defined as
\[
(L_f\lambda)(p)=(f(p))(\lambda).
\]
In the local coordinates,
\[
(L_f)(x_1,x_2,\cdots,x_n)=\Bigl(\frac{\partial \lambda}{\partial x_1}\;\frac{\partial \lambda}{\partial x_2}\; \cdots \; \frac{\partial\lambda}{\partial x_n}\Bigl)
\pmatrix{f_1\cr f_2\cr \vdots\cr f_n}.
\]

The set of smooth vector fields on a manifold $N$, denoted by $V(N)$, is a vector space over $\RR$.
The vector space $V(N)$ is a Lie algebra if a binary operation $V\times V \rightarrow V$, called a product
and denoted by $[\cdot,\cdot]$, is defined such that
\ben
\item[(i)]
$
[v,w]=-[w,v];
$
\item[(ii)]
$
[\alpha_1v_1+\alpha_2v_2,w]=\alpha_1[v_1,w]+\alpha_2[v_2,w];
$
\item[(iii)]
$
[v,[w,z]]+[w,[z,v]]+[z,[v,w]]=0.
$
\een

If the product $[\cdot,\cdot]$ is defined as
\[
([f,g](p))(\lambda)=(L_fL_g\lambda)(p)-(L_gL_f\lambda)(p),
\]
the set $V(N)$ with the product forms a Lie algebra.

The product $[f,g]$ in local coordinates is given by
\[
\pmatrix{
\frac{\partial g_1}{\partial x_1} & \frac{\partial g_1}{\partial x_2} & \cdots & \frac{\partial g_1}{\partial x_n}\cr
\frac{\partial g_2}{\partial x_1} & \frac{\partial g_2}{\partial x_2} & \cdots & \frac{\partial g_2}{\partial x_n}\cr
\vdots &\vdots &\ddots &\vdots\cr
\frac{\partial g_n}{\partial x_1} & \frac{\partial g_n}{\partial x_2} & \cdots & \frac{\partial g_n}{\partial x_n} }
\pmatrix{f_1\cr f_1\cr \vdots \cr f_n}
-
\pmatrix{
\frac{\partial f_1}{\partial x_1} & \frac{\partial f_1}{\partial x_2} & \cdots & \frac{\partial f_1}{\partial x_n}\cr
\frac{\partial f_2}{\partial x_1} & \frac{\partial f_2}{\partial x_2} & \cdots & \frac{\partial f_2}{\partial x_n}\cr
\vdots &\vdots &\ddots &\vdots\cr
\frac{\partial f_n}{\partial x_1} & \frac{\partial f_n}{\partial x_2} & \cdots & \frac{\partial f_n}{\partial x_n} }
\pmatrix{g_1\cr g_2\cr \vdots \cr g_n}
\]
\[
=
\frac{\partial g}{\partial x} f- \frac{\partial f}{\partial x}g.
\]

The repeating product is possible. To avoid the notation of \newline $[f, [f, \cdots [f, g]\cdots ]]$ for recursive operation, define
\[
ad^k_f g(x)=[f, ad^{k-1}_f g] (x)
\]
for $k\geq 1$, where $ad^0_f g(x)=g(x)$.

\section{Distributions}

A distribution $D$ on a manifold $N$ is a map which assigns to each $p\in N$ a linear
subspace $D(p)$ of the tangent space $T_pN$.
If for each $p\in N$ there exists a neighborhood $U$ of $p$ and a set of smooth
vector fields $X_i$, $i\in I$, such that
\[
D(q)=\mbox{span} \{X_i(q), \; i\in I\}, \quad q\in U.
\]
The dimension of a distribution $D$ at $p\in N$ is the dimension of the subspace $D(p)$.
A distribution is constant dimensional if the dimension of $D(p)$ does not depend on the
point $p\in N$.

Let $D$ be a constant dimensional distribution of dimension $k$. Then around any $p\in M$
there exist $k$ independent vector fields $X_1,X_2,\cdots,X_k$ such that
\[
D(q)=\mbox{span}\{X_1(q),X_2(q),\cdots,X_k(q)\}.
\]
The vector fields $X_1, X_2, \cdots, X_k$ are called the local generators of $D$. Every vector
field $X\in D$ can be represented by
\[
X(q)=\sum_{i=1}^k \alpha_i(q) X_i(q)
\]
for some smooth function $\alpha_i$, $i=1,2,\cdots,k$.

A distribution $D$ is involutive if
\[
[X, Y]\in D
\]
for all $X\in D$ and $Y\in D$ .

A submanifold $P$ of $M$ is an integral manifold of a distribution $D$ on $M$ if
\[
T_qP=D(q), \forall q\in P.
\]

Let $X_1, X_2, \cdots, X_k$ be linearly independent vector fields with $[X_i,X_j]=0$, $1\leq i,j\leq k$.
Then there exist local coordinates such that
\[
X_i=\frac{\partial}{\partial x_i}, \quad 1\leq i\leq k.
\]
In other words, if $D$ is an involutive distribution of constant dimension $k$, then
there exist local coordinates $x_1,x_2,\cdots,x_n$ such that
\[
D=\mbox{span} \{\frac{\partial}{\partial x_1},\frac{\partial}{\partial x_2},\cdots,\frac{\partial}{\partial x_k}\}.
\]

%  \section{Nonlinear systems affine in control}
%
%  In the thesis, we study multiple input multiple output nonlinear systems affine in control,
%  \be
%  \label{nonsys0}
%  \left\{\begin{array}{rcl}\dot x &=& f(x)+\displaystyle\sum_{i=1}^m g_i(x)u_i,\cr
%                             y_i &=& h_i(x), \qquad 1\leq i\leq p. \end{array}\right.
%  \ee
%  The state
%  \[
%  x=\mm{x_1 \cr x_2 \cr \vdots \cr  x_n}
%  \]
%  is belong to an open set $U$ of $\RR^n$. There are $m$ inputs and $p$ outputs.
%  The mappings $f$, $g_1, g_2, \cdots, g_m$ are $n$ dimensional real valued mapping defined on $U$.
%  The mappings $h_1, h_2, \cdots, h_p$ are real valued mapping also defined on $U$.
%
%
%  \be
%  \label{nonsys}
%  \left\{\begin{array}{rcl}\dot x &=& f(x)+g(x)u,\cr
%                             y &=& h(x),\end{array}\right.
%  \ee
%  where $x\in\RR^n$, $u\in\RR^m$ and $y\in\RR^p$ are the state, input and output,
%  respectively.

\chapter{Normal Forms of Nonlinear Systems Affine in Control}

%\begin{abstract}
The nonlinear extensions of both finite and infinite zero structures
of linear systems have been well understood for single input single output
systems and have found many applications in nonlinear control
theory. The extensions of these notions to multiple input multiple
output systems have proven to be highly sophisticated.
Existing extensions either were
made under restrictive assumptions that not even square invertible
linear systems can satisfy or do not represent the nonlinear extensions
of the related linear system notions.
In this chapter, we propose constructive algorithms for decomposing a
nonlinear system that is affine in control. These algorithms require
modest assumptions on the system and apply to
general multiple input multiple output systems that do not necessarily
have the same number of inputs and outputs. They lead
to various normal form representations and
reveal the structure at infinity, the
zero dynamics and the invertibility properties, all of which
represent nonlinear extensions of relevant linear system structural
properties of the system they represent.

\section{Introduction}

The nonlinear analogues of linear system structural properties,
such as relative degrees (or infinite zero structure), zero dynamics
(or finite zero structure) and invertibility properties, have played
critical roles in recent literature on the analysis and control
design for nonlinear systems (see, {\em e.g.},
\cite{
byaf84,chgn03,conc99,diar07,foga94,hiio79m,isnc95,isnd81,jiau04,kane06,kaof05,khns02,liberzon2004output,manc96,nind90,pees96,
sett02,siam81,tegs94} and
the references therein for a sample of this literature).
The normal forms that are associated with these structural
properties, along with the basic tools like those reported in
\cite{esof92,isnc99,krna95,sags90,tetf95},
have enabled many major breakthroughs in nonlinear control theory.

Consider a multiple input multiple output (MIMO) nonlinear system affine in control
%with $m$ inputs $u=\col\{u_1,u_2,\cdots,u_m\}$,
%$p$ outputs $y=\col\{y_1,y_2,\cdots,y_p\}$,
%and the state $x=\col\{x_1,x_2,\cdots,x_n\}$ being defined in an open set $U$ of $\RR^n$,
\be
\label{nonsys_nfn}
\left\{\begin{array}{rcl}\dot x &=& f(x)+g(x)u,\cr
                           y &=& h(x),\end{array}\right.
\ee
where $x\in\RR^n$, $u\in\RR^m$ and $y\in\RR^p$ are the state, input and output,
respectively.
%In the equation, we have suppressed their dependency on time $t$ for
%notational brevity.
Let the mappings $f$, $g$ and $h$ be smooth in an open set $U\subset\RR^n$ containing
the origin $x=0$, with $f(0)=0$ and $h(0)=0$.

%and $f( 0 )=0$ and $h( 0 )=0$.
%In a neighborhood of a
%equilibrium point $ x=0 $,
%In the chapter, the mappings $f$, $g$ and $h$ are smooth
%with $f( 0 )=0$ and $g( 0 )=0$.
%as in \cite{isnc95,isnc99,nind90}.
%Note that they can also be analytic as in \cite{nind90}, or meromorphic as in \cite{conc07}.

A single input single output system,
{\em i.e.},
$m=p=1$ in (\ref{nonsys_nfn}), has a
relative degree $r$ at $ x=0 $ if
\be\label{siso1}
L_gL_f^k h(x)=0, \quad k<r-1,
\ee
in a neighborhood of $x=0$, and
\be\label{siso2}
L_gL_f^{r-1}h( 0 )\neq 0.
\ee
If system (\ref{nonsys_nfn}) has a relative degree $r$, then on
an appropriate set of coordinates in a neighborhood of  $ x=0$, it
takes the following normal form (see, {\em e.g.}, \cite{khns02}),
\be\label{basic}
\left\{\begin{array}{rcl}
\dot \eta &=& f_0(\eta,\xi),\cr
\dot\xi_i &=& \xi_{i+1},\quad i=1,2,\cdots,r-1,\cr
\dot\xi_{r} &=& a_1(\eta,\xi)
      +b_1(\eta,\xi)u,\cr
 y &=& \xi_1,
\end{array}\right.
\ee
where $\xi=\col\{\xi_1,\xi_2,\cdots,\xi_{r}\}$,
$b_1(0,0)\neq 0$, and
$
\dot\eta =f_0(\eta,0)
$
is the zero dynamics.
With a state feedback, this normal form reduces to
the zero dynamics cascaded with a {\em clean} chain of integrators linking
the input to the output.
Here by {\em clean} we mean that no other signal enters the middle of the chain.

Such a nice feature is extended to nonlinear systems with more
than one input output pairs. That is, a special class of square invertible nonlinear
system with $m=p>1$ can be transformed into the zero dynamics cascaded with
$m$ {\em clean} chains of integrators.
To do this, the notion of vector relative degree was introduced in \cite{saac89,byas91}.
System (\ref{nonsys_nfn})
with $m=p>1$
has a vector relative degree
$\{r_1,r_2,\cdots,r_m\}$ at $ x=0 $ if
\be\label{vrd1}
L_{g_j}L_f^k h_i(x)=0,\quad 0\leq k<r_i-1, \; 1\leq i,j\leq m,
\ee
in a neighborhood of $ x=0 $, and
\be\label{vrd2}
\det \{L_{g_j}L_f^{r_i-1} h_i( 0 )\}_{m\times m} \neq 0.
\ee
%where $A(x)=\{L_{g_j}L_f^{r_i-1} h_i(x)\}_{m\times m}$.
%\[
%A(x)=\pmatrix{
%L_{g_1}L_f^{r_1-1}h_1(x)&\cdots& L_{g_m}L_f^{r_1-1}h_1(x)\cr
%L_{g_1}L_f^{r_2-1}h_2(x)&\cdots& L_{g_m}L_f^{r_2-1}h_2(x)\cr
%\vdots&\ddots&\vdots\cr
%L_{g_1}L_f^{r_m-1}h_m(x)&\cdots& L_{g_m}L_f^{r_m-1}h_m(x)\cr
%}.
%\]
If system (\ref{nonsys_nfn}) has a vector relative degree $\{r_1,r_2,\cdots,r_m\}$ at $ x=0 $,
then with an appropriate change of coordinates, it can be described by
\be\label{vrd3}
\left\{\begin{array}{rcl}
\dot \eta &=& f_0(x)+g_0(x)u,\cr
\dot\xi_{i,j} &=& \xi_{i,j+1},\qquad j=1,2,\cdots,r_i-1,\cr
\dot\xi_{i,r_i} &=& a_i(x)+b_i(x)u,\cr
y_i &=& \xi_{i,1}, \qquad i=1,2,\cdots,m,
\end{array}\right.
\ee
which contains $m$ {\em clean} chains of integrators.
Moreover, if the distribution spanned by the column vectors of $g(x)$ is
involutive in a neighborhood of $ x=0 $,
a set of local coordinates can be selected such that $g_0(x)=0$.
The clean chains of integrators are called a prime form in \cite{mosi73} for linear systems,
and the necessary and sufficient geometric
conditions for the existence of prime forms for nonlinear systems is developed in \cite{marino1994equivalence}.
There is a large body of nonlinear systems and control literature based on
the form (\ref{vrd3}) (see e.g.,
%\cite{nind90,isnc95,khns02} and references therein).
\cite{atas99,huot03,husi97,jido98,lios02,oria03},
for a small sample).

The conditions for the existence of a vector relative degree, (\ref{vrd1}) and (\ref{vrd2}),
though similar to (\ref{siso1}) and (\ref{siso2})
in form, are not easy to be satisfied.
Simple change of coordinates in the output space could alter the property (\ref{vrd1}).
That is, (\ref{vrd1}) is satisfied only under certain output coordinates.
Consider the linear system $(A,B,C)$ from \cite{scgn99},
\be
%A=\left[ \begin {array}{rrrrr}
A=\mm{
1&1&-2&0&0\cr
0&5&-4&1&2\cr
0&1&0&0&1\cr
-2&0&-1&0&0\cr
0&0&0&-1&0
%\end {array} \right],\quad
},\quad
B=\mm{
0&0\cr
1&1\cr
0&0\cr
-1&1\cr
0&0
},
\]
\[
C=\mm{
0&1&-2&0&0\cr
0&1&-2&0&1
}.\label{exam_sch}
\ee
As shown in \cite{scgn99}, the system does not possess a vector relative degree.
If we apply an output transformation
\[
T_\rmo=\mm{ 1 &0\cr -1 &1},
\]
it can be verified that $(A,B,T_\rmo C)$ has a vector relative degree $\{1,2\}$.
In other words, the system (\ref{exam_sch}) meets the vector relative degree conditions
only under appropriate coordinates of the output space.

In general,
the vector relative degree is a rather restrictive structural property that
not even all square invertible linear systems, with the freedom
of choosing coordinates for the state, output and input spaces and state feedback,
could possess.
A square invertible linear system with $m=p>1$ in general can only
be transformed into the zero dynamics cascaded with $m$ chains of
integrators, with all but one chains containing output injection terms
(see \cite{saas87}). That is, there are interconnections between these chains.
For example, consider a linear system $(A,B,C)$ with
\be\label{exam1}
A=\mm{
     0&     0&     0&     0\cr
     \alpha&     0&     1&     0\cr
     0&     0&     0&     1\cr
     0&     0&     0&     0},\quad
B=\mm{
     1&     0\cr
     0&     0\cr
     0&     0\cr
     0&     1},\quad
C=\mm{
     1&     0&     0&     0\cr
     0&     1&     0&     0},\quad \alpha\neq 0.
\ee
The system contains two chains of integrators of lengths 1 and 3.
The parameter $\alpha$ represents an output injection term, which in turn represents the
interconnections between the two chains.
Such an interconnection cannot be removed through coordinate transformations and state feedback,
and thus
system (\ref{exam1}) cannot be represented by two {\em clean} chains of integrators.
In other words, even with the freedom of choosing coordinates and static state feedback, (\ref{exam1})
does not have a vector relative degree.
To see this, suppose  that there exist nonsingular coordinate transformations
$T_\rms$, $T_\rmi$ and $T_\rmo$ such that
\[
\tilde{A}=T_\rms^{-1}A T_\rms=\mm{
     *&     *&     *&     *\cr
     0&     0&     1&     0\cr
     0&     0&     0&     1\cr
     *&     *&     *&     *},\quad
T_\rms^{-1} BT_\rmi=B,\quad
T_\rmo^{-1}C T_\rms=C,
\]
which indicates that the system can be decoupled into two {\em clean} chains of integrators.
Denote $T_\rms=\{t_{i,j}\}_{4\times 4}$.
By $T_\rms B=BT_\rmi$ and $C T_\rms=T_\rmo C$, we obtain
$t_{1,3}=t_{1,4}=t_{2,1}=t_{2,3}=t_{2,4}=t_{3,1}=t_{3,4}=0$.
The $(2,1)$ entry of $AT_\rms-T_\rms \tilde{A}$ is $\alpha t_{1,1}=0$.
So, $t_{1,1}=0$, consequently, $T_\rms$ is singular. This is a contradiction.

Similarly,
\be\label{exam2}
A=\mm{
     0&     0&     0&     0\cr
     0&     0&     1&     0\cr
     0&     0&     0&     1\cr
     0&     0&     0&     0},\quad
B=\mm{
     1&     0\cr
     \alpha&     0\cr
     0&     0\cr
     0&     1},\quad
C=\mm{
     1&     0&     0&     0\cr
     0&     1&     0&     0},\quad \alpha\neq 0.
\ee
The system in (\ref{exam2})
does not have a vector relative degree even with the freedom of choosing coordinates and static state feedback.
The parameter $\alpha$ here represents an input coupling term between the two chains.

A major generalization of the form (\ref{vrd3}) was made in
\cite{byls88,isnc99,isnc95,scgn99}, where
MIMO square invertible systems are considered.
In \cite{isnc95},
with Zero Dynamics Algorithm,
a sequence of nested submanifolds
$M_0\supset M_1 \supset \cdots \supset M_k \supset \cdots=Z^\star$
are defined, and
system (\ref{nonsys_nfn}) is transformed into the form,
\be\label{zerody1}
\left\{\begin{array}{rcl}
\dot \eta &=& f_0(x)+g_0(x)u,\cr
\dot\xi_{i,j} &=& \xi_{i,j+1}+
\displaystyle            \sum_{l=1}^{i-1}
            \delta_{i,j,l}(x)v_l
            +\sigma_{i,j}(x)u, \quad j=1,2,\cdots,n_i-1,\cr
\dot\xi_{i,n_i} &=& v_i,\cr
 y_i &=& \xi_{i,1}, \qquad i=1,2,\cdots,m,
\end{array}\right.
\ee
where
$\sigma_{1,j}(x)=0,\; j=1,2,\cdots,n_1-1$,
and $\sigma_{i,j}(x)=0,\; i>1, \;j=1,2,\cdots,n_i-1$, in $Z^\star$,
and the static state feedback is given by $v_i=a_i(x)+b_i(x)u$, $i=1,2,\cdots,m$, with
the matrix $\col\{b_1(x),b_2(x),\cdots,b_m(x)\}$ being smooth and nonsingular.
In the algorithm,
the rank of certain matrices are assumed to be constant
on these nested submanifolds.
With some stronger assumptions imposed in the algorithm \cite{isnc99,scgn99},
{\em i.e.}, the rank of certain matrices were assumed to be constant
for all $x\in U$ (not just in these submanifolds),
one can have all $\sigma_{i,j}(x)=0$.
Moreover, if certain vector fields commute,
%distribution spanned by column vectors of $g(x)$ is involutive,
one can select coordinates such that $g_0(x)=0$. Thus,
system (\ref{zerody1}) becomes
\be\label{zerody2}
\left\{\begin{array}{rcl}
\dot \eta &=& f_0(x),\cr
\dot\xi_{i,j} &=& \xi_{i,j+1}+
  \displaystyle          \sum_{l=1}^{i-1}
            \delta_{i,j,l}(x)v_l, \quad j=1,2,\cdots,n_i-1,\cr
\dot\xi_{i,n_i} &=& v_i,\cr
 y_i &=& \xi_{i,1}, \qquad i=1,2,\cdots,m.
\end{array}\right.
\ee
The applications of the form (\ref{zerody2}) in solving the problem of
asymptotic stabilization,
disturbance decoupling, tracking and regulation can be found in \cite{isnc99} and the references therein.

The infinite zeros of a linear system can be defined either through
the root locus theory or as the Smith-McMillan zeros of the transfer function at infinity
\cite{ross70,hatf76}.
They can also be characterized in state-space \cite{mosi73,wolm79}.
%,desi83}.
On the other hand, the structure at infinity was introduced
for a certain class of nonlinear systems in \cite{isnf83}, and was further developed
for smooth systems or analytic systems in \cite{niza85}
and for meromorphic systems in \cite{flan86,mond88,beri89}.

In \cite{niza85} and \cite{nind90} (Chapter 9),
formal zeros at infinity are defined in terms of a set of geometric conditions.
In particular,
for system (\ref{nonsys_nfn}) defined on a smooth manifold $M$,
a sequence of
the locally controlled invariant
distributions
\[
\calD_0\supset\calD_1\supset\cdots\supset\calD_n
\]
in $\ker dh(x)$ are defined, where
\[
\begin{array}{lll}
\calD_0&=&TM, \\
\calD_{i+1}&=&
\bigl\{X\in V(M) |\, [f,X]\in \calD_i+G,\, [g_j,X]\in \calD_i+G, \;\; j=1,2,\cdots,m\bigl\} \rule{0pt}{15pt} \\
& & \;\; \cap\ker dh, \qquad i=1,2,\cdots,n-1, \rule{0pt}{15pt} \\
\end{array}
\]
and where $G=\span\{g_1,g_2,\cdots,g_m\}$,
$V(M)$ denotes the set of smooth vector fields on a smooth manifold $M$,
and $TM$ denotes the tangent bundle of $M$.
Under the assumption of the distributions $D_i$ and $D_i\cap G$ on $M$ having
constant dimensions, the formal zeros at infinity can be defined.
Formal zeros at infinity plays an important role in
the input output decoupling problem by static state feedback, in which
after possible relabeling of inputs, the control $u_i$ does not
influence the output $y_j$, $j\neq i$.
But this structure information does not show in a normal form
in the references \cite{niza85,nind90}.

In \cite{conc99,beri89}, %\cite{mond88,beri89},   %conc99,conc07},
a linear-algebraic strategy is developed based
on the use of vector spaces over the field of meromorphic functions.
As a counterpart to the above differential-geometric approach, the algebraic approach
considers system (\ref{nonsys_nfn}) with $f$, $g$ and $h$ being meromorphic.
Except for some singular points,
the two approaches lead to the same results as in \cite{respondek1990right},
in particular, the same notions of rank and structure at infinity.
The structure at infinity is related to a chain of subspaces
$\calE_0\subset\calE_1\subset\cdots\subset\calE_n$, where
\[
\calE_0=\span_\calK\{dx\},
\]
\[
\calE_i=\span_\calK\{dx, d\dot{y},\cdots,dy^{(i)}\},\; i=1,2,\cdots,n,
\]
and where $\calK$ are meromorphic functions.
The structure at infinity is then determined by
\[
\sigma_k=\mbox{dim}_\calK \frac{\calE_k}{\calE_{k-1}}.
\]
With a generalized state space transformation, a regular generalized state feedback
and a universal, additive output injection,
system (\ref{nonsys_nfn}) can be transformed into a canonical form,
which contains time derivatives of inputs and
shows the structure at infinity explicitly.
%Note that the transformation cannot be physically
%implemented, since output injection cannot be carried out to a given dynamic structure through
%its input and output channels.

As pointed out in \cite{isnc95}, if all $\delta_{i,j,l}(x)=0$,
the set of integers $\{n_1,n_2,\cdots,n_m\}$ in (\ref{zerody2})
corresponds to
the vector relative degree, which in this case, represents
the infinite zero structure if the system is linear.
These integers however are not related to the infinite zero structure
of linear systems when $\delta_{i,j,l}(x)\neq 0$, and thus
cannot be defined as the nonlinear extension of and
expected to play a similar role as infinite zeros.
To see this, consider the following linear system $(A,B,C)$,
\be\label{counter3}
A=
\mm{
     0&     1&     0&     0&     0\cr
     0&     0&     0&     0&     0\cr
     0&     0&     0&     1&     0\cr
     0&     0&     0&     0&     1\cr
     0&     0&     0&     0&     0},\quad
B=\mm{
     0&     0\cr
     1&     0\cr
     \alpha&     0\cr
     0&     0\cr
     0&     1},\quad
     \]
     \[
C=\mm{
     1&     0&     0&     0&     0\cr
     0&     0&     1&     0&     0},\quad
\alpha\neq 0, \label{linexam1}
\ee
which is in the form of (\ref{zerody2})
with $n_1=2$ and $n_2=3$.
However, by using the toolkit \cite{chls04,lils04,liu2005linear}, we can
find state, input and output transformations $T_\rms$, $T_\rmi$ and $T_\rmo$ such that
\[
T_\rms^{-1}AT_\rms=\mm{
     1&     0&     0&    -\alpha&     0\cr
    1/\alpha&     0&     1&     0&     0\cr
     0&     0&     0&     1&     0\cr
     0&     0&     0&     0&     1\cr
     0&     0&     0&     0&     0},\quad
T_\rms^{-1}BT_\rmi=\mm{
     1&     0\cr
     0&     0\cr
     0&     0\cr
     0&     0\cr
     0&     1},
\]
\[
T_\rmo^{-1} C T_\rms=\mm{
     1&     0&     0&     0&     0\cr
     0&     1&     0&     0&     0},
\]
\[
T_\rms=\mm{
     0&     1&     0&     0&     0\cr
 1/\alpha&     0&     1&     0&     0\cr
     1&     0&     0&     0&     0\cr
     0&     0&     0&    -\alpha&     0\cr
     0&     0&     0&     0&    -\alpha},\quad
T_\rmo=\mm{
     0&     1\cr
     1&     0},\quad
T_\rmi=\left[\begin{array}{rr}
     1/\alpha&     0\cr
     0&    -\alpha\end{array}\right].
\]
Thus, according to \cite{mosi73,saas87},
the system is invertible with
two infinite zeros $\{1,4\}$.
Therefore, the integers $n_1$ and $n_2$ in the form (\ref{zerody2})
does not generalize the notion of infinite zero structure of linear systems.

Invertibility of linear systems was first studied in
\cite{brtr65,saio69,siio69}.  In these references, inversion algorithms
and invertibility criteria are given.
Invertibility of nonlinear control systems was considered in
\cite{hiio79m,siam81}, which generalized the structure algorithm for linear systems \cite{siio69}.
Refs. \cite{nijmeuer1988dynamic,respondek1990right,respondek1988local}
carry out a systematic study of invertibility
of general nonlinear systems that are not necessarily affine in control.
The authors gave a list of equivalent conditions for right and
left invertibility for linear systems, and examined when and how these conditions
can be generalized to nonlinear systems. Based on \cite{siam81}, ref.
\cite{liberzon2004output} explicitly
constructs the left inverse of an affine output-input stable system.

Invertibility of nonlinear systems can also be determined by using the
structure algorithm in \cite{beri89,conc99,mond88,respondek1990right}.
In particular,  with a generalized state space transformation, a regular
generalized state feedback
and a universal, additive output injection,
system (\ref{nonsys_nfn}) can be transformed into a canonical form,
which contains time derivatives of inputs and
shows the structure at infinity and invertibility structures explicitly.
%Note that the transformation cannot be physically
%implemented, since output injection cannot be carried out to a given
%dynamic structure through
%its input and output channels.

A key feature of the normal forms is
that they represent a system in several interconnected subsystems.
These subsystems, along with the interconnections that exist among them,
lead us to a deeper insight into how control would take
effect on the system, and thus to the construction of control laws that
meet our design specifications.
The structure of a linear system
,characterized by a matrix triple $(A,B,C)$,
has been studied in great depth.
In 1973, Morse \cite{mosi73} showed that, under a group
of state, input and output transformations,
state feedback and output injection,
any matrix triple $(A,B,C)$ is uniquely characterized
by three lists of positive integers and a list of monic polynomials.
By identifying state variables in the structure algorithm in \cite{siio69},
Sannuti and Saberi \cite{saas87} explicitly constructed
state, input and output transformations that transform a general
MIMO system, not necessarily
square, into a so-called special coordinate basis form,
which displays all structural properties of
the system, including the finite and infinite zero structures
and invertibility properties.

Motivated by the many efforts reported in the nonlinear control
literature and a complete understanding and numerous applications
of the structural decomposition of linear systems,
we make an attempt to study structural properties of affine nonlinear systems
beyond the case of square invertible systems.
For a general nonlinear system (\ref{nonsys_nfn}) in the absence
of the vector relative degree assumption, we develop an algorithm, which is referred to as the
infinite zero structure algorithm and, under certain constant rank assumptions over $U$, results in
diffeomorphic state, input and output transformations and state feedback laws
under which the system can be represented in normal forms. In the special
case when the system is square and invertible, our normal forms take a form
similar to those in \cite{isnc99,scgn99}, but with an additional
property that allows the normal forms to reveal the nonlinear extension of infinite
zeros of linear systems. In addition, our development enhances
the existing results in some other ways. First, fewer assumptions are required.
Second, the resulting normal forms explicitly show invertibility structures and
nonlinear extension of invariant zeros. Third, our development
applies to general MIMO nonlinear systems that are not necessarily square.

The infinite zero structure algorithm will also be adapted to develop normal
forms that reveal system structural properties when the output is restricted to
zero. The adapted algorithm will be referred to as the zero output structure
algorithm. The assumptions required will also be in the form of constant ranks,
but in a sequence of nested subsets, rather than the more stringent constant
ranks on $U$ as required by the infinite zero structure algorithm. Our results
on zero output normal forms inherit the features pertaining to the infinite
zero structure algorithm and thus enhance the existing results on the
zero output normal forms in similar ways as the normal forms resulting from
the infinite zero structure algorithm.

These normal
forms include the ones
identified in \cite{byls88,isnc95,isnc99,scgn99} for square invertible systems as special cases.
In particular,
Under the milder assumptions on nonlinear systems,
and by carefully selecting new coordinates, simpler normal forms can be derived.
These normal forms not only reveal the infinite zero structure and
zero dynamics of the system, but also provide explicit information
on the system invertibility properties.
So far, the structure at infinity is related only to input-output decoupleable nonlinear systems.
In the chapter, we try to extend the concept of structure at infinity to input-output coupling nonlinear %systems.
In doing so, we introduce
the notions of
infinite zero of nonlinear systems.
The systems are not necessarily square.
We also explore the structural properties of nonlinear systems along the trajectory in which the output is %identically zero.
and introduce the notions of zero-structure at infinity, zero-invertibility
of nonlinear systems at
an equilibrium point
$ x=0 $.

The remainder of this chapter is organized as follows.
The infinite zero structure algorithm and the resulting normal forms are presented in Section \ref{algo_inf}.
The zero output structure algorithm and the resulting normal forms are given in Section \ref{algo_zerooutput}.
Section \ref{example_nfn} contains a few examples that
illustrate the main results of the chapter.
A brief conclusion to the chapter is drawn in Section \ref{conclusion_nfn}.
For clarity in the presentation, all proofs are given in the appendices.

%Throughout the chapter,
%for matrices $A_1, A_2, \cdots,$ and  $A_n$,
%all with a same number of columns, denote
%$
%\col\{A_1,A_2 \cdots, A_n\}
%=[A_1^\rmT\; A_2^\rmT\; \cdots\; A_n^\rmT]^\rmT.
%%=\mm{A_1\cr A_2\cr \vdots\cr A_n}.
%$

%For a differentiable function $h: \RR^n \rightarrow \RR$,
%$
%h(x)=\frac{\partial h(x)}{\partial x}=
%\left(\frac{\partial h}{\partial x_1}\;\;
%   \frac{\partial h}{\partial x_2}\;\;
%\cdots\;\; \frac{\partial h}{\partial x_n}\right).
%$
%With a slight abuse of notation,

\section{Normal Forms and Structure Properties of Nonlinear Systems}
\label{algo_inf}

In this section, we will find diffeomorphic
state, input and output transformations and static state
feedback laws under which system (\ref{nonsys_nfn}) can be represented in normal forms and discuss
about the intrinsic structural properties these normal forms reveal.
Similarly to many existing results (see, {\em e.g.}, \cite{isnc99,scgn99}), we rely on constant
rank assumptions over $U$. However, as will become clear,
our development here enhances the existing results in several ways. First, weaker assumptions
are required. Second, normal forms with simpler structure are resulted in, based on which
nonlinear extension of infinite zeros can be defined. Third, the resulting normal forms
explicitly show invertibility structures and nonlinear extension of invariant zeros.
Finally, our normal form development applies to nonlinear systems that are not necessarily square.

In particular, we first separate from the overall system
dynamics the dynamics associated with the infinite zeros, and then
carry out some further decomposition of the zero dynamics and the remaining dynamics.

%\subsection{Dynamics Associated with the Infinite Zero Structure}

\subsection{The Infinite Zero Structure Algorithm}
\label{infinite_zero_structure_algorithm}

Both our algorithm and the algorithm in \cite{isnc99,scgn99} involve repetitive
differentiations of the output and, under certain constant rank assumptions,
identification of functions to serve as new state variables.
What distinguishes our algorithm is how we identify the new state variables.
In each step of our algorithm, we identify not only
$\Theta_k(x)$, from which new state variables will be selected,
but also $\Omega_k(x)$, which contains $\Omega_{k-1}(x)$ and part of $\Theta_{k-1}(x)$,
in such a way that $L_g \Omega_k(x)$ is of full row rank
and
\[
\rank (L_g \Omega_k(x))=\rank(L_g \col\{\Theta_0(x),\Theta_1(x),\cdots,\Theta_{k-1}(x)\}).
\]
More specifically, we first identify $\Omega_k(x)$, then define $\Theta_k(x)$ to depend only on $\Theta_{k-1}(x)$ and $\Omega_k(x)$,
rather than on $\Theta_i(x)$, $i=1,2,\cdots,k-1$. Such an approach will be helpful
in selecting state variables that render the more informative normal forms.

Moreover, by choosing the function $\Theta_k(x)$ in such a way, we will be able to
carry out the algorithm with fewer constant rank assumptions than the
algorithm in \cite{isnc99,scgn99}, and more importantly, allow the algorithm
to be applicable to square but non-invertible systems and non-square systems.

We also will device criteria for the above repetitive procedure to stop. The times the
derivatives are taken on each output variable and which stopping
criterion is met determine the structure at infinity and the invertibility properties,
respectively.

{\bf Initial Step.}
Let $\Theta_0(x)=h(x)$, $\Omega_0(x)=\emptyset$, $\rho_0=0$ and $k=1$.

{\bf Step $k$.} We start with
$\Theta_{k-1}(x): U\rightarrow\RR^{p-\rho_{k-1}}$,
$\Omega_{k-1}(x): U\rightarrow\RR^{\rho_{k-1}}$,
where the matrix $L_g\Omega_{k-1}(x)$ has full row rank $\rho_{k-1}$.
Suppose that the following assumption holds.
%\[
%\langle \mbox{d} \Theta_{k-1}(x), f(x)+g(x) u\rangle
%=L_f \Theta_{k-1}(x) + L_g \Theta_{k-1}(x) u.
%\]

{\em Assumption $\calA_k$:
The matrix
$\mm{L_g\Omega_{k-1}(x)\cr L_g\Theta_{k-1}(x)}$
has constant rank $\rho_k$ for $x\in U$, and there exists
an $R_k \in \RR^{(\rho_k-\rho_{k-1})\times(p-\rho_{k-1})}$
such that the matrix
$
\mm{L_g\Omega_{k-1}(x)\cr L_g R_k \Theta_{k-1}(x)}
$
is of full row rank $\rho_k$ for $x\in U$.
}

Let
$S_k\in\bbR^{(p-\rho_k)\times(p-\rho_{k-1})}$ be
such that
\be\label{RkSk}
\det\left(\mm{R_k\cr S_k}\right)\neq 0.
\ee
Denote
\be\label{Omega_k}
\Omega_k(x)
%=\mm{\Omega_{k-1}(x)\cr R_k\Theta_{k-1}(x)}
=\col\{R_1\Theta_0(x), R_2\Theta_1(x), \cdots, R_k\Theta_{k-1}(x)\}.
\ee
The matrix $L_g \Omega_k(x)$ has full row rank $\rho_k$ and
\[
\rank\left(\mm{L_g\Omega_k(x)\cr L_g S_k\Theta_{k-1}(x)}\right)=\rho_k.
\]
Thus, there exist unique smooth functions
\[
P_{k,l}(x): U\rightarrow \RR^{(p-\rho_k)\times(\rho_l-\rho_{l-1})},\;l=1,2,\cdots,k,
\]
such that
\be\label{Pkj}
L_g  S_k\Theta_{k-1}(x)
-\sum_{l=1}^k P_{k,l}(x)L_g R_l\Theta_{l-1}(x)=0.
\ee
Define
\begin{eqnarray}\label{Thetak}
\Theta_k(x)
= L_f    S_k\Theta_{k-1}(x)
-\sum_{l=1}^k P_{k,l}(x)L_f R_l\Theta_{l-1}(x).
\end{eqnarray}
If $\;k+\sum_{j=1}^k j\,(\rho_j-\rho_{j-1})<n$ and $\rho_{k}<\min\{p,m\}$,
then increase $k$ by $1$ and repeat the above step. Otherwise, go to Final Step.

{\bf Final Step.}
Let $k^\star=k$, we have
\be\label{endalgo}
{k^\star}+\sum_{j=1}^{k^\star} j\,(\rho_j-\rho_{j-1})=n \quad \mbox{or} \quad \rho_{k^\star}=\min\{p,m\}.
\ee
Let $m_\rmd=\rho_{k^\star}$ and
\[
n_\rmd=\sum_{j=1}^{k^\star} j\,(\rho_j-\rho_{j-1}).
\]
Denote the set $\rho=\{\rho_1,\rho_2,\cdots,\rho_{k^\star}\}$.
Define a set of integers $0<q_1\leq q_2\leq \dots\leq q_{m_\rmd}$ as
\[
q=\{q_1, q_2, \dots, q_{m_\rmd}\}=
\{\;\overbrace{1,\cdots,1}^{\rho_1-\rho_0},\;\;
\overbrace{2,\cdots,2}^{\rho_2-\rho_1},\;\;
\cdots,\;\;
\overbrace{k^\star,\cdots,^{^{}} k^\star}^{\rho_{k^\star}-\rho_{k^\star-1}}\; \}.
\]

{\bf End.}

\begin{defi}\em
System (\ref{nonsys_nfn}) is said to be regular, if Assumption $\calA_k$, $k=1,2,\cdots,k^\star$,
are satisfied.
\end{defi}

\subsection{Normal Forms}

We will base on the infinite zero structure algorithm to derive normal forms of system (\ref{nonsys_nfn}).
Denote
\be
v_k=L_f R_k\Theta_{k-1}(x)+L_g  R_k\Theta_{k-1}(x)u.\label{vk}
\ee
By (\ref{Pkj}) and (\ref{Thetak}),
\be
\label{Theii1}
\Theta_j(x)
=\frac{d}{dt} S_j\Theta_{j-1}(x)-\sum_{l=1}^j P_{j,l}(x)v_l,
\quad j=1,2,\cdots, k^\star.
\ee
For the notational brevity, denote $S_i S_{i-1}\cdots S_{j+1} S_j$ as $S_{i\leftrightarrow j}$,
with $ S_{i\leftrightarrow j}=1$ for $j>i$.
We first define the new states representing the dynamics of $i(\rho_i-\rho_{i-1})$
integrators, which connect the input $L_g R_i\Theta_{i-1}(x)u$ to the
output $ R_i S_{i-1\leftrightarrow 1}y$,
\begin{eqnarray*}
\zeta_{i,j}&=&R_i S_{i-1\leftrightarrow j}\Theta_{j-1}(x),\quad j=1,2,\cdots,i,\\
\zeta_i&=&\col\{\zeta_{i,1},\zeta_{i,2},\cdots,\zeta_{i,i}\},\\
&=&\col\{
R_i S_{i-1\leftrightarrow 1}\Theta_0(x),
R_i S_{i-1\leftrightarrow 2}\Theta_1(x),
\cdots,R_i\Theta_{i-1}(x)\},\quad i=1,2,\cdots,k^\star.
\end{eqnarray*}
Note that $\zeta_{i,j} : U\rightarrow\RR^{\rho_i-\rho_{i-1}}$.
In view of (\ref{vk}) and (\ref{Theii1}), we have,
\begin{eqnarray}
\label{dzki}
\dot{\zeta}_{i,j}&=&
\zeta_{i,j+1}+R_i S_{i-1\leftrightarrow j+1}
 \sum_{l=1}^j P_{j,l}(x)v_l,\quad j=1,2,\cdots,i-1,\\
\dot{\zeta}_{i,i}&=&v_i, \quad i=1,2,\cdots,k^\star. \nonumber\label{zkk}
\end{eqnarray}
Let
\begin{eqnarray}
\bar{\Phi}_\rmd(x) &=&\col\{\zeta_1,\zeta_2,\cdots,\zeta_{k^\star}\},\label{Phi1}\\
{\Gamma}_{\rmi\rmd}(x)&=&\col\{L_g R_1\Theta_{0}(x),
L_g R_2\Theta_{1}(x),\cdots,
L_g R_{k^\star}\Theta_{k^\star-1}(x)\}  \nonumber \\ &=&L_g\Omega_{k^\star}(x),\\ \label{Gammaid}
{\Gamma}_{\rmo\rmd}&=&\col\{R_1, R_2 S_{1},\cdots,
         R_{k^\star} S_{k^\star-1\leftrightarrow 1}\}.\label{Gammaod}
\end{eqnarray}
It is obvious that
$\bar{\Phi}_\rmd(x): U\rightarrow\RR^{n_\rmd}$,
${\Gamma}_{\rmi\rmd}(x) : U\rightarrow\RR^{m_\rmd\times m} $
and ${\Gamma}_{\rmo\rmd}\in \RR^{m_\rmd\times p}$.
To construct a new set of coordinates, we need the following assumption.

{\em Assumption $\calB$:
The matrix $d\bar{\Phi}_\rmd(x)$ is of full row rank for $x\in U$.}

Note that Assumption $\calB$ is automatically satisfied if
$P_{k,l}(x),\;l=1,2,\cdots,k,\; \newline k=1,2,\cdots,k^\star$, in the infinite zero structure algorithm
are independent of $x$.

\begin{lemma}
\label{independent}
Suppose that system (\ref{nonsys_nfn}) is regular, and that
$P_{k,l}(x),\;l=1,2,\cdots,k,\; k=1,2,\cdots,k^\star$, in the infinite zero structure algorithm are constant matrices.
Then, $d\bar{\Phi}_\rmd(x)$ is of full row rank for $x\in U$.
\end{lemma}
{\bf Proof:} See Appendices~\ref{p_independent}.
\hfill $\square$

By the infinite zero structure algorithm, we know that
$\Gamma_{\rmi\rmd}(x)=L_g\Omega_{k^\star}(x)$ is of full row rank.
Note that $R_i S_{i-1 \leftrightarrow 1}$, $i=1,2,\cdots,{k^\star}$, are the coefficients
in $\zeta_{i,1}=R_i S_{i-1\leftrightarrow 1}\Theta_0(x)$.
Under Assumption $\calB$, $\Gamma_{\rmo\rmd}$ is of full row rank.
In what follows, we augment the state variables $\zeta_{i,j}$'s with
$n-n_\rmd$ additional state variables to form a full set of state variables
for the system.
%We need $n-n_\rmd$ additional state variables.
Similarly, we also need
to augment the input variables $v_i$'s and the output variables
$\zeta_{i,1}$'s with
$m-m_\rmd$ additional
input variables and $p-m_\rmd$ output variables
to form a full input
vector and output vector, respectively.

Note that $\zeta_{i,j}$ contains $\rho_i-\rho_{i-1}$ states, and thus
$\zeta_{i,j}$, $j=1,2,\cdots,i$,
define $\rho_i-\rho_{i-1}$ chains containing a total of $i(\rho_i-\rho_{i-1})$ integrators.
If $\rho_i-\rho_{i-1}>1$, we introduce the permutation matrix
$\Xi(\rho_i-\rho_{i-1},i)$ to reorder the states such that
each chain contains $(\rho_i-\rho_{i-1})$ integrators
and corresponds to only one input and one output,
where $\Xi(s,t)\in\RR^{st\times st}$ with
\[
\Xi(s,t)=\bigl[e_1\;\; e_{t+1}\;\; \cdots\; e_{(s-1)t+1}|\;
e_2\;\; e_{t+2}\; \cdots\;
e_{(s-1)t+2}\; |\; 
\]
\[
\hspace{60mm}\cdots\;
|\; e_s e_{t+s}\;\; \cdots\; e_{st}\bigl]^\rmT
\]
and $e_l$ being the $l$th column of the identity matrix $I_{st}$.
Define
\[
\xi={\Phi}_{\rmd}(x)
=\Upsilon\bar{\Phi}_{\rmd}(x)\hspace{90mm}
\]
\be
=\col\{
\Xi(\rho_1-\rho_0,1)\zeta_{1},
\Xi(\rho_2-\rho_1,2)\zeta_{2},\cdots,
\Xi(\rho_{k^\star}-\rho_{{k^\star}-1},{k^\star})\zeta_{k^\star}
\},\label{xi_zeta}
\ee
\[
\xi=\col\{\xi_1,\xi_2,\cdots,\xi_{m_\rmd}\},\quad
\xi_i=\col\{\xi_{i,1},\xi_{i,2},\cdots,\xi_{i,q_i}\},
\]
where
\[
\Upsilon=
\blkdiag\{\Xi(\rho_1-\rho_0,1),\Xi(\rho_2-\rho_1,2),\cdots,\Xi(\rho_{k^\star}-\rho_{{k^\star}-1},{k^\star})\}.
\]
Note that if $\rho_i-\rho_{i-1}\leq 1$ for $0\leq i\leq {k^\star}$, then
$\Upsilon=I$.
Define a new set of coordinates,
\[
\pmatrix{\eta\cr \xi}
=\pmatrix{\Phi_\rme(x)\cr \Phi_{\rmd}(x)}=\Phi(x), \quad
\pmatrix{u_\rme\cr u_{\rmd}}=\Gamma_\rmi(x)u
=\mm{\Gamma_{\rmi\rme}(x)\cr \Gamma_{\rmi\rmd}(x)}u,\quad
\]
\be
\label{trans1a}
\label{trans1b}
\label{trans1c}
\pmatrix{y_\rme\cr y_{\rmd}}=\Gamma_\rmo y
=\mm{\Gamma_{\rmo\rme}\cr \Gamma_{\rmo\rmd}}y, \hspace{50mm}
\ee
where
${\Phi}_\rme(x) : U\rightarrow\RR^{n-n_\rmd}$ is smooth
and such that
${\Phi}(x)$ is a diffeomorphism on $x\in U$,
$\Gamma_{\rmi\rme}(x): U\rightarrow\RR^{(m-m_\rmd)\times m}$
is smooth
and such that the matrix $\Gamma_\rmi(x)$ is nonsingular,
and $\Gamma_{\rmo\rme}\in\RR^{(p-m_\rmd)\times p}$
is such that the constant matrix $\Gamma_\rmo$ is nonsingular.

The variables $\xi$, $u_\rmd$ and $y_\rmd$ correspond to the structure at infinity,
and the variables $\eta$, $u_\rme$
and $y_\rme$ represent the additional state, input and output variables, respectively,
to form complete sets of state, input and output variables.

Denote
\[
\col\{a_1(x),a_2(x),\cdots,a_{m_\rmd}(x)\}=L_f\Omega_{k^\star}(x),\;
\]
\[
\col\{b_1(x),b_2(x),\cdots,b_{m_\rmd}(x)\}=L_g\Omega_{k^\star}(x).
\]
Let $\delta_{i,i,l}(x)$ be smooth functions and
\[
\delta_i(x)=\mm{
\delta_{i,1,1}(x)&\delta_{i,1,2}(x)&\cdots& \delta_{i,1,i-1}(x)\cr
\delta_{i,2,1}(x)&\delta_{i,2,2}(x)&\cdots& \delta_{i,2,i-1}(x)\cr
\vdots&\vdots& \ddots & \vdots\cr
\delta_{i,q_i-1,1}(x)&\delta_{i,q_i-1,2}(x)&\cdots& \delta_{i,q_i-1,i-1}(x)\cr
0&0&\cdots&0
}, \quad i=1,2,\cdots,m_\rmd.
\]
Define
\[
\col\{\mm{\delta_1(x) &0},\mm{\delta_2(x)&0},\cdots,\delta_{m_\rmd}(x)\}\hspace{20mm}
\]
\be
\hspace{30mm}=\Upsilon
\col\{\mm{\mu_1(x)&0},\mm{\mu_2(x)&0},\cdots,\mu_{{k^\star}}(x)\}, \label{defdelta}
\ee
where
\begin{eqnarray*}
\mu_i(x)&=&\blkdiag\{
R_i S_{i-1\leftrightarrow 2},
R_i S_{i-1\leftrightarrow 3},
\cdots,
R_i,
I_{\rho_i-\rho_{i-1}}\}\\
&&\mm{
P_{1,1}(x) &0 & \cdots& 0\cr
P_{2,1}(x) &P_{2,2}(x)&\cdots &0\cr
\vdots&\vdots&\ddots&\vdots\cr
P_{i-1,1}(x) & P_{i-1,2}(x) & \cdots & P_{i-1,i-1}(x)\cr
0&0&\cdots&0}.
\end{eqnarray*}
We have the following result.
\begin{theo}
\label{theo-form1}
Suppose that system (\ref{nonsys_nfn}) is regular, and that Assumption $\calB$ holds.
Let $q=\{q_1, q_2, \cdots,$ $q_{m_\rmd}\}$
be as obtained in the infinite zero structure algorithm.
Then there exist a set of coordinates in $U$,
i.e., diffeomorphic state, input and output transformations,
%as defined in (\ref{trans1a}) - (\ref{trans1c})
such that
the system takes the following form,
\be\label{form1}
\left\{\begin{array}{rcl}
\dot \eta &=& f_\rme(x)+g_{\rme}(x)u_\rme
           + \displaystyle\sum_{l=1}^{m_\rmd}
            \varphi_l(x)v_{\rmd,l},\cr
\dot\xi_{i,j} &=& \xi_{i,j+1}+ \displaystyle\sum_{l=1}^{i-1}
            \delta_{i,j,l}(x)v_{\rmd,l},\quad j=1,2,\cdots,q_i-1,\cr
\dot\xi_{i,q_i} &=& v_{\rmd,i},\cr
y_\rme &=& h_\rme (x),\cr
y_{\rmd,i} &=& \xi_{i,1}, \qquad i=1,2,\cdots,m_\rmd,
\end{array}\right.
\ee
where $v_{\rmd,i}=a_i(x)+b_i(x)u$, $i=1,2,\cdots,m_\rmd$, with \newline
$\col\{b_1(x),b_2(x),\cdots,b_{m_\rmd}(x)\}$ being nonsingular for $x\in U$,
and
\be\label{key1}
\delta_{i,j,l}(x)=0,\quad \mbox{for } \; j<q_l, \; i=1,2,\cdots,m_\rmd.
\ee
\end{theo}

By (\ref{dzki}), the dynamics $\dot \zeta_{i,j}$ does not relate to the state feedbacks $v_l$ with $l>j$.
Inequality (\ref{key1}) follows from this fact.
Indeed, (\ref{key1}) can be combined into the from (\ref{form1}) by replacing
$\sum_{l=1}^{i-1}\delta_{i,j,l}(x)v_{\rmd,l}$
with $\sum_{l=1,\; q_l\leq j}^{i-1}\delta_{i,j,l}(x)v_{\rmd,l}$.

\begin{rema}
The results of Theorem \ref{theo-form1} are applicable to general MIMO
systems that are not necessarily square. For square and
invertible systems, normal form (\ref{form1}) is in the same form as the
one derived in \cite{isnc99,scgn99},
where no vector relative degree assumption is required either. However, normal form
(\ref{form1}) possesses an extra property (\ref{key1}) (see Example \ref{exx1}).
As will be seen, such a property plays a key role in defining the nonlinear
extension of the infinite zeros of linear systems.
\end{rema}

In what follows, we further simplify the normal form in Theorem~\ref{theo-form1}.

{\it Assumption $\calC$ } :
There exists a $\Gamma_{\rmi\rme}(x)$ in (\ref{trans1b}) such that
the distribution spanned by the column vectors of
$
g_\rmd(x)=g(x)\Gamma_\rmi^{-1}(x)
\mm{0\cr I_{m_\rmd}}
$
is involutive.

\begin{theo}
\label{theo-form2}
Suppose that the conditions in Theorem~\ref{theo-form1} and
{\it Assumption $\calC$} are satisfied. Then there exists
a set of %global
coordinates in $U$
%as defined in (\ref{trans1a}) - (\ref{trans1c})
such that the system takes
the form of Theorem~\ref{theo-form1} with
\be\label{key3}
\varphi_l(x)=0,\quad l=1,2,\cdots,m_\rmd.
\ee
\end{theo}

{\bf Proof:} See Appendices~\ref{p_theo-form2}.
\squ

Let us
apply the infinite zero structure algorithm to a linear system $(A,B,C)$, {\em i.e.}, system (\ref{nonsys_nfn}) with
\[
f(x)=Ax,\quad g(x)=B, \quad h(x)=Cx.
\]
It is obvious that Assumptions $\calA_k$, $k=1,2,\cdots,k^\star$,  $\calB$ and $\calC$ automatically hold.
%Then, we obtain a normal form as a special case of the form given in Theorem~\ref{theo-form2}.

\begin{theo}\label{linear}
Consider a linear system $(A,B,C)$.
There exist nonsingular state, input and output transformations, and a state feedback,
such that the system takes the form,
\be\label{form1_linear}
\left\{\begin{array}{rcl}
\dot \eta &=& A_{11}\eta+A_{12}\col\{\xi_{1,1},\xi_{2,1},\cdots,\xi_{m_\rmd,1}\}+B_1u_\rme,\cr
\dot\xi_{i,j} &=& \xi_{i,j+1}+
\displaystyle \sum_{l=1}^{i-1}
            \delta_{i,j,l} v_{\rmd,l},\quad j=1,2,\cdots,q_i-1,\cr
\dot\xi_{i,q_i} &=& v_{\rmd,i},\cr
y_\rme &=& C_1\eta,\cr
y_{\rmd,i} &=& \xi_{i,1}, \quad i=1,2,\cdots,m_\rmd,
\end{array}\right.
\ee
where
$(A_{11},B_1,C_1)$ does not contain dynamics that is
simultaneously controllable and observable,
$v_{\rmd,i}=a_i\col\{\eta,\xi\}+b_i u$, $i=1,2,\cdots,m_\rmd$, with
$\col\{b_1,b_2,\cdots,b_{m_\rmd}\}$ being nonsingular, and
$
\delta_{i,j,l}=0,\;\; \mbox{for } \; j<q_l, \; i=1,2,\cdots,m_\rmd.
$
\squ
\end{theo}
%{\bf Proof:} See Appendices~\ref{linear-form1}. \squ

The form (\ref{form1_linear}) can be achieved by some additional state transformation
on the linear counterpart of (\ref{form1}).
In (\ref{form1_linear}), the dynamics of $\eta$ depends only on $\eta$ and $\xi_{i,1}$, $i=1,2,\cdots,m_\rmd$, and
$y_\rme$ only depends on $\eta$.
It can be verified that the finite zeros are given
by the simultaneously uncontrollable and unobservable dynamics of $(A_{11},B_1,C_1)$.
The infinite zeros are $\{q_1,q_2,\cdots,q_{m_\rmd}\}$.
The system is
left invertible if $u_\rme$ is absent,
right invertible if $y_\rme$ is absent,
invertible if both $u_\rme$ and $y_\rme$ are absent,
and degenerate if both $u_\rme$ and $y_\rme$ are present.

Some remarks on the infinite zero structure algorithm and normal forms are given as follows.

\begin{rema}
In the infinite zero structure algorithm, there is only one constant rank assumption
in each step, while
in the constrained dynamics algorithm \cite{nind90} and
the zero dynamics algorithm \cite{isnc95,isnc99,scgn99},
each step involves two constant rank assumptions.
However, in the infinite structure algorithm, to construct a new set of coordinates,
Assumption $\calB$ is needed.
Assumption $\calB$ automatically holds if certain matrices are constant (see Lemma~\ref{independent}).
\end{rema}

\begin{rema}
In the structure algorithm, the smooth matrix valued functions $P_{k,l}(x)$, $l=1,2,\cdots,k$, can be found as follows.
By (\ref{Pkj}),
\[
L_g S_k\Theta_{k-1}(x)(L_g\Omega_k(x))^\rmT\hspace{30mm}
\]
\[
\hspace{20mm}
=\mm{P_{k,1}(x)&P_{k,2}(x)&\cdots& P_{k,k}(x)} L_g\Omega_k(x)(L_g\Omega_k(x))^\rmT.
\]
The matrix $L_g\Omega_k(x)$ is of full row rank,
thus
\[
\det(L_g\Omega_k(x)(L_g\Omega_k(x))^\rmT)\neq 0.
\]
Therefore,
\[
\mm{P_{k,1}(x)&P_{k,2}(x)&\cdots& P_{k,k}(x)}\hspace{20mm}
\]
\[
\hspace{20mm}
=[L_g S_k\Theta_{k-1}(x)][L_g\Omega_k(x)]^\rmT
[L_g\Omega_k(x)(L_g\Omega_k(x))^\rmT]^{-1}.
\]
\end{rema}

\begin{rema}
Suppose $dh(x)$ is of full row rank and $g(x)$ is of full column rank
in $U$,
then we can stop repeat Step $k$ in the infinite zero structure algorithm and go to Final Step
if
\[
[\max\{m,p\}-\rho_k-1]+k+\sum_{j=1}^k j(\rho_j-\rho_{j-1})=n,
\]
rather than
\[
k+\sum_{j=1}^k j(\rho_j-\rho_{j-1})=n.
\]
This will lead to fewer steps in the algorithm.
\end{rema}

\begin{rema}
Consider $m=p$.
If we further assume $P_{k,j}(x)$, $j=1,2,\cdots,k,\;$ $k=1,2,\cdots,{k^\star}$
are independent of $x$, we obtain
the structure algorithm of Chapter 5 in \cite{isnc95}.
\end{rema}

\begin{rema}
By Lemma~\ref{independent},
we do not request that the matrices $d \Theta_i(x)$, $i=1,2,\cdots,k^\star$ or their combinations
have constant rank in $U$. By the infinite zero structure algorithm, we always can find
$d \Phi_\rmd(x)$, {\em i.e.}, linear combinations of
$d \Theta_i(x)$, $i=1,2,\cdots,k^\star$, has constant rank. And thus
$d \Phi_\rmd(x)$ can be used as part of the new state coordinate.
\end{rema}

\begin{rema}\label{contin-k}
The infinite zero structure algorithm stops at Step $k^\star$ when (\ref{endalgo}) is satisfied.
Carrying on the algorithm further would not increase $\rho_k$. That is,
$\rho_k=\rho_{k^\star}$,
for $k>k^\star$.
This can be seen in two cases.
Case 1: $\rho_{k^\star}=\min\{p,m\}$.
Suppose there exists a $k^\circ>k^\star$ such that
$\rho_{k^\circ}>\rho_{k^\star}=\min\{p,m\}$.
Then, by the algorithm, $L_g\Omega_{k^\circ}(x)$ is a $\rho_{k^\circ}\times m$ full row rank matrix
and $p>\rho_{k^\circ}$.
This is a contradiction.
Case 2:
$k^\star+\sum_{j=1}^{k^\star} j(\rho_j-\rho_{j-1})=n$.
Suppose there exists a $k^\circ>k^\star$ such that
$\rho_{k^\circ}>\rho_{k^\star}$. Then,
$\sum_{j=1}^{k^\circ} j (\rho_j-\rho_{j-1})>k^\star(\rho_{k^\circ}
-\rho_{k^\star})+\sum_{j=1}^{k^\star} j (\rho_j-\rho_{j-1}) \geq n$.
However, it can be easily verified that
$\col\{d\zeta_1,d\zeta_2,\cdots,d\zeta_{k^\circ}\}_{x=0}$
is a $ (\sum_{j=1}^{k^\circ} j (\rho_j-\rho_{j-1}))\times n$ matrix with a full row rank.
This is also a contradiction.
\end{rema}

\begin{rema}
As observed in \cite{isnc95}, it is in general difficult to construct
a set of coordinates such that (\ref{key3}) is satisfied.
It entails the solution of a system of $n-n_\rmd$ partial differential equations.
However, in the special case that $\varphi_l(x)$, $l=1,2,\cdots,m_\rmd$, in (\ref{form1})
are independent of $x$,
{\em i.e.}, $\varphi_l(x)=\varphi_l$ is
a constant, by renaming the state variable
\[
\bar \eta=\eta-\sum_{l=1}^{m_\rmd} \varphi_l \xi_{l,q_l},
\]
the term $\sum_{l=1}^{m_\rmd}$ $\varphi_l(x)v_{\rmd,l}$ in (\ref{form1}) disappears
under the new set of coordinates.
\end{rema}

\begin{rema}
The variables $\zeta_{i,j},j=1,2,\cdots,i,\,i=1,2,\cdots,{k^\star}$, constitute
all the states associated with the structure at infinity. Note that
for some $i$ with $\rho_i=\rho_{i-1}$, $\zeta_{i,j}$ is not defined.
%The dynamic equations that govern these states,
%(\ref{dzki}) and (\ref{zkk}), represent the system dynamics at infinity.
For each $i=1,2,\cdots,{k^\star}$, the states
$\zeta_{i,j}$, $j=1,2,\cdots,i$, form $\rho_i-\rho_{i-1}$ chains of integrators, and
each chain contains $i$ integrators.
However, except for the smallest $i=i_0$ such that $\rho_{i_0}>0$,
%and $i_0+1$ if $\hat \rho_{i_0+1}\neq 0$,
in which
$\zeta_{i_0,j}$
% and $\zeta_{i_0+1,j}$'s respectively
form $\rho_{i_0}$
%and $\hat \rho_{i_0+1}$
{\em clean} chains of
integrators that link the transformed inputs
% $v_{i_0}$
% and $v_{i_0+1}$
to the
%respective
transformed outputs $\zeta_{i_0,1}$,
for each remaining $i$ with $\rho_i\neq \rho_{i-1}$,
the equations governing the states $\zeta_{i,j}$ represent chains of
$i$ integrators with the previous
transformed inputs $v_l$ ($l<i$)
injected into the integrators with $j\geq q_l$.
\end{rema}

\subsection{Infinite Zeros}

We now extend the linear system notion of infinite zeros to nonlinear systems.
%To motivate this extension,
Consider the normal form in Theorem~\ref{linear}.
The set $q=\{q_1,q_2,\cdots,q_{m_\rmd}\}$ as obtained in the infinite zero structure algorithm coincides with
the infinite zeros of this linear system as defined in
\cite{mosi73,saas87}. This motivates the following definition.

\begin{defi}\label{infinitezero}\em
Suppose that the nonlinear system (\ref{nonsys_nfn}) is regular.
The infinite zeros of
the system are the set of integers
$q=\{q_1, q_2, \cdots, q_{m_\rmd}\}$ as identified in the infinite zero structure algorithm.
\end{defi}

%Thus, for a more specific reference, we refer to the structure algorithm in Section~\ref{infinite_zero_structure_algorithm}
%as the {\em infinite zero structure algorithm}.

Roughly speaking, each integer $q_i$ in the set $q$
represents a chain of integrators of length $q_i$
connecting an input and output pair.
%Lemma~\ref{RS_inv} has shown that $q$ is independent of the choice of $R_k$ and $S_k$ in the algorithm.
We will further justify Definition~\ref{infinitezero} as an extension
of the linear system notion of infinite zeros to nonlinear systems by showing that
the set $q$ is invariant under diffeomorphic state,
input and output transformations, static state feedback and output injection.

Consider a diffeomorphic state transformation $z=\Phi(x)$ in $U$, we have
\be\label{nonsys2}
\left\{\begin{array}{rcl}
\dot{z}&=&\check{f}(z)+\check{g}(z)u,\cr
y&=&\check{h}(z),
\end{array}\right.
\ee
where
\[
\check{f}(z)=\Bigl[  \frac{\partial \Phi}{\partial x} f(x)\Bigl]_{x=\Phi^{-1}(z)},\quad
\check{g}(z)=\Bigl[  \frac{\partial \Phi}{\partial x} g(x)\Bigl]_{x=\Phi^{-1}(z)},\quad
\]
\[
\check{h}(z)=[ h(x)]_{x=\Phi^{-1}(z)}. \hspace{60mm}
\]
Following the infinite zero structure algorithm, it is easy to verify the following result.
\begin{lemma}
\label{diffe}
If system (\ref{nonsys_nfn}) is regular, then
system (\ref{nonsys2}) is regular too.
Moreover, both systems have the same infinite zeros.
\end{lemma}

We also have the following result.
\begin{lemma}
\label{invariant}
The infinite zeros of system (\ref{nonsys_nfn})
are invariant under
\begin{enumerate}
\item input transformation $\check{u}=\Gamma_\rmi(x)u$ with $\Gamma_\rmi(x): U\rightarrow \RR^{m\times m}$
being smooth and nonsingular;
\item output transformation
$\check{y}=\Gamma_\rmo\, y$ with $\Gamma_\rmo\in\RR^{p\times p}$ being nonsingular;
\item static state feedback $\check{u}=u-K(x)$ with $K(x): U\rightarrow \RR^m$ being smooth; and
\item output injection, {\em i.e.},
\be\label{nonsys4}
\left\{\begin{array}{rcl}\dot x &=& f(x)+F(x)h(x)+g(x)u,\\
                           y &=& h(x),\end{array}\right.
\ee
where $F(x): U\rightarrow\RR^{n\times p}$ is smooth.

\end{enumerate}
\end{lemma}

{\bf Proof:}
See Appendices~\ref{p_invariant}. \hfill $\square$

\subsection{Invertibility and Zero Dynamics}

Equations (\ref{endalgo}) in the infinite zero structure algorithm indicates the invertibility property of the system.
\begin{lemma}
System (\ref{nonsys_nfn}) is
left invertible if $\rho_{k^\star}=m<p$,
right invertible if $\rho_{k^\star}=p<m$,
invertible if $\rho_{k^\star}=m=p$,
and degenerate if $\rho_{k^\star}<\min\{m,p\}$.
\end{lemma}

Equivalently,
the system in Theorem~\ref{theo-form1} is
left invertible if $u_\rme$ is absent,
right invertible if $y_\rme$ is absent,
invertible if both $u_\rme$ and $y_\rme$ are absent,
and degenerate if both $u_\rme$ and $y_\rme$ are present.

%\subsection{The Zero Dynamics}

In \cite{isnc95}, the zero dynamics of a nonlinear system is defined for a square invertible nonlinear system.
Let $M$ be a smooth connected submanifold of $U$.
The manifold $M$ is said to be locally controlled invariant at $ x=0 $ if there
exist a smooth mapping $u: M \rightarrow \RR^m$ and a neighborhood $U^\circ$ of $x=0$ such
that $M$ is locally invariant under the vector field $f(x)+g(x)u(x)$.
A zero output submanifold in a neighborhood of $x=0$
for the nonlinear system (\ref{nonsys_nfn}) is a smooth
connected submanifold $M$,
which is locally controlled invariant at $x=0$
and
for each $x\in M$, $h(x)=0$.
Suppose $Z^\star$ is the locally maximal zero output submanifold with
$
\span\{g( 0 )\}\cap T_{ x=0 }Z^\star=0,
$
where $T_{ x=0 }Z^\star$ represents the tangent space to $Z^\star$ at $x=0$.
Then, there exists a unique smooth mapping $u^\star : Z^\star \rightarrow \RR^m$ such that
the vector field
$
f^\star(x)=f(x)+g(x)u^\star(x)
$
is tangent to $Z^\star$. The pair $(Z^\star, f^\star)$ is called the zero dynamic of (\ref{nonsys_nfn}).

The global version of $Z^\star$ for a square invertible nonlinear system
is defined in \cite{isnc99,scgn99} as a controlled invariant smooth embedded submanifold of $\RR^n$.

Here, we want to use the form (\ref{form1}) to derive the zero dynamics of general nonlinear system in $U$.
In particular,
%suppose $f(0)=0$ and $h(0)=0$ and $0\in U$.
let $y_\rmd=0$ in (\ref{form1}).
It then follows from the dynamic equations that $\xi=0$ and $v_{\rmd,i}=0$, $i=1,2,\cdots,m_\rmd$.
Consequently, the remaining dynamics reduces to
\be\label{abc}
\left\{\begin{array}{rcl}
\dot \eta &=& f_{\rme}(\eta, 0)+g_{\rme}(\eta, 0) u_\rme,\cr
y_\rme &=& h_\rme(\eta, 0).\cr
\end{array} \right.
\ee
Let $\calC_0$ be the smallest distribution that is invariant for (\ref{abc})
and contains the distribution spanned by
the column vectors of $g_{\rme}(\eta, 0)$, and $d\calO$ be
the smallest codistribution that is invariant for (\ref{abc}) and contains
the codistribution spanned by the row vectors of $dh_\rme(\eta, 0)$.
Note that the distribution $\calC_0$ characterizes local strong accessibility and
the codistribution $d\calO$ characterizes local observability.
The subsystem (\ref{abc}) does not contain any subspace that is both
strong locally accessible (by $u_\rme$) and locally observable (through $y_\rme$).
Otherwise, the infinite zeros are no longer $q=\{q_1, q_2, \dots, q_{m_\rmd}\}$.
Thus by \cite{isnc95,nind90}, we have the following result.

\begin{lemma}\label{abcdecompose}
Consider system (\ref{abc}). % which is obtained by letting $\xi=0$ in (\ref{form1}).
Assume that the distributions
$\calC_0$, $\ker d\calO$ and $\calC_0+\ker d\calO$ of (\ref{abc}) each has a constant dimension.
Then there exist a set of coordinates $\hat{z}=\col\{z_\rma,z_\rmb,z_\rmc\}$ such that
(\ref{abc}) takes the form
\be\label{zrmd0_3}
\left\{\begin{array}{rcl}
\dot{z}_{\rma}&=&f_{\rma}(z_{\rma},z_{\rmb}),\cr
\dot{z}_{\rmb}&=&f_{\rmb}(z_{\rmb}),\cr
\dot{z}_{\rmc}&=&f_{\rmc}(z_{\rma}, z_{\rmb}, z_{\rmc})
+g_{\rmc\rme}(z_{\rma},
z_{\rmb}, z_{\rmc})u_{\rme},\cr
y_{\rme} &=& h_{\rme\rmb}(z_{\rmb}),
\end{array}\right.
\ee
with
$\calC_0=\span\{\frac{\partial}{\partial z_\rmc}\}$
and
$\ker d\calO=\span\{\frac{\partial}{\partial z_\rma},\frac{\partial}{\partial z_\rmb}\}$.
\hfill $\square$
\end{lemma}

%Roughly spoken, zero dynamics is the dynamics that is neither controllable nor observable
%in the subsystem (\ref{abc}).
%Suppose that under a new set of state variables
%$({z}_{\rma},{z}_{\rmb},{z}_{\rmc})$, the system
%(\ref{abc}) is written in the following form (see, for example,
%\cite{nind90}),
%where the state $z_{\rma}$ is neither controllable nor observable,
%$z_{\rmb}$ is observable but not controllable, and $z_\rmc$ is controllable
%but not observable. Consequently,

The decomposition (\ref{zrmd0_3})
%, which reflects the system invertibility properties,
allows us to decompose normal form (\ref{form1}) into
four distinct subsystems (see Example \ref{exx2}) as we can do in a
linear system (\cite{mosi73,saas87}). In a generalization to the notion of invariant
zero of  linear systems \cite{saas87}, the dynamics
$\dot{z}_{\rma}=f_{\rma}(z_{\rma},0)$
is referred to as the zero dynamics of system (\ref{nonsys_nfn}).
The case of $m=p=m_\rmd=\rho_{k^\star}$ has been
studied in \cite{isnc95,scgn99}.
In this case, $y_\rme$ and $u_\rme$ are absent from (\ref{form1}),
and $\dot \eta = f_{\rme}(\eta,0)$
is directly obtained
as the zero dynamics of system (\ref{nonsys_nfn}).

\subsection{Normal Forms of Square Invertible Systems}

We now consider the normal forms of system (\ref{nonsys_nfn}) with
$m=p=\rho_{k^\star}=m_\rmd$, {\em i.e.}, a square invertible system, which has been considered
in \cite{isnc95,isnc99,scgn99}. In this case, $y_\rme$ and $u_\rme$ do not exist and we have the
following result, as a corollary to Theorems~\ref{theo-form1} and \ref{theo-form2}.
\begin{coro}\label{square_nonlin}
Suppose that a square invertible system (\ref{nonsys_nfn})
has infinite zeros $q=\{q_1, q_2, $ $\cdots, q_{m}\}$,
Assumption $\calB$ holds, and
the distribution spanned by the column vectors of $g(x)$ is involutive.
Under a new set of coordinates,
%defined by (\ref{trans1a}) - (\ref{trans1c}),
the system takes the form,
\be\label{form1-square_nonlin}
\left\{\begin{array}{rcl}
\dot \eta &=& f_\rme(x)
           + \displaystyle \sum_{l=1}^{m}
            \varphi_l(x) v_l,\cr
\dot\xi_{i,j} &=& \xi_{i,j+1}+ \displaystyle \sum_{l=1}^{i-1}
            \delta_{i,j,l}(x)v_l,\quad j=1,2,\cdots,q_i-1,\cr
\dot\xi_{i,q_i} &=& v_i,\cr
 y_i &=& \xi_{i,1}, \quad i=1,2,\cdots,m,
\end{array}\right.
\ee
where
$v_i=a_i(x)+b_i(x)u$ with
$\col\{b_1(x),b_2(x),\cdots,b_m(x)\}$ being nonsingular,
and
\be\label{key1-square}
\delta_{i,j,l}(x)=0,\;\; \mbox{ for }\; j<q_l,\; i=1,2,\cdots,m.
\ee
If, in addition, the distribution spanned by the column vectors of $g(x)$ is involutive for $x\in U$,
then there exist a set of coordinates such that
\[
\varphi_l(x)=0, \quad l=1,2,\cdots,m.
\]
\squ
\end{coro}

Note that the form given in Corollary~\ref{square_nonlin}
is the same as (\ref{zerody2}) except for the additional structural property (\ref{key1-square}).
The $\dot\xi_{i,j}$ equation in (\ref{form1-square_nonlin}) displays a triangular structure of the
control inputs that enter the system. Property (\ref{key1-square}) imposes an additional
structure within each chain of integrators on how control inputs enter the system.
With this additional structural property,
the set $q=\{q_1, q_2, \cdots, q_{m_\rmd}\}$ represents infinite zeros
when the system is linear.
Note that the property (\ref{key1-square})
which can be deduced from (\ref{dzki}),
is a key feature which the form (\ref{zerody1}) resulting from the algorithm
in \cite{byls88,isnc95,isnc99,scgn99} does not possess.
To see the significance of property (\ref{key1-square}), we transfer
system (\ref{counter3}) into the normal form (\ref{form1-square_nonlin}),
\[
\dot{\tilde{x}}=
\mm{
0&0&0&0&0\cr
0&0&1&0&0\cr
0&0&0&1&0\cr
0&0&0&0&1\cr
0&0&0&0&0}\tilde{x}+
\mm{
1&0\cr
0&0\cr
1/\alpha &0\cr
0&0\cr
0&1}
\pmatrix{v_{1}\cr v_{2}},\quad
\]
\[
\tilde{y}
=\mm{1&0&0&0&0\cr 0&1&0&0&0}\tilde{x},
\]
by using the following state and output transformations and state feedback,
\vspace{-4pt}
\[
\tilde{x}=\mm{
0&0&1&0&0\cr
1&0&0&0&0\cr
0&1&0&0&0\cr
0&0&0&-1/\alpha&0\cr
0&0&0&0&-1/\alpha}x,\quad
\tilde{y}=\mm{0&1\cr 1&0} y,\quad
\]
\[
\pmatrix{v_{1}\cr v_{2}}=
\pmatrix{
-\alpha\tilde{x}_4+\alpha u_1\cr
u_2
}.
\]
It is obvious that $q_1=1$ and $q_2=4$, which coincide with the infinite zeros of
this linear system (see, {\em e.g.}, \cite{mosi73,saas87}).

\begin{rema}
The normal form in Corollary~\ref{square_nonlin}
can be further simplified by using the method in \cite{isnc99,scgn99}.
Define the vector fields
$
Y_{k,j}(x), \; 1\leq j\leq m, \; 1\leq k\leq q_j.
$
If these vector fields commute, then there exist a set of coordinates such that
the dynamics of $\eta$ in Corollary~\ref{square_nonlin} simplifies to
$
\dot \eta = f_\rme(\eta,\xi_{1,1},\xi_{2,1},\cdots,\xi_{m,1}).
$
\end{rema}

  Next, we follow the method in \cite{isnc99,scgn99} to further simplify
  the dynamic of $\eta$ in Corollary~\ref{square_nonlin}. Define
  \begin{eqnarray*}
  \tilde{f}(x)&=&f(x)-g(x)\col\{b_1(x),b_2(x),\cdots,b_m(x)\}^{-1}
  \col\{a_1(x),a_2(x),\cdots,a_m(x)\},\\
  \tilde{g}(x)&=&g(x)\col\{b_1(x),b_2(x),\cdots,b_m(x)\}^{-1},
  \end{eqnarray*}
  and let
  \[
  Y_m^k(x)=(-1)^{k-1}\ad_{\tilde{f}}^{k-1}\tilde{g}_{m}(x), \quad 1\leq k\leq q_m,
  \]
  and for $1\leq j\leq m-1$, $1<k\leq q_j$,
  \begin{eqnarray*}
  Y_j^1(x)&=&\tilde{g}_j(x)-\sum_{l=j+1}^{m_\rmd} \sum_{i=2}^{q_l} \delta_{l,q_l-i+1,j}(x)Y_l^i(x),\\
  Y_j^k(x)&=&(-1)^{k-1}\ad_{\tilde{f}}^{k-1}Y_j^1(x).
  \end{eqnarray*}

  {\em Assumption $\calD$ : %There exist a local coordinates defined by (\ref{trans1a}) - (\ref{trans1c}) such that
  The vector fields
  $
  Y_j^k(x), \; 1\leq j\leq m, \; 1\leq k\leq q_j,
  $
  commute, {\em i.e.},
  \[
  [Y_i^s,Y_j^k]=\frac{\partial Y_j^k}{\partial x} Y_i^s -\frac{Y_i^s}{\partial x} Y_j^k=0,\quad
  1\leq i,j\leq m,\;\; 1\leq s\leq q_i,\;\; 1\leq k\leq q_j.
  \]
  }

  The following result is immediate from \cite{isnc99}.

  \begin{theo}\label{square_nonlin_strict_2}
  Suppose that the conditions of Corollary~\ref{square_nonlin}
  and Assumption $\calD$ hold,
  then the dynamics of $\eta$ in Corollary~\ref{square_nonlin} can be simplified to
  $
  \dot \eta = f_\rme(\eta,\xi_{1,1},\xi_{2,1},\cdots,\xi_{m,1}).
  $
  \end{theo}

\section{Normal Forms of Nonlinear Systems Relating to the Zero Output}
\label{algo_zerooutput}

In determining the zero dynamics, a normal form representation of the nonlinear system is also given in Chapter 6 of \cite{isnc95},
which displays structure information along one special output trajectory, the zero output.
Two constant rank assumptions are made in the nested submanifolds $M_k$,
$k=1,2,\cdots,k^\star$.

Here, we will show that the infinite zero structure algorithm can be adapted for the same problem.
In particular, for system (\ref{nonsys_nfn}), we will introduce
{\it Assumption ${\bar\calA}_k$}, $k=1,2,\cdots,k^\star$,
in the nested subsets $M_k$,
$k=1,2,\cdots,k^\star$, rather than for all $x\in U$.
Because the nested subsets $M_k$, $k=1,2,\cdots,k^\star$, are related to the zero output,
we refer to the resulting
algorithm as the zero output structure algorithm.
%It reveals structure properties of nonlinear system with the constraint that the output is identically zero.
%We call it as a structure algorithm, rather than zero dynamics algorithm,
%because it also reveals some system structure properties.

{\bf Zero Output Structure Algorithm}

{\bf Initial Step.}
Let
$\Theta_0(x)=h(x)$, $\Omega_0(x)=\emptyset$, $\rho_0=0$ and $k=1$.
% and repeat the following procedure.

{\bf Step $k$.} We start with
$\Theta_{k-1}(x): U\rightarrow \RR^{p-\rho_{k-1}}$ and
$\Omega_{k-1}(x): U\rightarrow \RR^{\rho_{k-1}}$,
where the matrix $L_g\Omega_{k-1}(x)$ has full row rank $\rho_{k-1}$ in $M_k\cap O_k$,
with
$M_k=\{x: \; \Theta_i(x)=0,\; i=0,1,\cdots,k-1.\}$ and
$O_k$ being a neighborhood of $x=0$.

{\em Assumption $\bar \calA_k$:}
{\em The matrix
$\mm{L_g\Omega_{k-1}(x)\cr L_g\Theta_{k-1}(x)}$
has a constant rank $\rho_k$ in $M_k\cap O_k$, and there exists an
$R_k\in\RR^{(\rho_k-\rho_{k-1})\times(p-\rho_{k-1})}$ such that
\be\label{rbHk2}
\rank \left(\mm{L_g\Omega_{k-1}(x)\cr L_g R_k \Theta_{k-1}(x)}\right)=\rho_k,
\quad \forall\; x\in (M_k\cap O_k)^c,
\ee
where $(M_k\cap O_k)^c$ is the connected component of $M_k\cap O_k$ containing  $x=0$.}

Suppose that {\em Assumption $\bar \calA_k$} is satisfied.
Let $S_k
%\in\bbR^{(p-\rho_k)\times(p-\rho_{k-1})}
$ and $\Omega_k(x)$ be as in (\ref{RkSk})-(\ref{Omega_k}).
Thus, the matrix $L_g \Omega_k(x)$ has full row rank $\rho_k$ for $x\in (M_k\cap O_k)^c$, and
\[
\rank\left(\mm{L_g\Omega_k(x)\cr L_g S_k\Theta_{k-1}(x)}\right)=\rho_k,
\quad \forall\; x\in (M_k\cap O_k)^c.
\]
%{\em i.e.}, the rows of $L_g  S_k\Theta_{k-1}(x)$ are in the space spanned by
%the rows of $L_g R_i\Theta_{i-1}(x), i=1,2,\cdots, k$, in $N_k$.
Therefore, there exist smooth functions
$P_{k,l}(x)\in\RR^{(p-\rho_k)\times(\rho_l-\rho_{l-1})},l=1,2,\cdots,k$, such that
\be\label{Pkj2}
L_g  S_k\Theta_{k-1}(x)
-\sum_{l=1}^k P_{k,l}(x)L_g R_l\Theta_{l-1}(x)-W_k(x)=0,
\ee
where $W_k(x)$ is a matrix valued smooth function with $W_k(x)=0$ in $(M_k\cap O_k)^c$.
Denote
\be
v_k=L_f R_k\Theta_{k-1}(x)+L_g  R_k\Theta_{k-1}(x)u,\label{vk2}
\ee
and define
\begin{eqnarray}\label{Thetak2}
\Theta_k(x)
&=&\displaystyle\frac{d}{dt} S_k\Theta_{k-1}(x)
    -\sum_{l=1}^k P_{k,l}(x)v_l-W_k(x)u\nonumber\\
&=& L_f \displaystyle  S_k\Theta_{k-1}(x)
-\sum_{l=1}^k P_{k,l}(x)L_f R_l\Theta_{l-1}(x).
\end{eqnarray}
If $k+\sum_{j=1}^k j(\rho_j-\rho_{j-1})<n$ and $\rho_{k}<\min\{p,m\}\;$,
then increase $k$ by $1$ and repeat the above step. Otherwise, go to Final Step.

{\bf Final Step.}
The same as the final step in the infinite zeros structure algorithm in Section~\ref{infinite_zero_structure_algorithm}.

{\bf End.}

\begin{defi}\em
The point $x=0$ is said to be a regular point of system (\ref{nonsys_nfn}) if
Assumption $\bar\calA_k$, $k=1,2,\cdots, k^\star$,
in the zero output structure algorithm
are satisfied.
\end{defi}

Note that in Step $k$,
the choice of the matrices $R_k$ and $S_k$, which satisfy (\ref{rbHk2}) and (\ref{RkSk}),
are not unique.
\begin{lemma}
\label{preinv}
Suppose that $\check{R}_i$ and $\check{S}_i$,
$i=1,2,\cdots,k^\star$, are different choices
yielding $\check{\Theta}_i(x)$, $\check{\Omega}_i(x)$, $\check{M}_i$ and $\check{O}_i$.
Then,
\be\label{pkinv2}
\check{M}_i=M_i,\quad
\check{\Theta}_i(x)=\sum_{l=1}^{i-1}Q_{i,l}(x)\Theta_l(x)+T_i(x)\Theta_i(x)+V_i(x),
%i=1,2,\cdots,k^\star,
\ee
%\[
%d\check{\Theta}_i(x)=Q_i(x)dG_i(x)+T_i(x)d\Theta_i(x)+Q_i(x)
%\]
%\[
%L_f\check{\Theta}_i(x)=Q_i(x)L_fG_i(x)+T_i(x)L_f\Theta_i(x)+Q_i(x)f(x)
%\]
%\[
%L_g\check{\Theta}_i(x)=Q_i(x)L_g G_i(x)+T_i(x)L_g\Theta_i(x)+Q_i(x)g(x)
%\]
where $T_i(x)$ is
a nonsingular matrix valued smooth function, and $V_i(x)$ is
smooth with
$V_i(x)=0$ in $M_i \cap \check{O}_i\cap O_i $.
\end{lemma}
{\bf Proof:} See Appendices~\ref{p_preinv}. \squ

The following result follows directly from Lemma~\ref{preinv}.

\begin{lemma}
\label{RS_inv}
%The integers $\{\rho_1,\rho_2,\cdots,\rho_{k^\star}\}$,
The set $\rho$,
and hence the set $q$,
%=\{q_1,q_2,\cdots,q_{m_\rmd}\}$,
as identified in
the zero output structure algorithm are invariant with respect to
the choice of matrices $ R_i$ and $ S_i,\; i=1,2,$ $ \cdots, k^\star$.
\end{lemma}

Define $\bar\Phi_\rmd(x)$, $\Gamma_\rmi(x)$ and $\Gamma_\rmo$ as in (\ref{Phi1})-(\ref{Gammaod}).
We have the following crucial result.

\begin{lemma}
\label{pindependent2}
Let $x=0$ be a regular point of system (\ref{nonsys_nfn}).
Then, $d\bar{\Phi}_\rmd(0)$,
${\Gamma}_{\rmI\rmd}(0)$ and $\Gamma_{\rmO\rmd}$ are of full row rank.
\end{lemma}
{\bf Proof:} See Appendices~\ref{p_independent2}.
\hfill $\square$

\begin{theo}
\label{theo-form2-nested}
Consider system (\ref{nonsys_nfn}).
Suppose that $x=0$ is a regular point.
%and that the matrix $\bar\Phi_\rmd(x)$ is of full row rank.
Let $q=\{q_1, q_2, \cdots, q_{m_\rmd}\}$
be as obtained in the zero output structure algorithm.
There exist a set of coordinates,
i.e., diffeomorphic state, input and output transformations,
%as defined in (\ref{trans1a}) - (\ref{trans1c})
such that
the system assumes the following form,

\be\label{form2-nested}
\left\{\begin{array}{rcl}
\dot \eta &  =  & f_\rme(x) + g_{\rme}(x)u_\rme
             +  \displaystyle\sum_{l=1}^{m_\rmd}
            \varphi_l(x)v_{\rmd,l},\cr
  \dot\xi_{i,j} &  =  & \xi_{i,j+1}+\displaystyle\sum_{l=1}^{i-1}
            \delta_{i,j,l}(x)v_{\rmd,l}+\sigma_{i,j}(x)u,\quad j=1,2,\cdots,q_i-1,\cr
   \dot\xi_{i,q_i} &  =  & v_{\rmd,i},\cr
y_\rme &  =  & h_\rme (x),\cr
y_{\rmd,i} &  =  & \xi_{i,1}, \qquad i=1,2,\cdots,m_\rmd,
\end{array}\right.
\ee
where $\sigma_{1,j}(x)=0,\; j=1,2,\cdots,n_1-1$,
and $\sigma_{i,j}(x)=0,\; i>1, \;j=1,2,\cdots,n_i-1$, in $(M_j\,\cap\,O_j)^c$,
$v_{\rmd,i}=a_i(x)+b_i(x)u$, $i=1,2,\cdots,m_\rmd$, with
$\col\{b_1 (x),b_2(x),\cdots,b_{m_\rmd}(x)\}$ being nonsingular,
and
\be\label{key2}
\delta_{i,j,l}(x)=0,\quad \mbox{for } \; j<q_l, \; i=1,2,\cdots,m_\rmd.
\ee
\end{theo}

We next consider system (\ref{nonsys_nfn}) with
$m=p=m_\rmd=\rho_{k^\star}$ in the zero output structure algorithm, which has been considered
in \cite{isnc95,isnc99,scgn99}. In this case, $y_\rme$ and $u_\rme$ do not exist and we have the
following result.

\begin{coro}\label{square_nonlinear_nested}
Suppose that the conditions in Theorem~\ref{theo-form2-nested} hold with
$m=p=m_\rmd=\rho_{k^\star}$.
Then, there exist a set of local coordinates
such that
the system takes the form,
\be\label{form1-square2}
\left\{\begin{array}{rcl}
\dot \eta &=& f_\rme(x)
           +\displaystyle\sum_{l=1}^{m}
            \varphi_l(x) v_l,\cr
\dot\xi_{i,j} &=& \xi_{i,j+1}+\displaystyle\sum_{l=1}^{i-1}
            \delta_{i,j,l}(x)v_l+\sigma_{i,j}(x)u,\quad j=1,2,\cdots,q_i-1,\cr
\dot\xi_{i,q_i} &=& v_i,\cr
 y_i &=& \xi_{i,1}, \quad i=1,2,\cdots,m,
\end{array}\right.
\ee
where $\sigma_{1,j}(x)=0,\; j=1,2,\cdots,n_1-1$,
and $\sigma_{i,j}(x)=0,\; i>1, \;j=1,2,\cdots,n_i-1$, in $(M_j\,\cap \, O_j)^c$,
$v_i=a_i(x)+b_i(x)u$ with\newline
$\col\{b_1(x),b_2(x),\cdots,b_m(x)\}$ being nonsingular,
and
\be\label{key1-square2}
\delta_{i,j,l}(x)=0,\;\; \mbox{ for }\; j<q_l,\; i=1,2,\cdots,m.
\ee
The submanifold $Z^\star$ is given as
$
Z^\star=\{x\in U\, : \, \xi_{i,j}(x)=0, \, j=1,2,\cdots, q_i,\, i=1,2,\cdots,m.\}.
$
\squ
\end{coro}

\begin{rema}
Corollary~\ref{square_nonlinear_nested} is the same as the result in Chapter 6 of \cite{isnc95},
except that there is
property (\ref{key1-square2}) here.
\end{rema}

\begin{rema}
The zero output structure algorithm requires milder regularity assumptions
than the infinite zero structure algorithm.
For example, consider
\[
f(x)=\mm{0\cr 0},\quad g(x)=\mm{1&0\cr 0&1},\quad h(x)=\mm{x_1\cr x_1x_2}.
\]
By the infinite zero structure algorithm, the system is not regular, since
\[
L_g h(x)=\mm{1&0\cr x_2& x_1}
\]
does not have a constant rank in a neighborhood of $x=0$.
However, by the zero output structure algorithm, $x=0$ is a regular point with $\rho=\{1,1\}$,
and the locally maximal zero output submanifold is $Z^\star=\{0\}$.
\end{rema}

\begin{rema}
We have similar results as in Lemmas~\ref{diffe}, \ref{invariant} and \ref{abcdecompose}
for zero output structure algorithm.
If the point $x=0$ of system (\ref{nonsys_nfn}) is regular,
then the point $\Phi( 0 )$ of system (\ref{nonsys2}) is regular too.
The set of the integers $q$ as identified in the zero output structure algorithm
are invariant under the state, input and output transformations, state feedback
and output injection as defined in Lemma~\ref{invariant}.
The zero dynamics can be computed similarly as in Lemma~\ref{abcdecompose}.
\end{rema}

\section{Examples}
\label{example_nfn}

Examples~\ref{exx1} and \ref{exx2}
illustrate the infinite zero structure algorithm,
and Example~\ref{exx3}
illustrates the zero output structure algorithm.
and Example~\ref{appp} is an application of our results to a practical system.

\begin{example}\label{exx1}\rm
Consider system (\ref{nonsys_nfn}) with
\[
f(x)=\pmatrix{x_3\cr x_5\cr x_1\cr x_1x_2\cr x_4},\quad
g(x)=\mm{0& 0\cr 0& 0\cr 1&x_3\cr 0&1\cr x_4&x_3 x_4},\quad
h(x)=\pmatrix{x_1\cr x_2},
\]
and $U=\{x: \, x_1<1\}$. We carry out the infinite zero structure algorithm as follows.

{\bf Initial Step.}
Let $\Theta_0(x)=h(x)$, $\Omega_0(x)=\emptyset$, $\rho_0=0$ and $k=1$.

{\bf Step 1.}
\[
L_f\Theta_0(x)=\pmatrix{x_3\cr x_5},\quad
L_g\Theta_0(x)=\mm{0&0\cr 0&0}.
\]
Hence, $\rho_1=0$. Let
\[
R_1=\emptyset, \quad S_1=I_2.
\]
Thus,
\[
v_1=\emptyset,\quad \Omega_1(x)=\emptyset,\quad
P_{1,1}(x)=\emptyset,\quad
\Theta_1(x)=\col\{x_3, x_5\}.
\]

{\bf Step 2.}
\[
L_f\Theta_1(x)=\pmatrix{x_1\cr x_4},\quad
L_g\Theta_1(x)=\mm{1 & x_3\cr x_4 & x_3x_4}.
\]
Hence, $\rho_2=1$. Let
\[
R_2=\mm{1 &0},\quad S_2=\mm{0&1}.
\]
Thus,
\[
v_2=x_1+\mm{1 & x_3}u,\quad
\Omega_2(x)=x_3,
\]
\[
P_{2,1}(x)=\emptyset,\quad
P_{2,2}(x)=x_4,\quad
\Theta_2(x)=x_4-x_1 x_4.
\]
{\bf Step 3.}
\[
L_f\Theta_2(x)=x_1x_2-x_3x_4-x_1^2 x_2,\quad
L_g\Theta_2(x)=\mm{0 & 1-x_1}.
\]
Hence, $\rho_3=2$. Let
\[
R_3=1,\quad
S_3=\emptyset.\quad
v_3=x_1x_2-x_3x_4-x_1^2 x_2+\mm{0 & 1-x_1}u.
\]

{\bf Final Step.} $k^\star=3$. $m_\rmd=2$. $n_\rmd=5$,
$\rho=\{0,1,2\}$.
$q=\{2,3\}.$

Let
\[
\col \{\xi_{1,1}, \xi_{1,2}, \xi_{2,1}, \xi_{2,2}, \xi_{2,3}\}
=\col \{x_1, x_3, x_2, x_5, x_4-x_1x_4\},
\]
\[
y_{\rmd,1}=y_1,\quad
y_{\rmd,2}=y_2.
\]
\[
v_{\rmd, 1}=x_1+\mm{1 & x_3}u, \quad
v_{\rmd, 2}=x_1x_2-x_3x_4-x_1^2 x_2+\mm{0 & 1-x_1}u.
\]
%Thus
%$
%x=\col \{\xi_{1,1}, \xi_{2,1}, \xi_{1,2}, \frac{\xi_{2,3}}{1-\xi_{1,1}}, \xi_{2,2}\}.
%$
The form (\ref{form1}) is given by
\[
\left\{\begin{array}{rcl}
\dot{\xi}_{1,1}&=&\xi_{1,2},\cr
\dot{\xi}_{1,2}&=&v_{\rmd,1},\cr
\dot{\xi}_{2,1}&=&\xi_{2,2},\cr
\dot{\xi}_{2,2}&=&\xi_{2,3}+ \frac{\xi_{2,3}}{1-\xi_{1,1}} v_{\rmd,1},\cr
\dot{\xi}_{2,3}&=&v_{\rmd,2},\\
\rule{0mm}{5mm}
y_{\rmd,1}&=&\xi_{1,1},\cr
y_{\rmd,2}&=&\xi_{2,1}.
\end{array}\right.
\]
The system is invertible with two infinite zeros of order $2$ and $3$.
The zero dynamic degenerates to the single point $x=0$. Note that $\delta_{2,1,1}(x)=0$,
{\em i.e.}, the term of $v_{\rmd,1}$ does not appear in the dynamic equation of $\dot\xi_{2,1}$.
\end{example}

\begin{example}\rm
\label{exx2}

Consider system (\ref{nonsys_nfn}) with
\[
f(x)=\pmatrix{-x_1+x_3\cr
x_2x_4\cr
-x_2x_4-x_2x_4^2\cr
-x_4
},
\quad
g(x)=\mm{
x_2&e^{-x_4}\cr
0&0\cr
0&e^{-x_4}\cr
1&0},
\quad
h(x)=\pmatrix{x_2\cr x_4}.
\]
The system is defined globally, {\em i.e.}, $U=\RR^4$.
%It is obvious that $f(0)=0$ and $h(0)=0$, {\em i.e.}, $x=0$.
We apply the infinite zero structure algorithm.

{\bf Initial Step.}
Let $\Theta_0(x)=h(x)$, $\Omega_0(x)=\emptyset$, $\rho_0=0$ and $k=1$.

{\bf Step 1.}
\[
L_f \Theta_0(x)=\pmatrix{x_2x_4\cr -x_4},\quad
L_g\Theta_0(x)=\mm{0&0\cr 1&0}.
\]
So, $\rho_1=1$. Let
\[
R_1=\mm{0 &1},\quad
S_1=\mm{1&0}.
\]
Thus,
\[
v_1=-x_4+u_1,\quad
\Omega_1(x)=x_4,\quad
P_{1,1}(x)=0,\quad
\Theta_1(x)=x_2x_4.
\]

{\bf Step 2.}
\[
L_f\Theta_1(x)=-x_2x_4+x_2x_4^2,\quad
L_g\Theta_1(x)=\mm{x_2& 0}.
\]
So, $\rho_2=1$.
Let
$
R_2=\emptyset,\;\; S_2=1.
$
Thus,
\[
v_2=\emptyset,\quad
\Omega_2(x)=x_4,\quad
P_{2,1}(x)=x_2,\quad
P_{2,2}(x)=\emptyset,\quad
\Theta_2(x)=x_2x_4^2.
\]

{\bf Step 3.}
\[
L_f\Theta_2(x)=-2x_2x_4^2+x_2x_4^3,\quad
L_g\Theta_2(x)=\mm{2x_2x_4& 0}.
\]
So, $\rho_3=1$.
Let
\[
R_3=\emptyset,\quad
S_3=1,\quad
v_3=\emptyset,\quad
\Omega_3(x)=x_4,
\]
\[
P_{3,1}(x)=2x_2x_4,\quad
P_{3,2}(x)=\emptyset,\quad
P_{3,3}(x)=\emptyset,\quad
\Theta_3(x)=x_2x_4^3.
\]

{\bf Final Step.} $k^\star=3$, $m_\rmd=1$, $n_\rmd=1$,
$\rho=\{1,1,1\}$ and $q=\{1\}$.
%Define $\zeta_{1,1}=x_4$, we have $\dot{\zeta}_{1,1}=v_1.$
It is obvious that
$
\Phi_\rmd(x)=x_4,\;\;
\Gamma_{\rmi\rmd}(x)=\mm{1& 0},
\;\;
\Gamma_{\rmo\rmd}(x)=\mm{0 &1}.
$
Let
\[
\Gamma_\rmi(x)=\mm{0& e^{-x_4}\cr 1&0},\quad
\Gamma_\rmo=\mm{1&0\cr 0&1}.
\]
Thus, $g_\rmd(x)=\col\{x_2,0,0,1\}$.
Find $\Phi_\rme$ such that $d\Phi_\rme\, g_\rmd(x)=0$, {\em i.e.},
\[
\displaystyle \frac{\partial \Phi_\rme}{\partial x_1}x_2+\frac{\partial \Phi_\rme}{\partial x_4}=0.
\]
We obtain $\Phi_\rme(x)=\col\{x_1-x_2x_4,x_2,x_3\}$.
Let
$
\col\{\eta_1, \eta_2, \eta_3, \xi\}%=\Phi(x)
=\col\{x_1-x_2x_4, x_2, x_3, x_4\}.
$
Thus,
\[
\left\{\begin{array}{rcl}
\dot \eta_1 &=& -\eta_1+\eta_3-\eta_2\xi^2+u_\rme,\cr
\dot \eta_2 &=& \eta_2\xi,\cr
\dot \eta_3 &=& -\eta_2\xi-\eta_2\xi^2+u_\rme,\cr
\dot \xi &=& v_\rmd,\cr
y_\rme &=& \eta_2,\cr
y_\rmd &=& \xi,
\end{array}\right.
\]
with
$u_\rme=e^{-x_4}u_2$, and $v_\rmd=-x_4+u_1$.
We take the following further transformation on $\eta$. Let $z_\rma=\eta_1-\eta_3$,
$z_\rmb=\eta_2$ and $z_\rmc=\eta_2+\eta_3$.
%\[
%\left(\begin{array}{c}z_\rma\cr z_\rmb\cr z_\rmc\end{array}\right)
%   =\left[\begin{array}{rrr}1&0&-1\cr 0&1&0\cr 0&1&1\end{array}\right]\eta.
%\]
Then, the system takes the following form
\[
\left\{\begin{array}{rcl}
\dot z_\rma &=& -z_\rma+z_\rmb \xi,\cr
\dot z_\rmb &=& z_\rmb \xi,\cr
\dot z_\rmc &=& -z_\rmb \xi^2 +u_\rme,\cr
\dot \xi &=& v_\rmd,\cr
y_\rme &=& z_\rmb,\cr
y_\rmd &=& \xi,
\end{array}\right.
\]
with $\col\{z_\rma,z_\rmb,z_\rmc,\xi\}=
\col\{x_1-x_3-x_2x_4,x_2,x_2+x_3,x_4\}$.
The zero dynamics is
$
\dot z_\rma=-z_\rma.
$
It is also clear from the normal form above that
the system has an infinite zero of order $1$ and is not invertible.
\end{example}

\begin{example}\rm
\label{exx3}
Consider system (\ref{nonsys_nfn}) with
\[
f(x)=\pmatrix{x_3\cr x_4\cr x_3x_4\cr x_1x_3x_4},\quad
g(x)=\mm{1 &x_1\cr x_1&x_2\cr x_2&-x_3\cr x_3&1},\quad
h(x)=\pmatrix{x_1\cr x_2}.
\]
It is obvious that $f(0)=0$ and $h(0)=0$. We carry out the zero output structure algorithm as follows.

{\bf Initial Step.}
Let $\Theta_0(x)=h(x)$, $\Omega_0(x)=\emptyset$, $\rho_0=0$ and $k=1$.

{\bf Step 1.}
Let $M_1=\{x:\; x_1=x_2=0\}$,
\[
L_f\Theta_0(x)=\pmatrix{x_3\cr x_4},\quad
L_g\Theta_0(x)=\mm{1 &x_1\cr x_1&x_2}.
\]
Hence, $\rho_1=1$. Let $R_1=\mm{1 &0}$, $S_1=\mm{0&1}$. Thus,
\[
v_1=x_3+u_1+x_1u_2,\quad
\Omega_1(x)=x_1,\quad
P_{1,1}(x)=x_1,
\]
\[
W_1(x)=\mm{0& x_2-x_1^2},\quad
\Theta_1(x)=x_4-x_1x_3.
\]

{\bf Step 2.}
Let $M_2=\{x:\, x_1=x_2=x_4=0\}$.
\[
L_f\Theta_1(x)=-x_3^2,\quad
L_g\Theta_1(x)=\mm{-x_1x_2 & 1}.
\]
Hence, $\rho_2=2$. Let $R_2=1$, $S_2=\emptyset$. Thus,
\[
v_2=-x_3^2-x_1x_2u_1+u_2,\quad
\Omega_2(x)=\col\{x_1, x_4-x_1x_3\}.
\]

{\bf Final Step.} $k^\star=2$, $m_\rmd=2$, $n_\rmd=3$,
$\rho=\{1,2\}$ and $q=\{1,2\}.$

The distribution spanned by the column vectors of $g(x)$ is not involutive.
Define
\[
\col\{\eta, \xi_{1,1}, \xi_{2,1}, \xi_{2,2}\}
=\col\{x_3, x_1, x_2, x_4-x_1x_3\},
\]
\[
v_{\rmd,1}=x_3+u_1+x_1u_2
\]
and
\[
v_{\rmd,2}=-x_3^2-x_1x_2u_1+u_2.
\]
%hence, $x=\col\{\xi_{1,1}, \xi_{2,1}, \eta, \xi_{2,2}+\xi_{1,1}\eta\}.$
In the region $\{x: x_1^2+x_2^2<1\}$,
\begin{eqnarray*}
\dot{\eta}&=&x_3x_4+\mm{x_2 &-x_3}u\cr
&=&x_3x_4+\frac{-x_2x_3-x_3^3}{1+x_1^2x_2}+\frac{1}{1+x_1^2x_2}
\mm{x_2-x_1x_2x_3& -x_1x_2-x_3}
\pmatrix{v_{\rmd,1}\cr v_{\rmd,2}},
\end{eqnarray*}
and thus
the form (\ref{form2-nested}) is given by
\[
\left\{\begin{array}{rcl}
\dot{\eta}&=&
(\xi_{2,2}+\xi_{1,1}\eta)\eta+\frac{-\xi_{2,1}\eta-\eta^3}
{1+\xi_{1,1}^2\xi_{2,1}}\\
&&+\frac{1}{1+\xi_{1,1}^2\xi_{2,1}}\mm{\xi_{2,1}-\xi_{1,1}\xi_{2,1}\eta
& -\xi_{1,1}\xi_{2,1}-\eta}
\pmatrix{v_{\rmd,1}\cr v_{\rmd,2}},\cr
\dot{\xi}_{1,1}&=&v_{\rmd,1},\cr
\dot{\xi}_{2,1}&=&\xi_{2,2}+\xi_{1,1}v_{\rmd,1}-\mm{0& \xi_{2,1}-\xi_{1,1}^2}u,\cr
\dot{\xi}_{2,2}&=&v_{\rmd,2},\\
\rule{0mm}{5mm}
y_1&=&\xi_{1,1},\cr
y_2&=&\xi_{2,1}.
\end{array}\right.
\]
By letting $y=0$, we obtain the zero dynamics
\[
\dot{\eta}=-\eta^3.
\]
\end{example}

\begin{example}\rm
\label{appp}
Consider the dynamics of an underactuated vehicle \cite{foga94}.
%which has
%more degrees of freedom to be
%controlled than independent control inputs.
The model is given by
\[
M\dot{\nu}+C(\nu)\nu+D(\nu)\nu+g(\mu)=\mm{\tau\cr 0},\quad \dot{\mu}=J(\mu)\nu,
\quad y=\nu,
\]
where $\mu\in\RR^{n_1}$ denotes the position and orientation,
$\nu\in\RR^p$ denotes velocities to be controlled and $\tau\in\RR^m$
denotes control forces and moments
with $p>m$, $g(\mu)$ is the gravitation and buoyancy vector, and
the inertia matrix $M$ is constant, symmetric, nonsingular and positive definite.
Denote
$\col\{\nu_1, \nu_2\}=M\nu$, with $\nu_1\in\RR^m$.
The model is given by
\be\label{vehic1e}
\left\{\begin{array}{rcl}
\dot{\mu}&=&J(\mu)\nu,\\
\dot{\nu}_2&=&-\mm{0& I_{p-m}}C(\nu)\nu
    -\mm{0& I_{p-m}}D(\nu)\nu-\mm{0&I_{p-m}}g(\mu),\\
\dot{\nu}_1&=&-\mm{I_m& 0} C(\nu)\nu
    -\mm{I_m& 0} D(\nu)\nu-\mm{I_m& 0} g(\mu)
    +\tau, \\
y_\rme &=& \nu_2,\\
y_\rmd &=& \nu_1,
\end{array}\right.
\ee
where $\nu$ is evaluated by $M^{-1}\col\{\nu_1,\nu_2\}$.
Note that (\ref{vehic1e}) is already in the form of (\ref{form1}), with
$\eta=\col\{\mu,\nu_2\}$, $\xi=\nu_1$,
$v_\rmd=\tau-[I_m\; 0] C(\nu)\nu-[I_m\; 0] D(\nu)\nu-[I_m\; 0] g(\mu)$,
and $u_\rme$ nonexistent.
Thus, the system is left invertible and has $m$ infinite zeros of order $1$.
To determine its zero dynamics, we consider the subsystem (\ref{abc}) as follows,
\[
\left\{\begin{array}{rcl}
\dot{\mu}&=&0,\cr
\dot{\nu}_2&=&-[0\;\; I_{p-m}] g(\mu),\cr
y_\rme &=& \nu_2.
\end{array}
\right.
\]
The zero dynamics depends heavily on $[0\; I_{p-m}] g(\mu)$.
If $\mu$ is observed through $\nu_2$, then, there is no zero dynamics.
Otherwise, if $g(\mu)=0$, then zero dynamics is given by $\dot{\mu}=0$.
It is interesting to note that it has been shown in \cite{pees96} that $g(\mu)$ is important
for the stabilizability of underactuated vehicles.
This difficulty can also be seen in the form (\ref{vehic1e}).
In the absence of $g(\mu)$, neither $y_\rme$ nor $y_\rmd$ contains any information of
the state $\mu$.
%The normal forms with infinite zeros and zero dynamic also
%show such a property from another point of view.

\end{example}

\section{Summary of the Chapter}
\label{conclusion_nfn}

We have presented constructive algorithms for decomposing an
affine nonlinear system into its normal form representations.
Such algorithms generalize the existing results in several ways.
They require less restrictive assumptions on the system
and apply to general MIMO systems that do not necessarily
have the same number of inputs and outputs.
The resulting normal forms reveal various nonlinear extensions of
linear system structural properties.
These algorithms and the resulting normal forms
are thus expected to facilitate the solution of several nonlinear
control problems.
% Initial results have shown that these normal
% forms simplify the conventional backstepping design
% and motivate a new backstepping design procedure that is able to
% stabilize systems on which the conventional backstepping is not
% applicable
% \cite{lisn08}.

%\bigskip

%\appendices
%\section{Appendices}
\section{Proofs}

\subsection{Proof of Lemma~\ref{independent}}
\label{p_independent}

Denote
\[
ad_f\, g=\mm{ad_f g_1& ad_f g_2& \cdots & ad_f g_m},\;
\]
and
\[
\langle d\Theta, g\rangle =\mm{\langle d\Theta, g_1\rangle & \langle d\Theta, g_2\rangle & \cdots & \langle d\Theta, g_m\rangle}.
\]
We want to show that all row vectors in the following list are linearly independent:

\[
\left.\begin{array}{rlllll}
d\zeta_1:\quad & d R_1 \Theta_0  &&&\cr
d\zeta_2:\quad &  d R_2  S_{1}\Theta_0  & d R_2 \Theta_1  &\qquad\quad&\cr
\vdots\qquad &\quad\vdots&\qquad\vdots&\quad\ddots &&\cr
d\zeta_k:\quad & d R_k  S_{k-1 \leftrightarrow 1}\Theta_0  & d R_k  S_{k-1 \leftrightarrow 2}\Theta_1 & \quad\cdots & d R_k \Theta_{k-1} &\cr
d\zeta_{k+1}:\quad & d R_{k+1} S_{k\leftrightarrow 1}\Theta_0  & d R_{k+1} S_{k \leftrightarrow 2}\Theta_1 &\quad\cdots &
d R_{k+1} S_{k}\Theta_{k-1} &
d R_{k+1}\Theta_k.
\end{array}\right.
\]

The rows of $dR_1\Theta_0$ are linearly independent since the matrix $L_g R_1\Theta_0$ is of full row rank.
Next, we show that the rows of $dR_1\Theta_0$, $d R_2  S_{1}\Theta_0$ and $d R_2 \Theta_1$ are linearly independent.
%the case with $m=2$, $\rho_1=1$ and $\rho_2=2$.
To do this, consider
\be\label{adfg1}
\mm{d R_1 \Theta_0\cr
d R_2  S_{1}\Theta_0\cr  d R_2 \Theta_1}
\mm{g & ad_f g}
=
\mm{L_g R_1 \Theta_0 & \langle d R_1 \Theta_0,\, ad_f g\rangle \cr
L_g R_2  S_{1}\Theta_0 & \langle d R_2  S_{1}\Theta_0,\, ad_f g\rangle\cr
L_g R_2 \Theta_1 & \star}.
\ee
%is of full row rank for $x\in U$.
By row operation, the right hand side of (\ref{adfg1}) can be transformed to
\[
\mm{L_g R_1 \Theta_0 & \langle d R_1 \Theta_0,\, ad_f g\rangle \cr
0 & \langle d R_2  S_{1}\Theta_0,\, ad_f g\rangle - R_2P_{1,1} \langle d R_1 \Theta_0,\, ad_f g\rangle \cr
L_g R_2 \Theta_1 & \star}.
\]
Considering
$
\langle d\phi, ad_f g\rangle=L_f\langle d\phi, g\rangle -\langle dL_f\phi, g\rangle
$,
we have
\[
\langle d R_2  S_{1}\Theta_0,\, ad_f g\rangle - R_2P_{1,1} \langle d R_1 \Theta_0,\, ad_f g\rangle \hspace{70mm}
\]
\begin{eqnarray*}
&=&L_f\langle d R_2  S_{1}\Theta_0,\,  g\rangle-
\langle dL_f R_2  S_{1}\Theta_0,\, g\rangle
-R_2P_{1,1} L_f \langle d R_1 \Theta_0,\,  g\rangle\\
&& \hspace{70mm}+ R_2P_{1,1} \langle dL_f R_1 \Theta_0,\, g\rangle\\
&=&L_g R_2 \Theta_1.
\end{eqnarray*}
Therefore, (\ref{adfg1}) is of full row rank for $x\in U$. Hence the row vectors
$d R_1 \Theta_0, d R_2  S_{1}\Theta_0$, and $d R_2 \Theta_1$ are linearly independent.

Similarly, the row vectors of $\col\{d\zeta_1(x),d\zeta_2(x)\,\cdots,d\zeta_{k+1}(x)\}$
are linearly independent.

\subsection{Proof of Theorem~\ref{theo-form2}}
\label{p_theo-form2}

%We follow the proof method in \cite{isnc95}.
Note that
\be\label{dbarH}
d \Omega_{k^\star}(x) g_\rmd(x)=\Gamma_{\rmi\rmd}(x)\Gamma_\rmi(x)^{-1}\mm{0\cr I_{m_\rmd}}
=I_{m_\rmd}.
\ee
Thus, the column vectors of $g_\rmd(x)$ are linearly independent for $x\in U$.
By Frobenius' Theorem, there exists
$n-m_\rmd$ real-valued functions $\lambda_1(x),\lambda_2(x),\cdots,\newline \lambda_{n-m_\rmd}$ such that
the rows of $d\lambda_1(x),d\lambda_2(x),\cdots,d\lambda_{n-m_\rmd}$
are linearly independent and
\be\label{Delta_g_d}
\col\{d\lambda_1(x),d\lambda_2(x),\cdots,d\lambda_{n-m_\rmd}\}\, g_\rmd(x)=0.
\ee
Thus, $g_\rmd(x)$ spans the kernel space of $\col\{d\lambda_1(x),d\lambda_2(x),\cdots,d\lambda_{n-m_\rmd}\}$.
Suppose that $\nu(x): U\rightarrow \RR^n$ satisfies
\[
\mm{d\Phi_\rmd(x)\cr
\col\{d\lambda_1(x),d\lambda_2(x),\cdots,d\lambda_{n-m_\rmd}\}
}\nu(x)=0.
\]
Considering (\ref{Delta_g_d}) and $\col\{d\lambda_1(x),d\lambda_2(x),\cdots,d\lambda_{n-m_\rmd}\}\nu(x)=0$,
we have $\nu(x)=g_\rmd(x)\varpi(x)$, where $\varpi(x): U\rightarrow \RR^{m_\rmd}$.
Thus, $d\Phi_\rmd(x) g_\rmd(x)\varpi(x)=0$.
In view of (\ref{dbarH}) and the fact that $d\Omega_{k^\star}(x)$ is formed from some rows of $d\Phi_\rmd(x)$,
$d\Phi_\rmd(x) g_\rmd(x)$ has
full column rank, implying that $\varpi(x)=0$ and hence $\nu(x)=0$.
Therefore, the space spanned by the row vectors of $d\Phi_\rmd(x)$ and $\col\{d\lambda_1(x),d\lambda_2(x),\cdots,d\lambda_{n-m_\rmd}\}$ has
dimension $n$.
Selecting $n-n_\rmd$ elements from $\lambda_1(x),\lambda_2(x),\cdots,\lambda_{n-m_\rmd}$ to form $\Phi_\rme(x)$ in (\ref{trans1a}), and by (\ref{Delta_g_d}),
we have
$
L_{g_\rmd}\Phi_\rme(x)=0
$
for $x\in U$.
Consequently,
$
\dot \eta=d \Phi_\rme(x) f(x)+d \Phi_\rme(x) g(x)u
$
with
\[
d \Phi_\rme(x) g(x)u
=d \Phi_\rme(x) g(x)\Gamma_\rmi^{-1}(x)\mm{u_\rme\cr u_\rmd}
=d \Phi_\rme(x) g(x)\Gamma_\rmi^{-1}(x)\mm{I_{m-m_\rmd} \cr 0}u_\rme,
\]
which leads to (\ref{key3}).

%\section{Proof of Theorem~\ref{linear}}
%\label{linear-form1}

\subsection{Proof of Lemma~\ref{invariant}}
\label{p_invariant}

1) and 2) are obvious from the infinite zero structure algorithm.

3) Apply the infinite zero structure algorithm to the closed-loop system,
\[
\left\{\begin{array}{rcl}\dot x &=& f(x)+g(x)K(x)+g(x)\check{u},\\
                           y &=& h(x).\end{array}\right.
\]
Let $\check{\Theta}_0(x)=h(x)=\Theta_0(x)$. Then,
\[
\frac{d}{dt}\check{\Theta}_0(x)
=L_{(f+gK)} {\Theta}_0(x) + L_g {\Theta}_0(x) u
=L_f \Theta_0(x)+[L_g \Theta_0(x)]K(x)+L_g \Theta_0(x)u.
\]
Letting $\check{R}_1=R_1$ and $\check{S}_1=S_1$, we have
$
[-P_{1,1}(x)\check{R}_1+\check{S}_1]L_g\check{\Theta}_0(x)=0
$.
%with $\check{P}_{1,1}(x)={P}_{1,1}(x)$,
By (\ref{Pkj}) and (\ref{Thetak}),
\[
\check{\Theta}_1(x)=
[-P_{1,1}(x)\check{R}_1+\check{S}_1]L_{(f+gK)}\check{\Theta}_0(x)
=\Theta_1(x).
\]

Similarly, letting $\check{R}_k=R_k$ and $\check{S}_k=S_k$,
we obtain
$
\check{\Theta}_k(x)=\Theta_k(x).
$
Thus, $\check\rho_k=\rho_k$.

4)
Let $\check{\Theta}_0(x)=h(x)=\Theta_0(x)$. Then,
\begin{eqnarray*}
\frac{d}{dt}\check{\Theta}_0(x)
&=&L_{(f+Fh)} {\Theta}_0(x) + L_g {\Theta}_0(x) u\\
&=&L_f \Theta_0(x)+[d\Theta_0(x)]F(x)h(x)+L_g \Theta_0(x)u.
\end{eqnarray*}
Let $\check{R}_1=R_1$ and $\check{S}_1=S_1$. We have
\[
[-P_{1,1}(x)\check{R}_1+\check{S}_1]L_g\check{\Theta}_0(x)=0.
\]
%with $\check{P}_{1,1}(x)={P}_{1,1}(x)$,
By (\ref{Pkj}) and (\ref{Thetak}),
\[
\check{\Theta}_1(x)=
[-P_{1,1}(x)\check{R}_1+\check{S}_1]L_{(f+Fh)}\check{\Theta}_0(x)
=\Theta_1(x).
\]

Similarly, letting $\check{R}_k=R_k$ and $\check{S}_k=S_k$,
we obtain
$
\check{\Theta}_k(x)=\Theta_k(x),
$
Thus, $\check\rho_k=\rho_k$.

\subsection{Proof of Lemma~\ref{preinv}}
\label{p_preinv}

We first establish the following result.

\begin{lemma}
Let $G_i(x)=\col\{\Theta_0(x),\Theta_1(x),\cdots,\Theta_{i-1}(x)\}$, $i=1,2,\cdots,k^\star$.
We have
\goodbreak
\begin{eqnarray}
C_i(x)L_g G_i(x)&=&\col\{W_1(x),W_2(x),\cdots,W_i(x)\},\label{pre1}\\
C_i(x)L_f G_i(x)&=&\col\{\Theta_1(x),\Theta_2(x),\cdots,\Theta_i(x)\}\label{pre3},
\end{eqnarray}
where $C_1(x)=-P_{1,1}(x)R_1+S_1$,
\[
C_i(x)=\mm{C_{i-1}(x)&0\cr
-P_{i,\triangleleft}(x)\blkdiag\{R_1,R_2,\cdots,R_{i-1}\} & -P_{i,i}(x)R_i+S_i},
\]
and
\[
P_{i,\triangleleft}(x)=\mm{P_{i,1}(x) & P_{i,2}(x) & \cdots & P_{i,i-1}(x)}.
\]
Moreover, the rows of $C_i(x)$ form a basis of the solution space of the homogeneous linear equation
$\gamma L_gG_i(x)=0$ in $(M_i\cap O_i)^c$.
\end{lemma}
{\bf Proof:}
%We have $L_g\check{G}_1(x)=L_g G_1(x)=L_gh(x)$,
We carry out the proof by induction.
By Assumption $\calA_1$,
the matrix $L_g\Theta_0(x)$ has a constant rank $\rho_1$ in $(M_1\cap O_1)^c$.
By (\ref{Pkj2}), we have $C_1(x)L_g\Theta_0(x)=W_1(x)$.
Since $\col\{R_1,S_1\}$ is nonsingular,
$C_1(x)$ has full row rank $p-\rho_1$, and hence its rows form a basis of
the solution space of $\gamma L_g\Theta_0(x)=0$ in $(M_1\cap O_1)^c$.
By (\ref{Thetak2}), we have $\Theta_1(x)=C_1(x) L_f\Theta_0(x)$.

%Assume that for $i=j-1$, we have
Assume that
\[
C_{j-1}(x)L_g G_{j-1}(x)=\col\{W_1(x),W_2(x),\cdots,W_{j-1}(x)\},
\]
\[
C_{j-1}(x)L_f G_{j-1}(x)=\col\{\Theta_1(x),\Theta_2(x),\cdots,\Theta_{j-1}(x)\},
\]
and the rows of $C_{j-1}(x)$ form a basis of the solution space of
$\gamma L_gG_{j-1}(x)=0$ in $(M_{j-1}\cap O_{j-1})^c$.
By (\ref{Pkj2}) and (\ref{Thetak2}), we have
\[
\mm{-P_{j,\triangleleft}(x) & -P_{j,j}(x)R_{j}+S_{j}}
\mm{L_g\Omega_{j-1}(x)\cr L_g\Theta_{j-1}(x)}=W_j(x),
\]
\[
\mm{-P_{j,\triangleleft}(x) & -P_{j,j}(x)R_{j}+S_{j}}
\mm{L_f\Omega_{j-1}(x)\cr L_f\Theta_{j-1}(x)}=\Theta_{j}(x).
\]
Hence,
\[
\mm{-P_{j,\triangleleft}(x)\blkdiag\{R_1,R_2,\cdots,R_{j-1}\} & -P_{j,j}(x)R_{j}+S_{j}}\hspace{30mm}
\]
\[
\hspace{70mm}\mm{L_g G_{j-1}(x)\cr L_g\Theta_{j-1}(x)}=W_j(x),
\]
\[
\mm{-P_{j,\triangleleft}(x)\blkdiag\{R_1,R_2,\cdots,R_{j-1}\} & -P_{j,j}(x)R_{j}+S_{j}}\hspace{30mm}
\]
\[
\hspace{70mm}\mm{L_f G_{j-1}(x)\cr L_f\Theta_{j-1}(x)}=\Theta_{j}(x).
\]
Thus,
\[
C_j(x)L_g G_{j}(x)=W_j(x),
\]
\[
C_j(x)L_f G_{j}(x)=\col\{\Theta_1(x),\Theta_2(x),\cdots,\Theta_{j}(x)\}.
\]
The matrix $C_j(x)$ is of full row rank, since $-P_{j,j}(x)R_{j}+S_{j}$ is of full row rank.
The matrices $C_j(x)$ and $L_gG_{j}(x)$ have $\sum_{\ell=1}^j(p-\rho_\ell)$
and $p+\sum_{\ell=1}^{j-1}(p-\rho_\ell)$
rows, respectively, and the rank of $L_gG_{j}(x)$ is $\rho_j$.
Thus,
the rows of $C_j(x)$ form a basis of the solution space of
$\gamma L_gG_{j}(x)=0$ in $(M_j\cap O_j)^c$. \squ

Now we are ready to prove Lemma~\ref{preinv}.
We do it by induction.
Consider $i=1$.
According to the algorithm, $\check{\Theta}_0(x)=h(x)=\Theta_0(x)$,
thus, $\check{M}_1=M_1$.
The rows of $C_1(x)$ form a basis of
the solution space of the homogeneous linear equation $\gamma L_g\Theta_0(x)=0$ in $(\check{M}_1\cap\check{O}_1)^c$.
Similarly, the rows of $\check{C}_1(x)=-\check{P}_{1,1}(x)\check{R}_1+\check{S}_1$
span the same solution space in $(\check{M}_1\cap\check{O}_1)^c$.
%It is obvious that $N_1\cap\check{M}_1=\{x : \mbox { a neighborhood of } x^\circ  \mbox{ with } h(x)=0\}$.
Therefore,
\[
\check{C}_1(x)
=T_1(x)C_1(x)+\varpi_1(x),
\]
where the matrix $T_1(x): (O_k\cap\check{O}_k)\rightarrow \RR^{(p-\rho_1)\times (p-\rho_1)}$ is a nonsingular and smooth,
and $\varpi_1(x)$ is smooth with $\varpi_1(x)=0$ in $\check{M}_1\cap O_1 \cap \check{O}_1$.
Thus, by (\ref{Thetak2}),
\[
\check{\Theta}_1(x)=\check{C}_1(x) L_f\Theta_0(x)
=T_1(x)\Theta_1(x)+V_1(x),
\]
where $V_1(x)=\varpi_1(x) L_f\Theta_0(x)=0$ in $M_1\cap O_1 \cap \check{O}_1$.

Assume that,
for $i=1,2,\cdots,j-1$, %(\ref{pkinv1}) and
equations in (\ref{pkinv2}) are satisfied. That is,
\[
\check{M}_i=M_i,\quad \check{\Theta}_i(x)=\sum_{l=1}^{i-1}Q_{i,l}(x)\Theta_l(x)+T_i(x)\Theta_i(x)+V_i(x),
\quad i=1,2,\cdots,j-1,
\]
where $T_i(x)$ is nonsingular and $V_i(x)$ is smooth with $V_i(x)=0$ in $M_i\cap O_i \cap \check{O}_i$.
Thus,
\be\label{pkinv3}
\check{G}_j(x)
=E_j(x)G_j(x)+\col\{0,V_1(x),\cdots,V_{j-1}(x)\},
\ee
with
$\col\{0,V_1(x),\cdots,V_{j-1}(x)\}=0$ in
$M_{j-1}\cap O_{j-1} \cap \check{O}_{j-1}$, and $E_j(x)$ being nonsingular, where
\[
E_j(x)=
\mm{I_p&0&0&\cdots&0\cr
0&T_1(x)&0&\cdots&0\cr
0& Q_{2,1}(x)&T_2(x)&\cdots&0\cr
\vdots&\vdots&\vdots&\ddots&\vdots\cr
0&Q_{j-1,1}(x)&Q_{j-1,2}(x)&\cdots& T_{j-1}(x)}.
\]

By (\ref{pkinv3}), we know that
$G_j(x)=0$ is equivalent to $\check{G}_j(x)=0$. Thus,
$\check{M}_j=M_j$.
We also have
\be\label{pkinv4}
d\check{G}_j(x)
=E_j(x) d G_j(x)+D_j(x),
\ee
where
\[
D_j(x)=\sum_{l=1}^\varsigma G_{j,l}\frac{\partial E_{j,l}}{\partial x}+
\frac{\partial}{\partial x} \col\{0,V_1(x),\cdots,V_{j-1}(x)\},
\]
with
$
\col\{G_{j,1},G_{j,2},\cdots,G_{j,\varsigma}\}=G_j(x),\;
[E_{j,1},E_{j,2},\cdots,E_{j,\varsigma}]=E_j(x),
$
and
$
\varsigma=\sum_{\ell=1}^{j-1}(p-\rho_\ell).
$
It is obvious that $D_j(x)=0$ in $M_j\cap O_j \cap \check{O}_j$.

By (\ref{pre1}),
$C_j(x)L_gG_j(x)=0$ and $\check{C}_j(x)L_g\check{G}_j(x)=0$ in $M_j\cap O_j \cap \check{O}_j$.
The rows of $C_j(x)$ and $\check{C}_j(x)$
span the solution spaces of homogeneous linear equations
$\gamma L_gG_j(x)=0$ and
$\gamma L_g\check{G}_j(x)=0$ in $M_j\cap O_j \cap \check{O}_j$, respectively.
By (\ref{pkinv4}),
\[
L_g\check{G}_j(x)
=E_j(x) L_g G_j(x)+D_j(x)g(x),
\]
and thus,
\be\label{pkinv9}
\check{C}_j(x)E_j(x)=
F_j(x){C}_j(x)+\varpi_{j+1}(x),
\ee
where $F_j(x)$ is a nonsingular matrix valued smooth function,
and $\varpi_j(x)=0$ in $M_j\cap O_j \cap \check{O}_j$.
Denote
\[
F_j(x)=\mm{Y_j(x)&\bar{Y}_j(x)\cr Q_j(x)& T_j(x)}, \quad
\varpi_{j+1}(x)=\mm{\bar{\mu}_j(x)\cr \mu_j(x)},
\]
where $T_j(x)$ is a $(p-\rho_j)\times(p-\rho_j)$ matrix, and
$\mu_j(x)$ is smooth with
$\mu_j(x)=0$ in $M_j\cap O_j \cap \check{O}_j$.
Due to the structure of
$C_j(x)$ and $\check{C}_j(x)$, we know that
$\bar{Y}_j(x)=0$ in $M_j\cap O_j \cap \check{O}_j$.
Thus, $T_j(x)$ is nonsingular in a neighborhood of $x=0$,
which contains $M_j\cap O_j \cap \check{O}_j$.
By (\ref{pkinv9}),
\begin{eqnarray}
\label{pkinv11}
\mm{-\check{P}_{j,\triangleleft}(x)\blkdiag\{\check{R}_1,\check{R}_2,\cdots,\check{R}_{j-1}\}
& -\check{P}_{j,j}(x) \check{R}_j+\check{S}_j}E_j(x)\nonumber\\
=Q_j(x)\mm{C_j(x)&0}
+T_j(x) [ -P_{j,\triangleleft}(x)\blkdiag\{R_1,R_2,\cdots,R_{j-1}\} \nonumber\\
 -P_{j,j}(x) R_j+S_j ] +U_j(x).
\end{eqnarray}
And by (\ref{Thetak2}), we have
\[
\Theta_j(x)=
\mm{-P_{j,\triangleleft}(x)\blkdiag\{R_1,R_2,\cdots,R_{j-1}\} & -P_{j,j}(x) R_j+S_j}
L_f G_j(x),
\]
\[
\check{\Theta}_j(x)=
\mm{-\check{P}_{j,\triangleleft}(x)\blkdiag\{\check{R}_1,\check{R}_2,\cdots,\check{R}_{j-1}\}
& -\check{P}_{j,j}(x) \check{R}_j+\check{S}_j}
L_f \check{G}_j(x).
\]
Thus, multiplying (\ref{pkinv11}) to the right by $L_fG_j(x)$ and using (\ref{pre3}) and (\ref{pkinv4}),
we have
\begin{eqnarray*}
\check{\Theta}_j(x)&=&
Q_j(x) C_j(x)L_fG_{j-1}(x)+T_j(x)\Theta_j(x)+V_j(x)\\
&=&\sum_{l=1}^{j-1}Q_{j,l}(x) \Theta_l(x)+T_j(x)\Theta_j(x)+V_j(x),
\end{eqnarray*}
where $\mm{Q_{j,1}(x)&Q_{j,1}(x)&\cdots&Q_{j,j-1}(x)}=Q_j(x)$ and
\[
V_j(x)=\mu_j(x)L_fG_j+\hspace{80mm}
\]
\[
\hspace{10mm}
\mm{-\check{P}_{j,\triangleleft}(x)\blkdiag\{\check{R}_1,\check{R}_2,\cdots,\check{R}_{j-1}\}
& -\check{P}_{j,j}(x) \check{R}_j+\check{S}_j}D_j(x)f(x).
\]
Therefore, $V_j(x)=0$ in $M_j\cap O_j \cap \check{O}_j$.

\subsection{Proof of Lemma~\ref{pindependent2}}
\label{p_independent2}

By the infinite zero structure algorithm, we know that
$\Gamma_{\rmi\rmd}(x)=L_g\Omega_{k^\star}(x)$ is of full row rank.
Note that $R_i S_{i-1 \leftrightarrow 1}$, $i=1,2,\cdots,{k^\star}$, are the coefficients
in $\zeta_{i,1}$.
So if $d\bar{\Phi}_\rmd(0)$ is of full row rank, $\Gamma_{\rmo\rmd}$ is of full row rank.
Thus, we only need to prove that $d\bar{\Phi}_\rmd(0)$ is of full row rank.
We prove it by induction.

Recall that
$d\bar{\Phi}_\rmd(x)
=\col\{d\zeta_1,d\zeta_2,\cdots,d\zeta_{k^\star}\}$.
We first prove that the row vectors of $d\zeta_1(0)$, or
$d R_1\Theta_0(0)$, are linearly independent.
It follows directly from the fact that $L_g R_1\Theta_0(x)$ has full row rank.

%  \[
%  \begin{tabular}{|l|l|l|l|l|}
%  \hline
%  $w_{1,0}d R_1 \Theta_0$  &&&&\\
%  $w_{2,0}d R_2  S_{1}\Theta_0$  & $w_{2,1}d R_2 \Theta_1$  & $\qquad\quad$&&\\
%  $\quad\vdots$ & $\qquad\vdots$ & $\quad\ddots$ &&\\
%  $w_{k,0}d R_k  S_{k-1 \leftrightarrow 1}\Theta_0$  & $w_{k,1}d R_k  S_{k-1 \leftrightarrow 2}\Theta_1$
%  & $\quad\cdots$ & $w_{k,k-1}d R_k \Theta_{k-1}$ &\\
%  $w_{k+1,0}d R_{k+1} S_{k\leftrightarrow 1}\Theta_0$  & $w_{k+1,1}d R_{k+1} S_{k \leftrightarrow 2}\Theta_1$ &$\quad\cdots$ &
%  $w_{k+1,k-1}d R_{k+1} S_{k}\Theta_{k-1}$ &
%  $w_{k+1,k}d R_{k+1}\Theta_k$\\
%  \hline
%  \end{tabular}
%  \]

Assume that the rows of
$\col\{d\zeta_1(0),d\zeta_2(0)\,\cdots,d\zeta_k(0)\}$
are linearly independent. We want to prove that
the rows of $\col\{d\zeta_1(0),d\zeta_2(0)\,\cdots,d\zeta_{k+1}(0)\}$
are linearly independent.

Let
$w_{i,j}(x): U\rightarrow \RR^{1\times(\rho_i-\rho_{i-1})}$,  $\;j=0,1,\cdots,i-1$, and $\; i=1,2,\cdots,k+1$.
Define
\be
\beta(x)
=\sum_{l=1}^{k+1} w_{l,l-1}(x) d R_l\Theta_{l-1}(x)+
\sum_{j=1}^{k}\;\sum_{i=j+1}^{k+1} w_{i,j-1}(x) d R_i  S_{i-1 \leftrightarrow j}\, \Theta_{j-1}(x).
\ee
By (\ref{Pkj2}),
\begin{eqnarray*}
\beta(x)g(x)
&=&
\sum_{l=1}^{k+1} w_{l,l-1}(x) L_g R_l\Theta_{l-1}(x)\\
&& \hspace{15mm}+
\sum_{j=1}^{k}\;\sum_{i=j+1}^{k+1} w_{i,j-1}(x) R_i  S_{i-1 \leftrightarrow j+1} L_g S_j \Theta_{j-1}(x)\cr
&=&
\sum_{l=1}^{k+1} w_{l,l-1}(x) L_g R_l\Theta_{l-1}(x)+\sum_{j=1}^{k}\;\sum_{i=j+1}^{k+1}\cr
&&\qquad w_{i,j-1}(x) R_i  S_{i-1 \leftrightarrow j+1}
\Bigl(\sum_{l=1}^j P_{j,l}(x)L_g R_l\Theta_{l-1}(x)+W_j(x)\Bigl)\cr
&=&
w_{k+1,k}(x) L_g R_{k+1}\Theta_k(x)+ \sum_{l=1}^k \psi_l(x) L_g R_l\Theta_{l-1}(x)\cr
&&+\sum_{j=1}^{k}\;\sum_{i=j+1}^{k+1} w_{i,j-1}(x) R_i  S_{i-1 \leftrightarrow j+1}W_j(x),
\end{eqnarray*}
where
\[
\psi_l(x)=
w_{l,l-1}(x) +
\sum_{j=l}^{k}\;\sum_{i=j+1}^{k+1} w_{i,j-1}(x) R_i  S_{i-1 \leftrightarrow j+1} P_{j,l}(x),\quad l=1,2,\cdots,k.
\]

By (\ref{Thetak2}),
\begin{eqnarray}\label{betaf}
\beta(x)f(x)
&=&
\sum_{l=1}^{k+1} w_{l,l-1}(x) L_f R_l\Theta_{l-1}(x)+ \sum_{j=1}^{k}\;\sum_{i=j+1}^{k+1}\cr
&& \hspace{10mm} w_{i,j-1}(x) R_i  S_{i-1 \leftrightarrow j+1} L_f S_j \Theta_{j-1}(x)\cr
&=&
\sum_{l=1}^{k+1} w_{l,l-1}(x) L_f R_l\Theta_{l-1}(x)+\sum_{j=1}^{k}\;\sum_{i=j+1}^{k+1}\cr
&& \hspace{10mm} w_{i,j-1}(x) R_i  S_{i-1 \leftrightarrow j+1}\Bigl[\Theta_j(x)+ \sum_{l=1}^j P_{j,l}(x)L_f R_l\Theta_{l-1}(x)\Bigl]\cr
&=&
w_{k+1,k}(x) L_f R_{k+1}\Theta_k(x)+\sum_{l=1}^k \psi_l(x) L_f R_l\Theta_{l-1}(x)\cr
&&
\hspace{10mm}+\sum_{j=1}^{k}\;\sum_{i=j+1}^{k+1} w_{i,j-1}(x) R_i  S_{i-1 \leftrightarrow j+1}\Theta_j(x).
\end{eqnarray}
Let
\be\label{q4}
\beta(x)=0.
\ee
Thus,
$
\beta(x)g(x)=0.
$
Since the matrix \newline
$\col\{L_g R_1\Theta_0(x), L_g R_2\Theta_1(x), \cdots,$
$L_g R_{k+1}\Theta_k(x)\}$ is of full row rank in $(M_k\cap O_k)^c$, we have
\be\label{Wl}
w_{k+1,k}(x)=0,\quad \psi_l(x)=0, \quad l=1,2,\cdots, k,
\ee
in $(M_k\cap O_k)^c$.
Thus, by (\ref{q4}) and (\ref{Wl}),
\be\label{q5}
\sum_{j=1}^{k}\;\sum_{i=j+1}^{k+1} w_{i,j-1}(0) d R_i  S_{i-1 \leftrightarrow j}\Theta_{j-1}(0)=0.
\ee
By (\ref{betaf}),
\begin{eqnarray}
d\Bigl[\beta(x) f(x)-
\sum_{l=1}^{k+1} \psi_l(x) L_f R_l\Theta_{l-1}(x)\Bigl]\hspace{30mm}\nonumber\\
=
\sum_{j=1}^{k}\;\sum_{i=j+1}^{k+1} w_{i,j-1}(x) d R_i  S_{i-1 \leftrightarrow j+1}\Theta_j(x)
 +\sum_{l=1}^k\Theta_l(x)\Pi_l(x),
\end{eqnarray}
where $\Pi_l(x)$, $l=1,2,\cdots,k$, are matrix valued functions of $w_{i,j}(x)$. Consider $\Theta_i( 0 )=0$.
We have
\be\label{q3}
\sum_{j=1}^{k}\;\sum_{i=j+1}^{k+1} w_{i,j-1}(0) d R_i  S_{i-1 \leftrightarrow j+1}\Theta_j(0)
=0.
\ee
By (\ref{q5}) and (\ref{q3}), we have
$w_{i,j}(0)=0$, for $j=0,1,\cdots,i-1$ and $i=1,2,\cdots,k+1$.
In conclusion, the row vectors of $\col\{d\zeta_1(0),d\zeta_2(0)\,\cdots,d\zeta_{k+1}(0)\}$
are linearly independent.

\chapter{Backstepping Design Procedure}

In Chapter 3, we developed a structural decomposition for
multiple input multiple output
nonlinear systems that are affine in control but
otherwise general. In this chapter we exploit the properties of such
a decomposition for the purpose of solving the stabilization problem.
In particular, this %structural
decomposition simplifies
the conventional backstepping design and motivates a new backstepping design
procedure that is able to stabilize some systems on which the conventional
backstepping is not applicable.
An numerical example also shows that different backstepping procedure lead to
different control performance.

\section{Introduction and Problem Statement}

Consider a nonlinear system of the form
\be
\label{nonsys_rev}
\left\{\begin{array}{rcl}\dot x &=& f(x)+g(x)u,\cr
                           y &=& h(x),\end{array}\right.
\ee
where $x\in\RR^n$, $u\in\RR^m$ and $y\in\RR^p$ are the state, input and
output, respectively, and the mappings $f$, $g$ and $h$ are
smooth % in a neighborhood of a point $0$
with $f(0)=0$ and $h(0)=0$.

In Chapter 3, we study the
structural properties of affine-in-control nonlinear systems beyond
the case of square invertible systems. We propose an algorithm that
identifies a set of integers that are equivalent to the infinite
zero structure of linear systems and leads to
a normal form representation that corresponds to these integers as well
as to the system invertibility structure.
This new normal form representation takes the following form
\be
\label{form1a}
\!\left\{\begin{array}{rcl}
\dot \eta &=& f_\rme(\eta,\xi)+g_\rme
    (\eta,\xi)u_\rme,\cr
\dot\xi_{i,j} &=& \xi_{i,j+1}+\displaystyle\sum_{l=1}^{i-1}
            \delta_{i,j,l}(x)v_{\rmd,l},\quad j=1,2,\cdots,q_i-1,\cr
\dot\xi_{i,q_i} &=& v_{\rmd,i},\cr
y_\rme &=& h_\rme (\eta, \xi),\cr
y_{\rmd,i} &=& \xi_{i,1}, \quad i=1,2,\cdots,m_\rmd,
\end{array}\right.
\ee
where $q_1\leq q_2\leq\cdots\leq q_{m_\rmd}$,
$\xi_i=\{\xi_{i,1},\xi_{i,2},\cdots,\xi_{i,q_i}\}$, $i=1,2,\cdots,m_\rmd$,
$\xi=\{\xi_1,\xi_2,\cdots,\xi_{m_\rmd}\}$,
$v_{\rmd,i}=a_i(x)+b_i(x)u$,
with the matrix $\col\{b_1(x),b_2(x),\cdots,$ $b_{m_\rmd}(x)\}$ being
of full row rank and smooth, and
$
\delta_{i,j,l}(x)=0,\; \mbox{for } \; j<q_l, \; i=1,2,\cdots,m_\rmd.
$

We note here that $m_\rmd$ is the largest integer for which the
system assumes the above form.
The system is left invertible if $u_\rme$ is non-existent, right invertible
if $y_\rme$ is non-existent, and invertible if both are non-existent.
In the case that the system is square and invertible, {\em i.e.},
the system that was considered in \cite{isnc95,isnc99},
$m=p=m_\rmd$ and the parts containing $y_\rme$ and $u_\rme$
drop off. Thus, the normal form (\ref{form1a}) simplifies to
\be
\label{form1-square_stabilization}
\left\{\begin{array}{rcl}
\dot \eta &=& f_\rme(\eta,\xi),\cr
\dot\xi_{i,j} &=& \xi_{i,j+1}+\displaystyle\sum_{l=1}^{i-1}
            \delta_{i,j,l}(x)v_l,\quad j=1,2,\cdots,q_i-1,\cr
\dot\xi_{i,q_i} &=& v_i,\cr
 y_i &=& \xi_{i,1}, \quad i=1,2,\cdots,m,
\end{array}\right.
\ee
where $q_1\leq q_2\leq\cdots\leq q_m$, and
\be\label{key2_stabilization}
\delta_{i,j,l}(x)=0,\quad \mbox{for } \; j<q_l, \; i=1,2,\cdots,m.
\ee
We note that if $q_1=q_2=\cdots=q_m$
in (\ref{form1-square_stabilization}), then by the property (\ref{key2_stabilization}),
the system turns out to have uniform relative degrees $q_1$.
The $\dot\xi_{i,j}$ equation in (\ref{form1-square_stabilization}) displays a
triangular structure of the control inputs that enter the system.
The property (\ref{key2_stabilization}) imposes additional
structure within each chain of integrators
on how control inputs enter the system.
With this additional structural property,
the set of integers $\{q_1, q_2, \cdots, q_{m_\rmd}\}$ indeed represent
infinite zero structure when the system is linear.

Control design techniques and structural decompositions
of nonlinear systems
have been developed interweavingly.
The discovery of structural properties and the corresponding
normal form representation
of the system
motivates new control
designs. On the other hand, the desire for achieving more stringent
closed-loop performances for a larger class of systems entails the
exploitation of more intricate structural properties.
For example, various stabilization results have been
obtained in this process.
In this chapter, we would like to
revisit the problem of stabilization. We will show how the property (\ref{key2_stabilization})
simplifies the conventional backstepping
design and motivates a new backstepping design technique that is
able to stabilize some systems that cannot be stabilized by
the conventional backstepping technique.

\section{Review of the Backstepping Design Technique}
\label{preliminary}

In the section, we recall some results on the backstepping design
methodology \cite{isnc95,isnc99,krna95}.
We first recall the integrator backstepping, on which the recursive backstepping procedure
is develop.

\begin{lemma}\cite{isnc95}\label{base1}
Consider
\be\label{inte}
\left\{\begin{array}{rcl}\dot \eta &=& f(\eta,\xi),\cr
                          \dot \xi &=& u,\end{array}\right.
\ee
where $(\eta,\xi)\in\RR^n\times R$ and $f(0,0)=0$.
Suppose there exists a smooth real-valued function
$
\xi=v^\star(\eta),
$
with $v^\star(0)=0$, and a smooth real-valued function $V(\eta)$, which is positive
definite and proper, such that
\[
\frac{\partial V}{\partial \eta} f(\eta,v^\star(\eta))<0,\quad \forall \eta\not=0.
\]
Then, there exists a smooth static feedback law $u=u(\eta,\xi)$ with $u(0,0)=0$, and
a smooth real-valued function $W(\eta,\xi)$, which is positive definite and proper, such that
\be\label{negative}
\frac{\partial W}{\partial \eta} f(\eta,\xi) + \frac{\partial W}{\partial \xi} u(\eta,\xi)<0,
\quad \forall (\eta,\xi)\not=0.
\ee
That is, $u=u(\eta,\xi)$ globally asymptotically stabilizes (\ref{inte}) at its equilibrium $(\eta,\xi)=0$.
\end{lemma}

The negative definiteness property in (\ref{negative}) can be replaced with a negative semi-definiteness property
along with a LaSalle's Invariance argument.

The backstepping design
method is readily applicable to systems that have vector relative
degrees and are represented in the following form,
\[
\left\{\begin{array}{rcl}
\dot \eta &=& f_0(x)+g_0(x)u,\cr
\dot\xi_{i,j} &=& \xi_{i,j+1},\qquad j=1,2,\cdots,r_i-1,\cr
\dot\xi_{i,r_i} &=& v_i,\cr
y_i &=& \xi_{i,1}, \qquad i=1,2,\cdots,m,
\end{array}\right.
\]
which contains $m$ clean chains of integrators.
Each of these chains is independently controlled
by a separate input.
Let us consider the following %standard
assumption.
\begin{assumption}
\label{assumption_stand}
The dynamics $\eta$ is driven only by $\xi_{i,1}$, $i=1,2,\cdots,m$, {\em i.e.},
\be
\label{doteta1}
\dot\eta=f_0(\eta,\xi_{1,1},\xi_{2,1},\cdots,\xi_{m,1}),
\ee
and there exist smooth functions
$\phi_{i,1}(\eta)$,
with $\phi_{i,1}(0)=0$, $i=1,2,\cdots,m$, such that
$
%\label{doteta2}
\dot\eta=f_0(\eta,\phi_{1,1}(\eta),\phi_{2,1}(\eta),\cdots,\phi_{m,1}(\eta))
$
is globally asymptotically stable at its equilibrium $\eta=0$.
\end{assumption}

Suppose that Assumption~\ref{assumption_stand} holds, then for the systems with vector relative degree,
it is straightforward to
design a globally asymptotically stabilizing feedback law
recursively, by viewing the next integrators as a new virtual input.
Such a design procedure is thus referred to as ``backstepping."

The technique of backstepping, however, cannot as easily been
implemented if the system does not have a vector relative degree.
An additional assumption,
which requires the coefficient functions $\delta_{i,j,l}$
in the following normal form
\be\label{zerody2_stabilization}
\left\{\begin{array}{rcl}
\dot \eta &=& f_0(x),\cr
\dot\xi_{i,j} &=& \xi_{i,j+1}+
  \displaystyle          \sum_{l=1}^{i-1}
            \delta_{i,j,l}(x)v_l, \quad j=1,2,\cdots,n_i-1,\cr
\dot\xi_{i,n_i} &=& v_i,\cr
 y_i &=& \xi_{i,1}, \qquad i=1,2,\cdots,m.
\end{array}\right.
\ee
to display a certain ``triangular" dependency on the state variables,
is needed \cite{isnc95,isnc99}. In what follows, we recall
from \cite{isnc99} such an additional assumption %on the normal form (\ref{zerody2})
and the backstepping design procedure that is implemented
under these assumptions.

\begin{assumption}
\label{tri1}
The functions $\delta_{i,j,l}$ depend only on variable
$\xi_{\ell_\rmp,\ell_\rmb}$, with
\ben
\item
$1\leq\ell_\rmp\leq m$ and $\ell_\rmb=1$; or,
\item
$\ell_\rmp\leq i-1$; or,
\item
$\ell_\rmp=i$ and $\ell_\rmb\leq j$.
\een
\end{assumption}

Under Assumptions~\ref{assumption_stand} and \ref{tri1}, a feedback law
\[
v_i=v_i(\eta;
%\xi_{1,1},\xi_{2,1},\cdots,\xi_{m,1};
\xi_1,\xi_2,
\cdots,\xi_i),\; i=1,2,\cdots,m,
\]
that globally stabilizes the whole system can be constructed from
$\phi_{i,1}(\eta), \;i=1,2,\cdots,m$, through a backstepping procedure \cite{isnc99}.
The procedure commences with the subsystem (\ref{doteta1}),
and is followed by backstepping $n_1$ times
through the variables in first chain of integrators to obtain
\[
v_1=v_1(\eta;
\xi_1;
\phi_{2,1}(\eta),\phi_{3,1}(\eta),\cdots,\phi_{m,1}(\eta)),
\]
and backstepping $n_2$ times through the variables in the second
chain of integrators to obtain
the feedback law
\[
v_2=v_2(\eta;
\xi_1,\xi_2;\phi_{3,1}(\eta),\phi_{4,1}(\eta),
\cdots,\phi_{m,1}(\eta)).
\]
This procedure is continued chain by chain for $i=1$ through $m$,
each backstepping $n_i$ times through $i$-th chain of integrators
to discover the feedback law
\[
v_i=v_i(\eta;\xi_1,\xi_2,\cdots,\xi_i;
\phi_{i+1,1}(\eta),\phi_{i+2,1}(\eta),
\cdots,\phi_{m,1}(\eta)).
\]
As the backstepping is implemented on the integrators chain by chain,
we will refer to the above backstepping procedure as the chain-by-chain
backstepping.

% \begin{rema}
% The above conventional chain-by-chain backstepping design
% procedure our new normal form (\ref{form1-square_stabilization}) - (\ref{key2_stabilization}),
% and its implementation on this new normal form is simpler than on
% the earlier normal form (\ref{zerody2_stabilization}) that is without Property (\ref{key2_stabilization}).
% \end{rema}

\section{Backstepping Design Procedures Revisited}
\label{backstepping}

In this section we focus on systems that are square invertible and
discuss about their stabilization by the backstepping technique.
We will first show that the conventional chain-by-chain backstepping design
technique as described in \cite{isnc99} and recalled in Section
\ref{preliminary} is applicable to our new normal form (\ref{form1-square_stabilization})-(\ref{key2_stabilization}),
and its implementation on this new normal form is simpler than on
the earlier normal form (\ref{zerody2_stabilization}).
We the propose a new
backstepping procedure which we refer to as the level-by-level
backstepping. In the level-by-level backstepping design procedure,
the backstepping is first implemented on the first integrators
of all chains and then on the second integrators of all chains, and so on.
We will show that the level-by-level backstepping will allow
the backstepping to be implemented on some systems for which the
chain-by-chain backstepping procedure is not applicable.
We will also show that the chain-by-chain backstepping and the
level-by-level backstepping can be mixed and implemented
on a same system
to allow stabilization of a larger class of systems.

\subsection{Conventional Chain-by-Chain Backstepping}

Since the normal form (\ref{form1-square_stabilization})-(\ref{key2_stabilization}) is a special case
of the normal form (\ref{zerody2_stabilization}), backstepping is applicable to
it.  As explained in \cite{isnc99}, the chain-by-chain backstepping
requires the system (\ref{form1-square_stabilization}) to satisfy
Assumptions~\ref{assumption_stand} and \ref{tri1}.
Under these two assumptions, the normal form (\ref{form1-square_stabilization})-(\ref{key2_stabilization})
is much simpler than the normal form (\ref{zerody2_stabilization}).
This simpler form makes the
implementation of the chain-by-chain backstepping simpler.
%We will illustrate this simplification by an example.

\begin{example}
\rm
A three input three output system in the form (\ref{zerody2})
with three chains of integrators of lengths $\{2,4,4\}$ and satisfying Assumption~\ref{tri1} will take the form (see Fig.~\ref{fig244_old}),
\be\label{equ_delta1}
\left\{\begin{array}{rcll}
\dot\eta &=&f_0(\eta,\xi_{1,1},\xi_{2,1},\xi_{3,1}),\cr
\dot\xi_{1,1} &=& \xi_{1,2},\cr
\dot\xi_{1,2} &=& v_1,\cr
\dot\xi_{2,1} &=& \xi_{2,2}+\delta_{2,1,1}(\eta,\xi_1,\xi_{2,1},\xi_{3,1})v_1,\cr
\dot\xi_{2,2} &=& \xi_{2,3}+\delta_{2,2,1}(\eta,\xi_1,\xi_{2,1},\xi_{2,2},\xi_{3,1})v_1,\cr
\dot\xi_{2,3} &=& \xi_{2,4}+\delta_{2,3,1}(\eta,\xi_1,\xi_{2,1},\xi_{2,2},\xi_{2,3},\xi_{3,1})v_1,\cr
\dot\xi_{2,4} &=& v_2,\cr
\dot\xi_{3,1} &=& \xi_{3,2}+\delta_{3,1,1}(\eta,\xi_1,\xi_2,\xi_{3,1})v_1+
\delta_{3,1,2}(\eta,\xi_1,\xi_2,\xi_{3,1})v_2,\cr
\dot\xi_{3,2} &=& \xi_{3,3}+\delta_{3,2,1}(\eta,\xi_1,\xi_2,
     \xi_{3,1},\xi_{3,2})v_1+\delta_{3,2,2}(\eta,\xi_1,\xi_2,\xi_{3,1},\xi_{3,2})v_2,\cr
\dot\xi_{3,3} &=& \xi_{3,4}+\delta_{3,3,1}(\eta,\xi_1,\xi_2,\xi_{3,1},
    \xi_{3,2},\xi_{3,3})v_1 \cr
    && \hspace{30mm}+\delta_{3,3,2}(\eta,\xi_1,\xi_2,\xi_{3,1},\xi_{3,2},\xi_{3,3})v_2,\cr
\dot\xi_{3,4} &=& v_3.
\end{array}\right.
\ee
\begin{figure}[ht]
\centerline{\includegraphics[width=115mm]{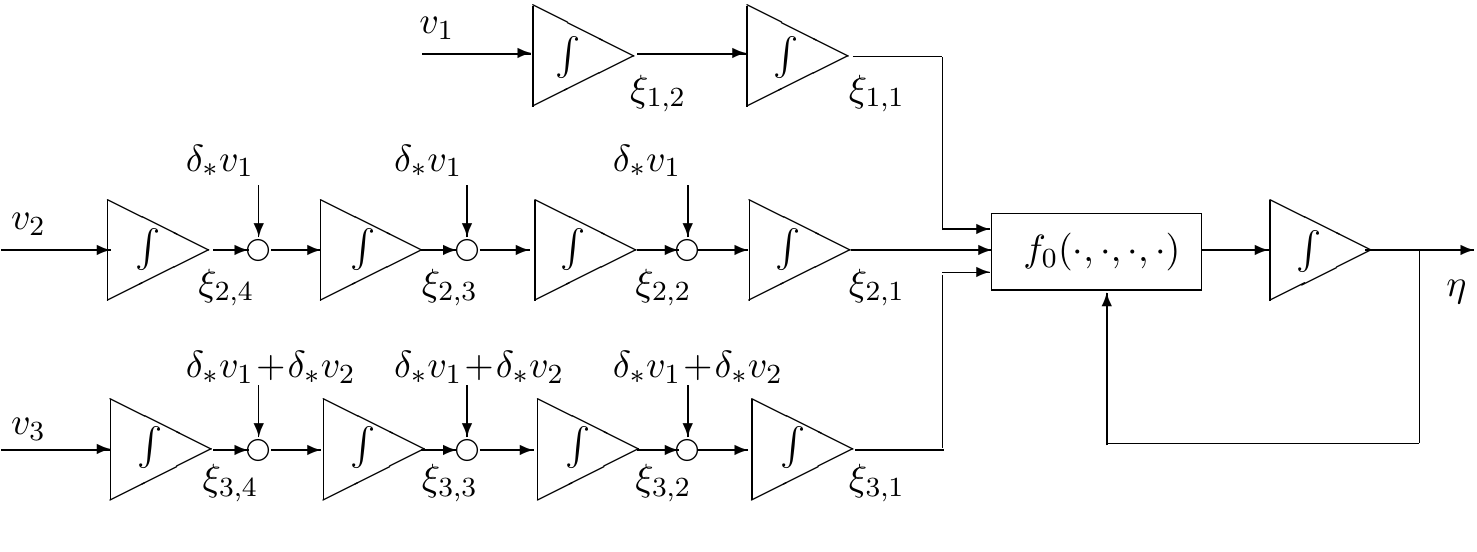}}
\caption{Nonlinear system in the form (\ref{zerody2_stabilization})
with three chains of integrators of
lengths $\{2,4,4\}$.}\label{fig244_old}
\end{figure}

On the other hand, under the same assumption,
the normal form (\ref{form1-square_stabilization})-(\ref{key2_stabilization}) would take the following simpler form
\be\label{equ_delta2}
\left\{\begin{array}{rcll}
\dot\eta &=&f_0(\eta,\xi_{1,1},\xi_{2,1},\xi_{3,1}),\cr
\dot\xi_{1,1} &=& \xi_{1,2},\cr
\dot\xi_{1,2} &=& v_1,\cr
\dot\xi_{2,1} &=& \xi_{2,2},\cr
\dot\xi_{2,2} &=& \xi_{2,3}+\delta_{2,2,1}(\eta,\xi_1,\xi_{2,1},\xi_{2,2},\xi_{3,1})v_1,\cr
\dot\xi_{2,3} &=& \xi_{2,4}+\delta_{2,3,1}(\eta,\xi_1,\xi_{2,1},\xi_{2,2},\xi_{2,3},\xi_{3,1})v_1,\cr
\dot\xi_{2,4} &=& v_2,\cr
\dot\xi_{3,1} &=& \xi_{3,2},\cr
\dot\xi_{3,2} &=& \xi_{3,3}+\delta_{3,2,1}(\eta,\xi_1,\xi_2,
     \xi_{3,1},\xi_{3,2})v_1,\cr
\dot\xi_{3,3} &=& \xi_{3,4}+\delta_{3,3,1}(\eta,\xi_1,\xi_2,\xi_{3,1},
    \xi_{3,2},\xi_{3,3})v_1,\cr
\dot\xi_{3,4} &=& v_3.
\end{array}\right.
\ee
Suppose that Assumption~\ref{assumption_stand} is satisfied.
The form (\ref{equ_delta2}) makes the
implementation of the chain-by-chain backstepping simpler,
due to the simpler structure (see Fig.~\ref{fig244}).

\begin{figure}[ht]
\centerline{\includegraphics[width=115mm]{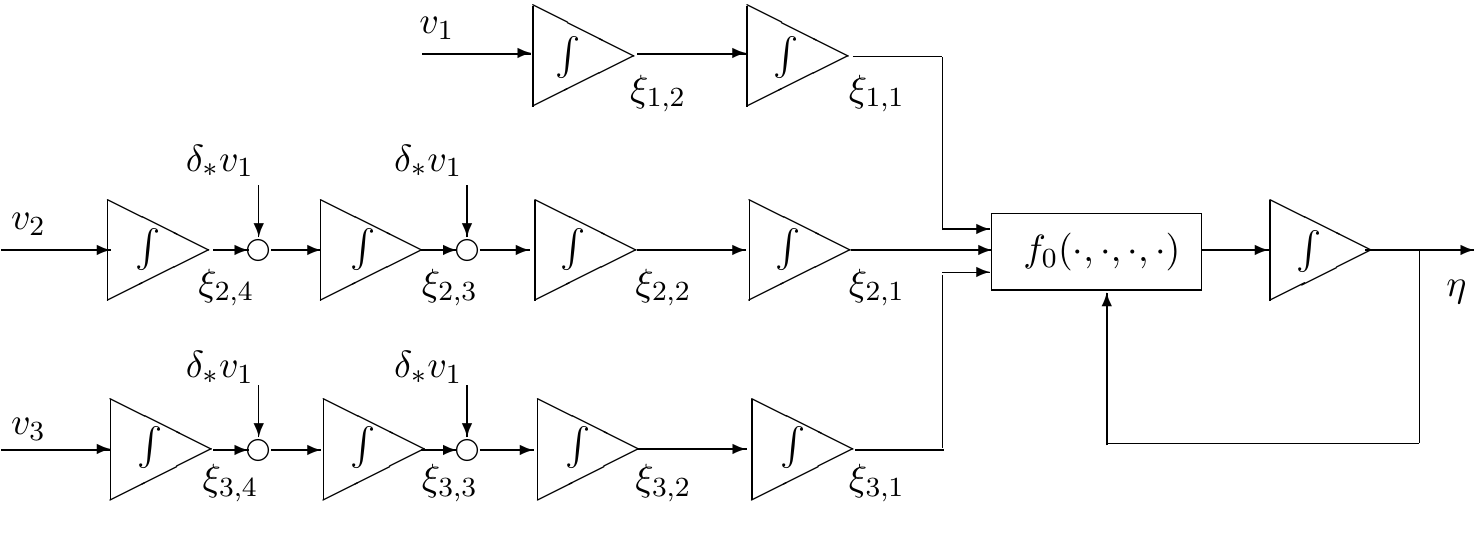}}
\caption{Nonlinear system in the form (\ref{form1-square_stabilization})-(\ref{key2_stabilization})
with three chains of integrators of
lengths $\{2,4,4\}$.}\label{fig244}
\end{figure}

\end{example}

\subsection{Level-by-Level Backstepping}

Let us call all $\xi_{i,1}$, {\em i.e.}, the ``leading" variables in
each chain of integrators which connect an input to an output,
the first level integrators, and call all $\xi_{i,2}$
the second level integrators, and so on.
As an alternative to the chain-by-chain backstepping, we here propose
to carry out the backstepping on all first level integrators, and then
repeat the procedure on all second level integrators until we reach to
last level of integrators. We will refer to such a backstepping
procedure as the level-by-level backstepping, in contrast with
the chain-by-chain backstepping procedure.

To make the level-by-level backstepping possible,
the coefficients $\delta_{i,j,l}$ in
the normal form (\ref{zerody2_stabilization}) with the property (\ref{key2_stabilization})
should satisfy the following assumption:

\begin{assumption}\label{tri2}
The functions $\delta_{i,j,l}$ depend only on variable
$\xi_{\ell_\rmp,\ell_\rmb}$, with
\ben
\item
$1\leq\ell_\rmp\leq m$ and $\ell_\rmb=1$; or,
\item
$\ell_\rmb\leq j-1$; or
\item
$\ell_\rmb=j$ and $\ell_\rmp\leq i$.
\een
\end{assumption}

We will say that
the coefficients $\delta_{i,j,l}$ in
the normal form (\ref{zerody2_stabilization})
have the chain-by-chain triangular dependency on state variables
if they satisfy Assumption~\ref{tri1}.
The coefficients $\delta_{i,j,l}$ in
the normal form (\ref{form1-square_stabilization}) - (\ref{key2_stabilization})
have the level-by-level triangular dependency on state variables
if they satisfy Assumption~\ref{tri2}.

Under Assumptions~\ref{assumption_stand} and \ref{tri2},
the level-by-level backstepping procedure
for the normal form (\ref{form1-square_stabilization}) with the property (\ref{key2_stabilization})
can be described as follows.
We will start with
\[
\dot\eta=f_0(\eta,\phi_{1,1}(\eta),\phi_{2,1}(\eta),\cdots,\phi_{m,1}(\eta)).
\]
After the first-level backstepping, we obtain the feedback laws
\[
v_i=v_i(\eta;\xi_{1,1},\xi_{2,1},\cdots,\xi_{i,1}),\quad i=1,2,\cdots,\alpha_1,
\]
where $\alpha_1$ is the number of chains that contain exactly one integrator,
{\em i.e.}, $n_1\!=\!n_2\!=\!\cdots\!=\!n_{\alpha_1}\!=\!1$.
For chains that contain more than one integrator, $\xi_{i,2}$ are viewed as virtual inputs,
and the desired $\xi_{i,2}$ are defined as
\[
\xi_{i,2}^\star=\phi_{i,2}(\eta;\xi_{1,1},\xi_{2,1},\cdots,\xi_{i,1}),\quad  i=\alpha_1+1,\alpha_1+2,\cdots,m.
\]
We next proceed with backstepping on the second level integrators.
After the second level backstepping, we obtain the feedback laws
\[
v_i=v_i(\eta;\xi_{1,1},\xi_{2,1},\cdots,\xi_{m,1};\xi_{\alpha_1+1,2},\xi_{\alpha_1+2,2},\cdots,\xi_{i,2}),
\quad
i=\alpha_1+1,\alpha_1+2,\cdots,\alpha_2,
\]
where $\alpha_2-\alpha_1$ is the number of chains that contain exactly two integrators,
{\em i.e.}, $n_{\alpha_1+1}=n_{\alpha_1+2}=\cdots=n_{\alpha_2}=2$.
For chains with lengths greater than $2$, the variables $\xi_{i,3}$ are viewed as virtual inputs,
and the desired $\xi_{i,3}$ are defined as
\[
\xi_{i,3}^\star=\phi_{i,3}(\eta; \xi_{1,1},\xi_{2,1}, \cdots,\xi_{m,1};\xi_{\alpha_1+1,2},\xi_{\alpha_1+2,2}
    \cdots,\xi_{i,2}),
\quad
\]
\[
\hspace{50mm}
i=\alpha_2+1,\alpha_2+2,\cdots,m.
\]
Continuing in this way, we finally obtain
\[
v_i=v_i(\eta; \xi_{1,1},\xi_{2,1},\cdots,\xi_{m,1};\xi_{\alpha_1+1,2},\xi_{\alpha_1+2,2},\cdots,\xi_{m,2};
\cdots;
\]
\[
\hspace{50mm}\xi_{\alpha_{n_m-1}+1,n_m},\xi_{\alpha_{n_m-1}+2,n_m},\cdots,\xi_{i,n_m}),
\]
for chains that contain $n_m$ integrators.

\begin{example}\rm
Consider a system in
the normal form (\ref{form1-square_stabilization})-(\ref{key2_stabilization})
 with three chains of integrators of lengths $\{2,4,4\}$. See Fig.~\ref{fig244}.
\be\label{equ_delta3}
\left\{\begin{array}{rcll}
\dot\eta &=& f_0(\eta,\xi_{1,1},\xi_{2,1},\xi_{3,1}),\cr
\dot\xi_{1,1} &=& \xi_{1,2},\cr
\dot\xi_{1,2} &=& v_1,\cr
\dot\xi_{2,1} &=& \xi_{2,2},\cr
\dot\xi_{2,2} &=& \xi_{2,3}+\delta_{2,2,1}(\eta,\xi_1,\xi_{2,1},\xi_{3,1},\xi_{2,2})v_1,\cr
\dot\xi_{2,3} &=& \xi_{2,4}+\delta_{2,3,1}(\eta,\xi_1,\xi_{2,1},\xi_{3,1},\xi_{2,2},\xi_{3,2},\xi_{2,3})v_1,\cr
\dot\xi_{2,4} &=& v_2,\cr
\dot\xi_{3,1} &=& \xi_{3,2},\cr
\dot\xi_{3,2} &=& \xi_{3,3}+\delta_{3,2,1}(\eta,\xi_1,\xi_{2,1},\xi_{3,1},
\xi_{2,2},\xi_{3,2})v_1,\cr
\dot\xi_{3,3} &=& \xi_{3,4}+\delta_{3,3,1}(\eta,\xi_1,\xi_{2,1},\xi_{3,1},\xi_{2,2},\xi_{3,2},\xi_{2,3},\xi_{3,3})v_1,\cr
\dot\xi_{3,4} &=& v_3.
\end{array}\right.
\ee

Clearly, Assumption~\ref{tri2} is satisfied, but Assumption~\ref{tri1} is not.
Consequently, the chain-by-chain backstepping cannot be implemented on this system.
In what follows, we will illustrate how to implement
the level-by-level backstepping on this system.

Let Assumption~\ref{assumption_stand} be satisfied, {\em i.e.},
there exist smooth functions $\phi_{i,1}(\eta)$,
with $\phi_{i,1}(0)=0$, $i=1,2,3$, such that
the equilibrium $\eta=0$ of the subsystem
\be\label{eta344_1}
\dot\eta=f_0(\eta,\phi_{1,1}(\eta),\phi_{2,1}(\eta),\phi_{3,1}(\eta))
\ee
is globally asymptotically stable.
The backstepping procedure starts with the subsystem (\ref{eta344_1}).
Now we consider the backstepping on the first level variables.
The variable $\xi_{1,1}$ can be viewed as the virtual input
of the subsystem $\dot\eta=f_0(\eta,\xi_{1,1},\phi_{2,1}(\eta),\phi_{3,1}(\eta))$,
and
the desired input is given by $\xi_{1,1}^\star=\phi_{1,1}(\eta)$.
To carry out the backstepping from $\xi_{1,1}$ to $\xi_{1,2}$, we consider the subsystem
\[
\left\{\begin{array}{rcl}
\dot\eta&=&f_0(\eta,\xi_{1,1},\phi_{2,1}(\eta),\phi_{3,1}(\eta)),\cr
\dot\xi_{1,1}&=&\xi_{1,2}.\cr
\end{array}\right.
\]
with $\xi_{1,2}$ as the virtual input.
By Lemma~\ref{base1}, this subsystem can be
globally asymptotically stabilized by a control of the form
\be\label{xi12}
\xi_{1,2}^\star=\phi_{1,2}(\eta,\xi_{1,1},\phi_{2,1}(\eta),\phi_{3,1}(\eta)).
\ee
The variable $\xi_{2,1}$ can be viewed as the virtual input
of the subsystem
\[
\left\{\begin{array}{rcl}
\dot\eta&=&f_0(\eta,\xi_{1,1},\xi_{2,1},\phi_{3,1}(\eta)),\cr
\dot\xi_{1,1}&=&\xi_{1,2}^\star,
\end{array}\right.
\]
and $\xi_{2,1}^\star=\phi_{2,1}(\eta)$ globally asymptotically stabilizes
its equilibrium \newline $\col\{\eta,\xi_{1,1}\}=0$.
To carry out the backstepping from $\xi_{2,1}$ to $\xi_{2,2}$, we
next look at the subsystem
\[
\left\{\begin{array}{rcl}
\dot\eta&=&f_0(\eta,\xi_{1,1},\xi_{2,1},\phi_{3,1}(\eta)),\cr
\dot\xi_{1,1}&=&\xi_{1,2}^\star,\cr
\dot\xi_{2,1} &=& \xi_{2,2},\cr
\end{array}\right.
\]
with $\xi_{2,2}$ as the virtual input. By Lemma~\ref{base1},
this subsystem can be
globally asymptotically stabilized by a control of the form
\be\label{xi22}
\xi_{2,2}^\star=\phi_{2,2}(\eta,\xi_{1,1},\xi_{2,1},\phi_{3,1}(\eta)).
\ee
Similarly, to backstep from $\xi_{3,1}$ to $\xi_{3,2}$, we consider
\[
\left\{\begin{array}{rcl}
\dot\eta&=&f_0(\eta,\xi_{1,1},\xi_{2,1},\xi_{3,1}),\cr
\dot\xi_{1,1}&=&\xi_{1,2}^\star,\cr
\dot\xi_{2,1} &=& \xi_{2,2}^\star,\cr
\dot\xi_{3,1} &=& \xi_{3,2},\cr
\end{array}\right.
\]
with $\xi_{3,2}$ as the virtual input. The subsystem can be
globally asymptotically stabilized by a control of the form
\be\label{xi32}
\xi_{3,2}^\star=\phi_{3,2}(\eta,\xi_{1,1},\xi_{2,1},\xi_{3,1}).
\ee
Thus, after the first level backstepping, the subsystem
\[
\left\{\begin{array}{rcl}
\dot\eta&=&f_0(\eta,\xi_{1,1},\xi_{2,1},\xi_{3,1}),\cr
\dot\xi_{1,1}&=&\xi_{1,2}^\star,\cr
\dot\xi_{2,1} &=& \xi_{2,2}^\star,\cr
\dot\xi_{3,1} &=& \xi_{3,2}^\star,\cr
\end{array}\right.
\]
can be written as
\be\label{etarmI}
\dot{\eta}_\rmI=f_\rmI(\eta_\rmI,\xi_{1,2},\xi_{2,2},\xi_{3,2}),
\ee
where $\eta_\rmI=\col\{\eta,\xi_{1,1},\xi_{2,1},\xi_{3,1}\}$,
and
$\xi_{1,2}$, $\xi_{2,2}$ and $\xi_{3,2}$ are the virtual inputs.
The equilibrium $\eta_\rmI=0$ of this subsystem (\ref{etarmI}) is
globally asymptotically stabilized by the virtual inputs
$\xi_{1,2}$, $\xi_{2,2}$ and $\xi_{3,2}$ as given by
(\ref{xi12}), (\ref{xi22}) and
(\ref{xi32}), respectively.

For the second level backstepping, consider
% the subsystem
\be\label{level2}
\left\{\begin{array}{rcl}
\dot{\eta}_\rmI&=&f_\rmI(\eta_\rmI,\xi_{1,2},\xi_{2,2},\xi_{3,2}),\cr
\dot\xi_{1,2} &=& v_1,\cr
\dot\xi_{2,2} &=& \xi_{2,3}+\delta_{2,2,1}(\eta,\xi_{1,1},\xi_{2,1},\xi_{3,1},\xi_{2,2})v_1,\cr
\dot\xi_{3,2} &=& \xi_{3,3}+\delta_{3,2,1}(\eta,\xi_{1,1},\xi_{2,1},\xi_{3,1},\xi_{2,2},\xi_{3,2})v_1,\cr
\end{array}\right.
\ee
and view $\xi_{2,3}$ and $\xi_{3,3}$ as its virtual inputs.
Following the same procedure as in the first level backstepping,
we find the controls of the form
\be\label{xi13}
\left\{\begin{array}{rcl}
v_1 &=& v_1(\eta_\rmI,\xi_{1,2},\xi_{2,2}^\star,\xi_{3,2}^\star),\cr
\xi_{2,3}^\star&=&\phi_{2,3}(\eta_\rmI,\xi_{1,2},\xi_{2,2},\xi_{3,2}^\star),\cr
\xi_{3,3}^\star&=&\phi_{3,3}(\eta_\rmI,\xi_{1,2},\xi_{2,2},\xi_{3,2})
\end{array}\right.
\ee
that globally asymptotically stabilize the equilibrium
\[
\eta_{\rmI\!\rmI}=\col\{\eta_\rmI,\xi_{1,2},\xi_{2,2},\xi_{3,2}\}=0
\]
of the subsystem (\ref{level2}).
In other word, the subsystem (\ref{level2})
can be written as
\[
\dot{\eta}_{\rmI\!\rmI}=f_{\rmI\!\rmI}(\eta_{\rmI\!\rmI},v_1,\quad
   \xi_{2,3},\xi_{3,3}),
\]
whose equilibrium $\eta_{\rmI\!\rmI}=0$
is globally asymptotically stabilized by the input $v_1$, and virtual inputs $\xi_{2,3}$ and $\xi_{3,3}$ given by
(\ref{xi13}).

Define
\be\label{xi13a}
\left\{\begin{array}{rcl}
\dot{\eta}_{\rmI\!\rmI} &=& f_{\rmI\!\rmI}(\eta_{\rmI\!\rmI},v_1,\xi_{2,3},\xi_{3,3}),\cr
\dot\xi_{2,3} &=& \xi_{2,4}+\delta_{2,3,1}(\eta,\xi_1,\xi_{2,1},\xi_{3,1},\xi_{2,2},\xi_{3,2},\xi_{2,3})v_1,\cr
\dot\xi_{3,3} &=& \xi_{3,4}+\delta_{3,3,1}(\eta,\xi_1,\xi_{2,1},\xi_{3,1},\xi_{2,2},\xi_{3,2},\xi_{2,3},\xi_{3,3})v_1,\cr
\end{array}\right.
\ee
on which we carry out the third level of backstepping to obtain
\be\label{xi13a_input}
\left\{\begin{array}{rcl}
\xi_{2,4}^\star&=&\phi_{2,4}(\eta_{\rmI\!\rmI},\xi_{2,3},\xi_{3,3}^\star),\\
\xi_{3,4}^\star&=&\phi_{3,4}(\eta_{\rmI\!\rmI},\xi_{2,3},\xi_{3,3}).
\end{array}\right.
\ee
The subsystem (\ref{xi13a})
can be defined as
\[
\dot{\eta}_{\rmI\!\rmI\!\rmI}=f_{\rmI\!\rmI}(\eta_{\rmI\!\rmI\!\rmI},v_1,\quad
   \xi_{2,4},\xi_{3,4}),
\]
whose equilibrium $\eta_{\rmI\!\rmI\!\rmI}=\col\{\eta_{\rmI\!\rmI},\xi_{2,3},\xi_{3,3}\}=0$
is globally asymptotically stabilized by the virtual inputs $\xi_{2,4}$ and $\xi_{3,4}$ given by
(\ref{xi13a_input}).

Finally, for the fourth level backstepping, we define
\[
\left\{\begin{array}{rcl}
\dot{\eta}_{\rmI\!\rmI\!\rmI} &=& f_{\rmI\!\rmI\!\rmI}(\eta_{\rmI\!\rmI\!\rmI},\xi_{2,4},\xi_{3,4}),\cr
\dot\xi_{2,4} &=& v_2,\cr
\dot\xi_{3,4} &=& v_3,
\end{array}\right.
\]
on which we carry out the last level of backstepping to obtain
\[
\left\{\begin{array}{rcl}
v_2&=&v_2(\eta_{\rmI\!\rmI\!\rmI},\xi_{2,4},\xi_{3,4}^\star),\\
v_3&=&v_3(\eta_{\rmI\!\rmI\!\rmI},\xi_{2,4},\xi_{3,4}).
\end{array}\right.
\]
The inputs $v_1$, $v_2$ and $v_3$
%as given by (\ref{v1}) and (\ref{v3})
globally asymptotically
stabilize the origin of the system
(\ref{equ_delta3}).
\end{example}

\begin{rema}
%The normal form (\ref{zerody2_stabilization}) with
The structural property (\ref{key2_stabilization})
makes the level-by-level backstepping possible.
It is not possible to implement the level-by-level backstepping
technique on the normal form (\ref{zerody2_stabilization}).
For example, consider a system in the form (\ref{zerody2_stabilization})
with two chains of integrators of lengths $\{2,3\}$, backstepping
the virtual input from $\xi_{2,1}=\phi_{2,1}(\eta)$
to $\xi_{2,2}$ by the dynamical equation
$
\dot\xi_{2,1}=\xi_{2,2}+\delta_{2,1,1}(\eta,\xi)v_1
$
is infeasible. At this stage, $v_1$ is not yet available.
\end{rema}

\subsection{Mixed Chain-by-Chain and Level-by-Level
Backstepping}

A system with a vector relative degree
is a special case of the systems (\ref{zerody2_stabilization})
with all $\delta_{i,j,l}=0$.
Thus, both chain-by-chain backstepping and level-by-level backstepping
can be implemented on it. Furthermore,
backstepping can be switched across chains and levels as long as
a variable of lower level in a chain is backstepped
earlier than variables of higher levels in the same chain.

In the absence of a vector relative degree,
the normal form (\ref{form1-square_stabilization}) with the property (\ref{key2_stabilization})
contains coefficient
functions $\delta_{i,j,l}$. The implementation of both
chain-by-chain and level-by-level backstepping require structural
dependency on state variables of $\delta_{i,j,l}$. Such structural
dependency constraint can be weakened by utilizing mixed chain-by-chain and level-by-level
backstepping.

\begin{example}\rm
Consider a system with three chains of integrators of lengths $\{2,4,4\}$.
\be
\left\{\begin{array}{rcll}
\dot\eta &=& f_0(\eta,\xi_{1,1},\xi_{2,1},\xi_{3,1}),\cr
\dot\xi_{1,1} &=& \xi_{1,2},\cr
\dot\xi_{1,2} &=& v_1,\cr
\dot\xi_{2,1} &=& \xi_{2,2},\cr
\dot\xi_{2,2} &=& \xi_{2,3},\cr
\dot\xi_{2,3} &=& \xi_{2,4}+\xi_{3,2}v_1,\cr
\dot\xi_{2,4} &=& v_2,\cr
\dot\xi_{3,1} &=& \xi_{3,2},\cr
\dot\xi_{3,2} &=& \xi_{3,3},\cr
\dot\xi_{3,3} &=& \xi_{3,4}+\xi_{2,4}v_1,\cr
\!\!\dot\xi_{3,4} &=& v_3.
\end{array}\right.
\ee

Let Assumption \ref{assumption_stand} be satisfied.
Due to the term $\xi_{3,2}v_1$, Assumption \ref{tri1} is not satisfied.
Similarly, because of the term $\xi_{2,4}v_1$, Assumption \ref{tri2} is not met.
%It is obvious that neither Assumption \ref{tri1} nor Assumption \ref{tri2} is satisfied.
As a result, neither the chain-by-chain nor the level-by-level backstepping
can be implemented on this system.
However, a mixed chain-by-chain and level-by-level backstepping
will successfully stabilize this system.
In particular, by Lemma~\ref{base1}, we can carry out backstepping in the
order of
$\{\xi_{1,1},\xi_{2,1},\xi_{3,1},\xi_{1,2},\xi_{2,2},\xi_{3,2},\xi_{2,3},\xi_{2,4},\xi_{3,3},$ $\xi_{3,4}\}$
to obtain %the stabilizing feedback laws
\[
\left.\begin{array}{rcl}
v_1&=&v_1(\eta,\xi_1),\cr
v_2&=&v_2(\eta,\xi_1,\xi_2,\xi_{3,1},\xi_{3,2}),\cr
v_3&=&v_3(\eta,\xi_1,\xi_2,\xi_3).
\end{array}\right.
\]
\end{example}

\begin{figure}[ht]
Chain-by-Chain Backstepping\newline
\centerline{\includegraphics[width=115mm]{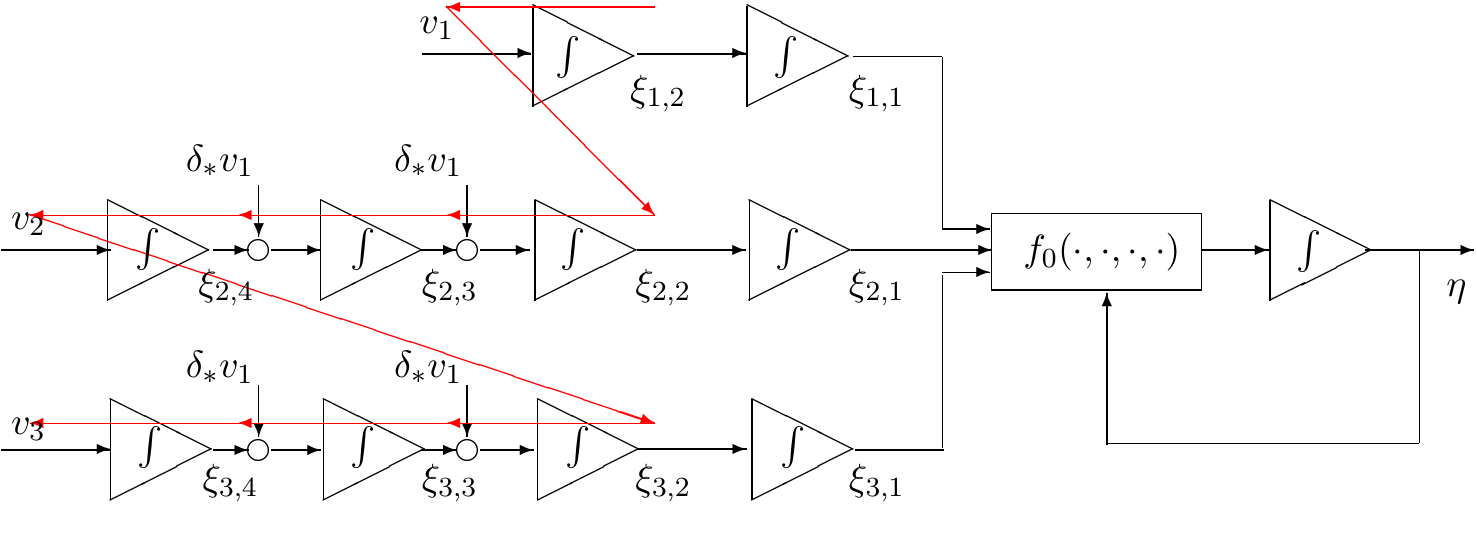}}
Level-by-Level Backstepping\newline
\centerline{\includegraphics[width=115mm]{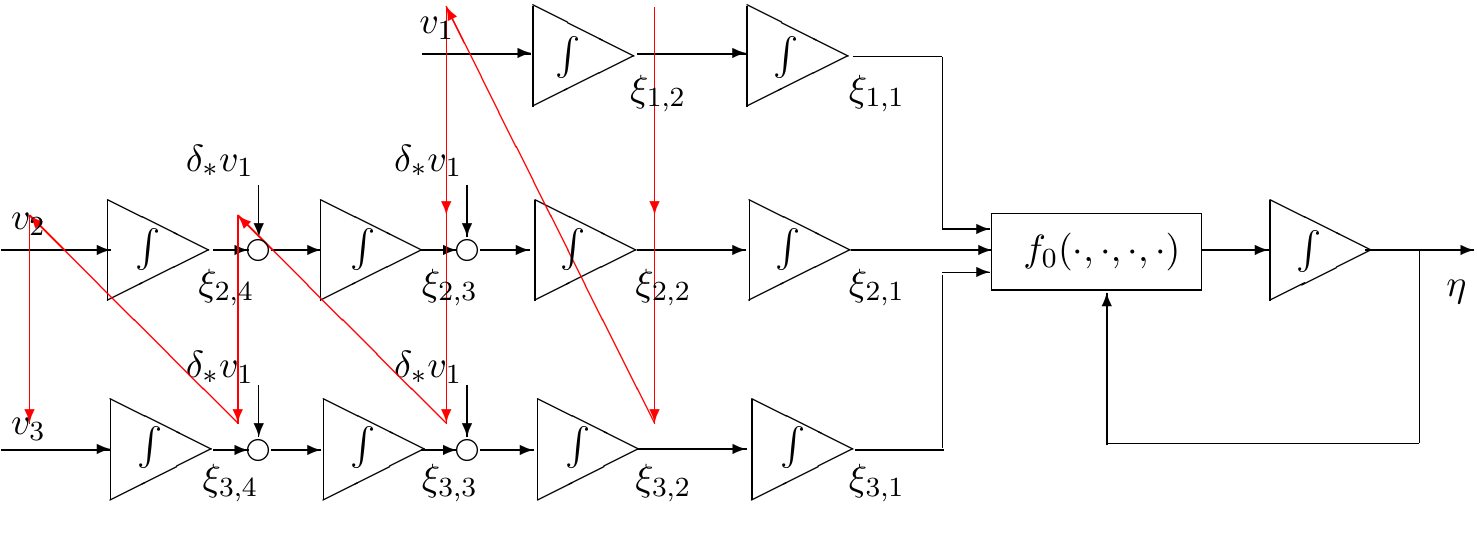}}
Mixed Chain-by-Chain and Level-by-Level Backstepping\newline
\centerline{\includegraphics[width=115mm]{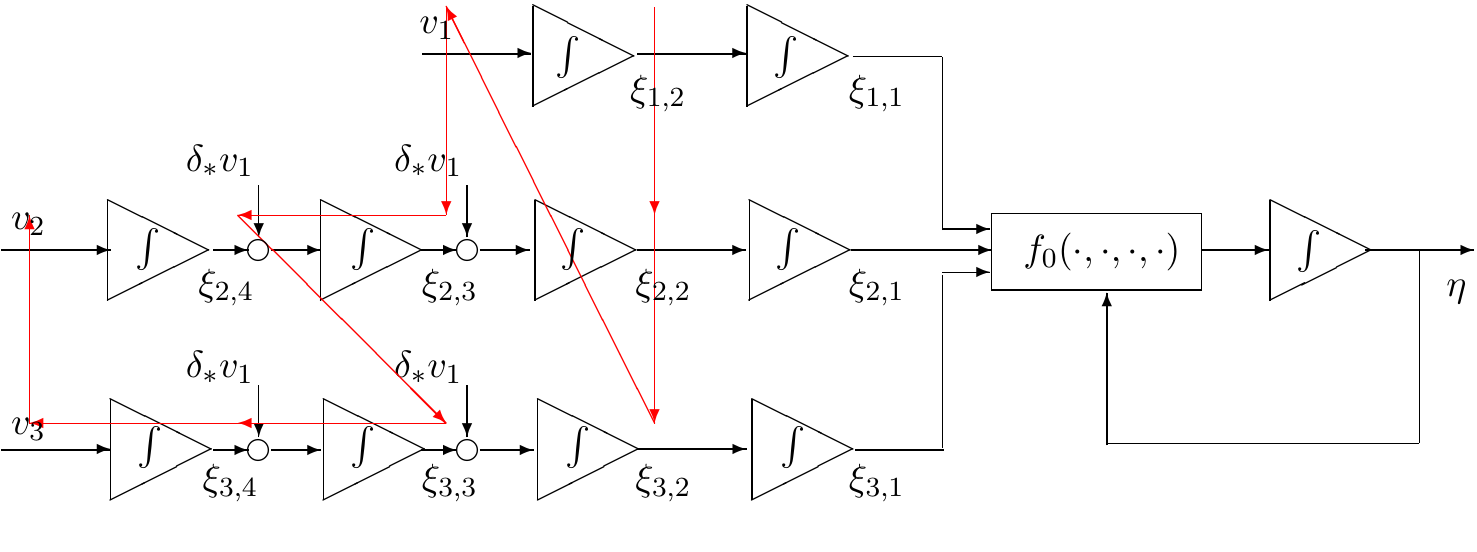}}
\caption{Backstepping procedures in nonlinear system with three chains of integrators of
lengths $\{2,4,4\}$}.
\label{all_fig_backstepping}
\end{figure}

Shown in Fig.~\ref{all_fig_backstepping} are backstepping procedures for the systems
in the normal form
(\ref{form1-square_stabilization})-(\ref{key2_stabilization})
with three chains of integrators of lengths $\{2,4,4\}$.

%\section{A General Result of Backstepping}

Motivated by the mixed chain-by-chain and level-by-level
backstepping, we give the following result, which includes all
the above backstepping procedures as special cases.

\begin{theo}\label{stab_general}
Consider a system in the form
\be\label{mix_form}
\left\{\begin{array}{rcl}
\dot \eta &=& f_0(\eta,\xi_{1,1},\xi_{2,1},\cdots,\xi_{m,1}),\cr
\dot\xi_{i,j} &=& \xi_{i,j+1}+ \sum_{l=1}^{i-1} \delta_{i,j,l}(\eta,\xi)v_l,\;\; j=1,2,\cdots,n_i-1,\cr
\dot\xi_{i,n_i} &=& v_i,\cr
y_i &=& \xi_{i,1},\;\; i=1,2,\cdots,m,
\end{array}
\right.
\ee
where $\xi=\col\{\xi_1,\xi_2,\cdots,\xi_m\}$, $\xi_i=\col\{\xi_{i,1},\xi_{i,2},\cdots,\xi_{i,n_i}\}$,\newline
$v=\col\{v_1,v_2,\cdots,v_m\}$, $n_1\leq n_2\leq \cdots \leq n_m$, and all functions
 $\delta_{i,j,l}$ and $f_0$ are smooth.
Then, there exists a feedback
$v=v(\eta,\xi)$ that globally asymptotically
stabilizes the system at $(\eta,\xi)=0$ if
\ben
\item
there exist
$\phi_{i,1}(\eta)$,
$i=1,2,\cdots,m$,
such that \newline
$\dot\eta=f_0(\eta,\phi_{1,1}(\eta),\phi_{2,1}(\eta),\cdots,\phi_{m,1}(\eta))$
is globally asymptotically stable at its equilibrium $\eta=0$; and
\item
there exists an ordered list $\kappa$ containing all variables of $\xi$
such that, for $j=1,2,\cdots,n_i-1$, $l=1,2,\cdots, i-1$, $i=1,2,\cdots,m$,
\ben
\item
$\xi_{i,j}$ appears earlier than $\xi_{i,j+1}$ in $\kappa$;
\item for $\delta_{i,j,l}\neq0$, the variables $\xi_l$ appear earlier than $\xi_{i,j}$
in $\kappa$;
\item $\delta_{i,j,l}$ depends only on $\eta$, $\xi_{\ell,1}$, $\ell=1,2,\cdots,m$, and
the variables that appear no later than $\xi_{i,j}$
in $\kappa$.
\een
\een
\end{theo}

The backstepping procedures can be carried out according to the order of $\kappa$.
In some cases, there exist more than one $\kappa$.
Backstepping in different orders lead to different dependency of controls
on state variables, which can be exploited to meet certain constraints
or performance requirement.

\begin{example}\rm
Consider a system in the form of
(\ref{mix_form}),
\[
\left\{\begin{array}{rcl}
\dot\eta &=& \eta+\xi_{1,1}+\xi_{2,1},\cr
\dot\xi_{1,1} &=& v_1,\cr
\dot\xi_{2,1} &=& \xi_{2,2}+\xi_{3,2}v_1,\cr
\dot\xi_{2,2} &=& v_2,\cr
\dot\xi_{3,1} &=& \xi_{3,2},\cr
\dot\xi_{3,2} &=& \xi_{3,3},\cr
\dot\xi_{3,3} &=& \xi_{3,4}+\xi_{2,2}v_1,\cr
\dot\xi_{3,4} &=& v_3.
\end{array}\right.
\]
The zero dynamics $\dot \eta=\eta$ can be stabilized by
\[
\xi_{1,1}^\star=0,\quad \xi_{2,1}^\star=-2\eta,\quad \xi_{3,1}^\star=0.
\]
Neither Assumption \ref{tri1} nor Assumption \ref{tri2} is satisfied.
So neither the chain-by-chain nor the level-by-level backstepping
can be carried out.
However, the system satisfies the conditions of Theorem~\ref{stab_general} with
the ordered list
\[
\kappa=\{\xi_{1,1},\xi_{3,1},\xi_{3,2},\xi_{2,1},\xi_{2,2},\xi_{3,3},\xi_{3,4}\}.
\]
By backstepping in the order of $\kappa$,
the stabilizing controller is given as
\[
\left\{\begin{array}{rcll}
v_1&=&-\eta,\cr
v_2&=&-6\eta-5\xi_{1,1}-6\xi_{2,1}-3\xi_{2,2}+\eta(-\xi_{3,1}+3\xi_{3,2})+\xi_{3,2}(\xi_{1,1}+\xi_{2,1}),\cr
v_3&=&-\xi_{3,1}-3\xi_{3,2}-5\xi_{3,3}-3\xi_{3,4}+(\eta+\xi_{1,1}+\xi_{2,1})(\beta-\eta\xi_{3,2})\cr
&&\; +\eta(\beta+\eta+3\xi_{1,1}+3\xi_{2,1}+4\xi_{2,2}-2\eta\xi_{3,2}-\eta\xi_{3,3}+2v_2),\cr
\end{array}\right.
\]
where $\beta=3\eta+2\xi_{1,1}+2\xi_{2,1}+2\xi_{2,2}-\eta\xi_{3,2}$,
the Lyapunov function
\[
V\!\!=\!\![\eta^2+\xi_{1,1}^2+\xi_{3,1}^2+\xi_{3,2}^2+(\xi_{2,1}+2\eta)^2
+\beta^2
+(\xi_{3,3}+\xi_{3,1}+\xi_{3,2})^2
+(\xi_{3,4}+\beta)^2]/2,
\]
and its derivative
\[
\dot V=-\eta^2-\xi_{3,2}^2-(\xi_{2,1}+2\eta)^2-\beta^2-(\xi_{3,3}+\xi_{3,1}+\xi_{3,2})^2-(\xi_{3,4}+\beta)^2\leq 0.
\]
\begin{figure}[ht!]
\centerline{
\includegraphics[width=100mm]{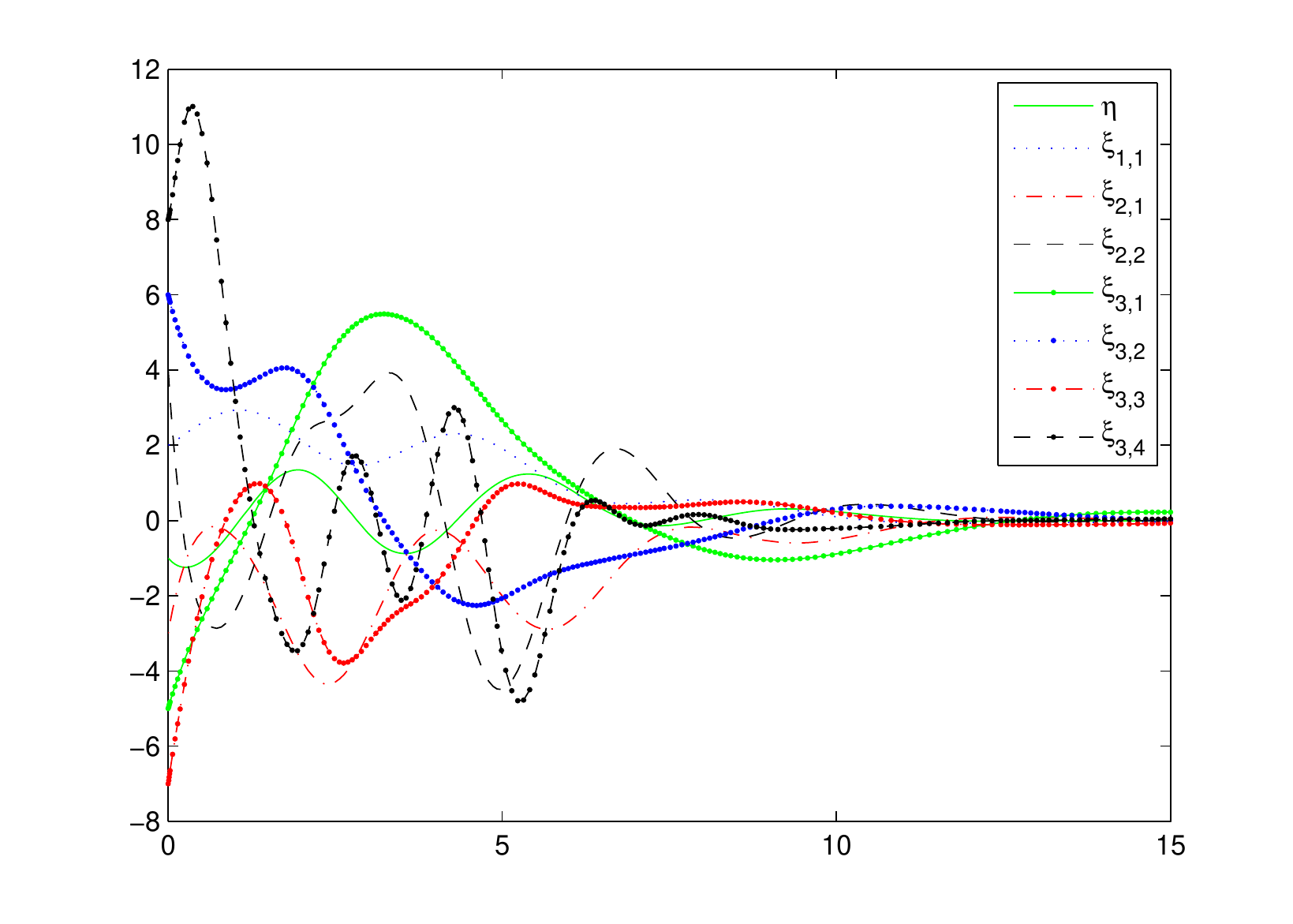}}
\centerline{
\includegraphics[width=100mm]{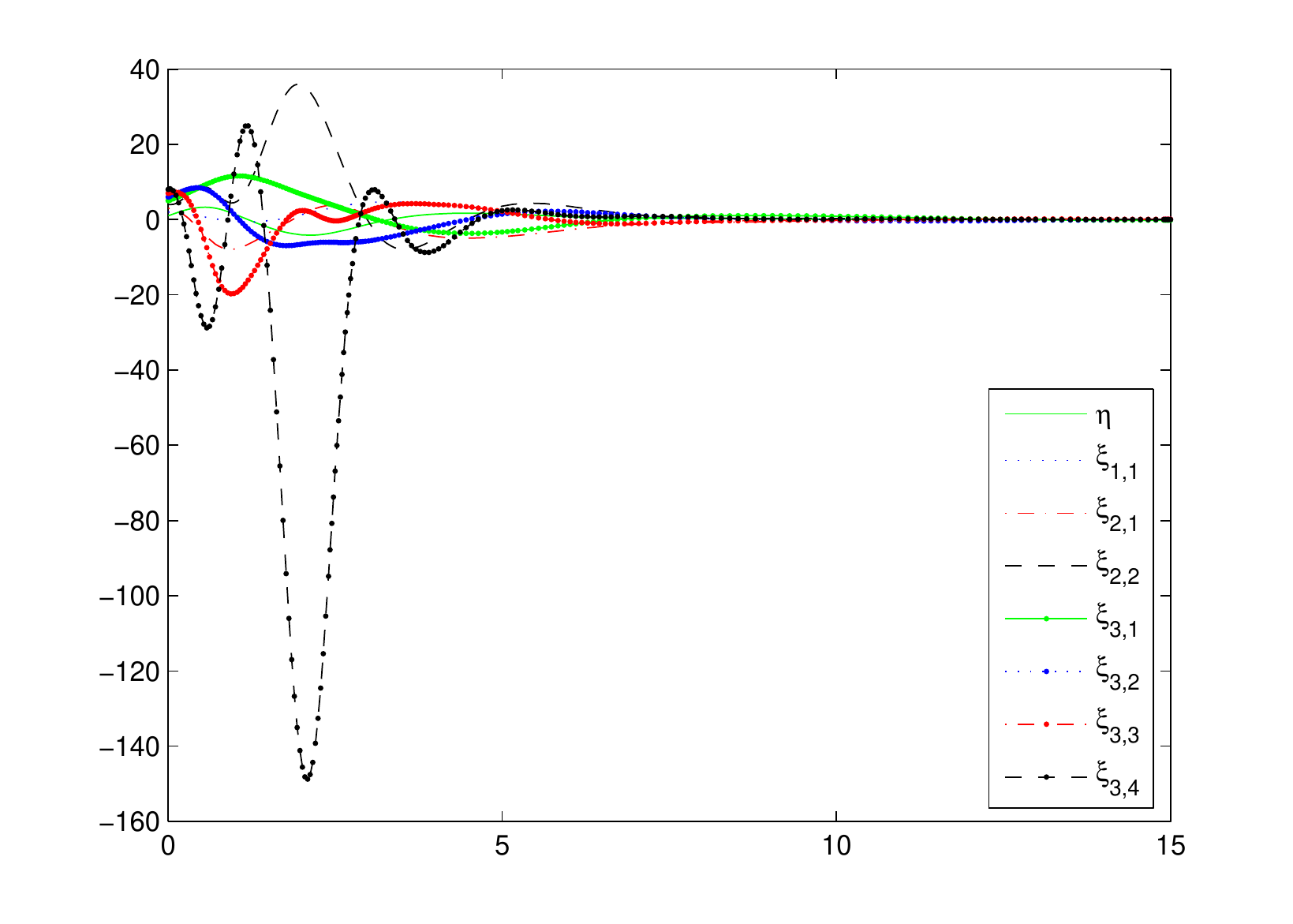}}
\vspace{-3mm}
\caption{State trajectories with two sets of initial conditions.
}
\label{mix124}
\end{figure}
It can be verified that no solution other than $(\eta(t),\xi(t))= 0$ can stay forever in $\{(\eta,\xi)\in\RR^8\; : \; \dot V(\eta,\xi)=0\}$.
Thus, by LaSalle theorem,
the closed-loop system is globally asymptotically stable at the origin.
Shown in Fig.~\ref{mix124} are some state trajectories of the closed-loop system with different initial conditions.

\end{example}

Next, we explain how different backstepping procedures leading to different control performance by a simple linear numerical example.

\begin{example}\rm
Consider the system (\ref{vrd3}) with a vector relative degree $\{2,2\}$,
\be
\left\{\begin{array}{rcl}
\dot \eta &=& \eta+\xi_{1,1}+\xi_{2,1},\cr
\dot\xi_{1,1} &=& \xi_{1,2},\cr
\dot\xi_{1,2} &=& v_1,\cr
\dot\xi_{2,1} &=& \xi_{2,2},\cr
\dot\xi_{2,2} &=& v_2,
\end{array}\right.
\ee
The zero dynamics here is linear for the convenience of backstepping. We certainly can use linear system tools to design controllers.

Here, in each step, integrator backstepping with $c=1$ is implemented (see Lemma 2.8 in \cite{krna95} for detail).
The backstepping starts with
\[
\xi_{1,1}^\star=-\eta,\quad \xi_{2,1}^\star=-\eta.
\]
We first carry out chain-by-chain backstepping in the order of \newline
$\{\xi_{1,1},\xi_{1,2},\xi_{2,1},\xi_{2,2}\}$,
and obtain
\be\label{u_chain}
\left\{\begin{array}{rcl}
v_1&=&-3\eta-5\xi_{1,1}-3\xi_{1,2},\cr
v_2&=&-11\eta-4\xi_{1,1}-6\xi_{1,2}-11\xi_{2,1}-3\xi_{2,2}.\cr
\end{array}\right.
\ee
The Lyapunov function of proving stability of the closed-loop system is
\begin{eqnarray*}
V&=&[\eta^2+(\xi_{1,1}+\eta)^2+(\xi_{1,2}+2\eta+2\xi_{1,1})^2+(\xi_{2,1}+\eta)^2\\
&&\hspace{40mm}+(\xi_{2,2}+8\eta+6\xi_{1,1}+2\xi_{1,2}+2\xi_{2,1})^2]/2
\end{eqnarray*}
and $\dot V=-2V$.

The level-by-level backstepping can be implemented in the order of\newline
$\{\xi_{1,1},\xi_{2,1},\xi_{1,2},\xi_{2,2}\}$
to arrive at
\be\label{u_level}
\left\{\begin{array}{rcl}
v_1&=&-5\eta-5\xi_{1,1}-3\xi_{1,2}-2\xi_{2,1},\cr
v_2&=&-9\eta-6\xi_{1,1}-2\xi_{1,2}-6\xi_{2,1}-2\xi_{2,2},\cr
\end{array}\right.
\ee
%the Lyapunov function
\[
V=[\eta^2+(\xi_{1,1}+\eta)^2+(\xi_{2,1}+\eta)^2+(\xi_{1,2}+2\eta+2\xi_{1,1})^2+(\xi_{2,2}+4\eta+2\xi_{1,1}+\xi_{2,1})^2]/2
\]
and $\dot V=-2V$.
\begin{figure}[ht!]
\centerline{
\includegraphics[width=120mm]{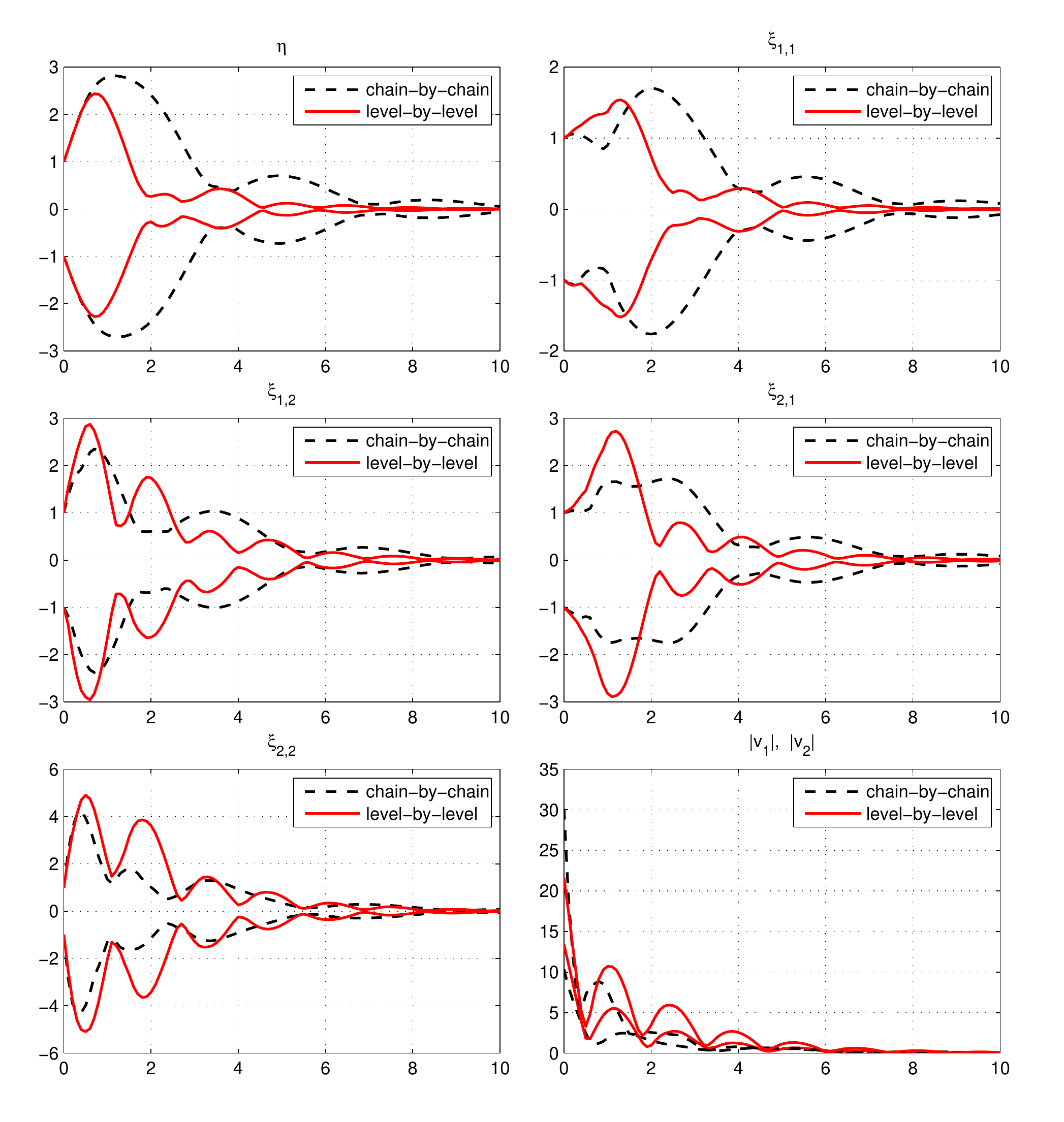}}
\vspace{-7mm}
\caption{Outline of state trajectories with different initial conditions.}
\label{vect_stab}
\end{figure}
\begin{figure}[ht!]
\centerline{
\includegraphics[width=120mm]{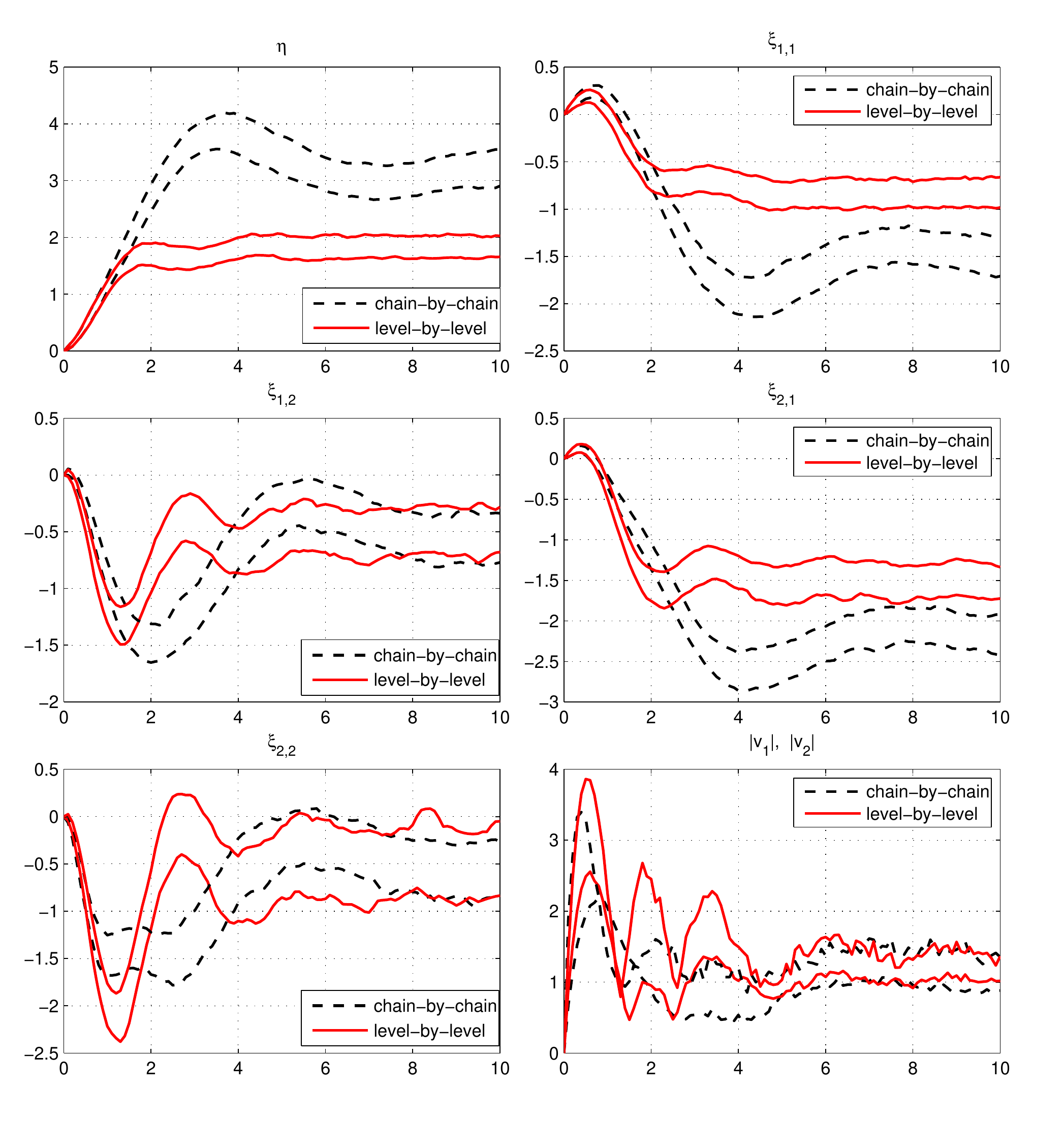}}
\vspace{-7mm}
\caption{Outline of state trajectories with disturbances.}
\label{vect_dis}
\end{figure}

Comparing the four feedback gains in (\ref{u_chain}) and (\ref{u_level}), $v_1$ obtained in the chain-by-chain backstepping is the smallest, and
$v_2$ obtained in the chain-by-chain backstepping is the largest, while $v_1$ and $v_2$ for level-by-level backstepping have average gains.
The input $v_1$ obtained in the chain-by-chain backstepping does not depend on $\xi_{2,1}$, while
$v_1$ obtained in the level-by-level backstepping does.
If there is disturbance in  $\xi_{2,1}$, both inputs of level-by-level backstepping can take account of it more directly.

Run simulation 1000 times with the initial conditions being uniformly distributed pseudo-random numbers in the interval $[-1, 1]$,
Figure.~\ref{vect_stab} shows outlines of each state trajectory, and outlines of absolute values of each input.
It indicates that state trajectories of level-by-level backstepping converge faster.

Consider adding the following disturbances to each state equation,
\[
w_i(t)=\mbox{rand}_{i,j}, \quad t\in[0.01j, \;0.01(j+1)],\quad  j\geq 0, \quad 1\leq i\leq 5,
\]
where $\mbox{rand}_{i,j}$ is an uniformly distributed pseudo-random number in the interval $[0, 1]$.
Run simulation 1000 times again with disturbance.
% and with the initial conditions being uniformly distributed pseudo-random numbers in the interval $[-1, 1]$,
Figure.~\ref{vect_dis} shows outlines of each state and absolute values of each input.
It indicates that the feedback law from level-by-level backstepping has a better disturbance rejection capability.

This example shows that backstepping in different orders lead to different dependency of controls
on state variables and different performance,
which can be exploited to meet certain constraints
or performance requirement.
In general, we could select a backstepping order with each input depend on more state variables.
In this case, each input would take full advantage of the state information, and thus the closed-loop systems tend to have a better performance.

\end{example}

\section{Summary of the Chapter}

We exploited the properties of a recently developed
structural decomposition for the stabilization of multiple input and multiple output
systems, and showed that this
decomposition simplifies the conventional chain-by-chain backstepping design and
motivates a new level-by-level backstepping design procedure that is able to
stabilize some systems for which the conventional backstepping procedure is
not applicable.
The chain-by-chain backstepping and level-by-level backstepping
can be combined to form a mixed backstepping design
technique. The enlarged class of systems that can be stabilized by this mixed backstepping
design procedure is characterized in the form of a theorem.

\chapter{Semi-global Stabilization for Nonlinear Systems}

In Chapter 3, we developed a structural decomposition for
multiple input multiple output
nonlinear systems that are affine in control but
otherwise general. Chapter 4 shows that this structural decomposition simplifies
the conventional backstepping design and allows a new backstepping design
procedure that is able to stabilize some systems on which the conventional
backstepping is not applicable.
In this chapter we further exploit the properties of such
a decomposition for the purpose of solving the semi-global stabilization problem
for minimum phase nonlinear systems without vector relative degrees.
By taking advantage of special structure of the decomposed system, we first apply
the low gain design to the part of system that possesses a linear dynamics.
The low gain design results in an augmented zero dynamics that is locally stable
at the origin with a domain of attraction that can be made arbitrarily large by
lowering the gain. With this augmented zero dynamics, backstepping design is
then apply to achieve semi-global stabilization of the overall system.

\section{Introduction and Problem Statement}

Consider the problem of semi-globally
stabilizing a nonlinear system of the affine-in-control form
\be
\label{nonsys_semi}
\left\{\begin{array}{rcl}\dot x &=& f(x)+g(x)u,\cr
                           y &=& h(x),\end{array}\right.
\ee
where $x\in\RR^n$, $u\in\RR^m$ and $y\in\RR^p$ are the state, input and
output, respectively, and the mappings $f$, $g$ and $h$ are
smooth % in a neighborhood of a point $0$
with $f(0)=0$ and $h(0)=0$. In a semi-global stabilization
problem, we are to construct, for any
given, {\em arbitrarily large}, bounded set of the state
space $\calX_0$, a smooth feedback law,
say $u=v_{\calX_0}(x)$, with $v(0)=0$, such that the closed-loop
system is asymptotically stable at the origin with
$\calX_0$ contained in the domain of attraction.

The non-local stabilization of nonlinear systems of the form
(\ref{nonsys_semi}) has been made possible by the structural
decomposition, in the form of various normal forms, of these systems.

In Chapter 3, we propose an algorithm that
identifies a set of integers that are equivalent to the infinite
zero structure of linear systems and leads to
a normal form representation that corresponds to these integers as well
as to the system invertibility structure.
In the case that the system is square and invertible, {\em i.e.},
the system that was considered in \cite{isnc95,isnc99,scgn99},
$m=p=m_\rmd$, the normal form simplifies to
\be
\label{form1-square_semi}
\left\{\begin{array}{rcl}
\dot \eta &=& f_0(\eta,\xi),\cr
\dot\xi_{i,j} &=& \xi_{i,j+1}\!+\!\displaystyle\sum_{l=1}^{i-1}
            \delta_{i,j,l}(x)v_l,\; j=1,2,\cdots,q_i\!-\!1,\!\cr
\dot\xi_{i,q_i} &=& v_i,\cr
 y_i &=& \xi_{i,1}, \qquad i=1,2,\cdots,m,
\end{array}\right.
\ee
where $q_1\leq q_2\leq\cdots\leq q_m$,
$
\xi_i=\col\{\xi_{i,1},\xi_{i,2},\cdots,\xi_{i,q_i}\},
$
\newline
$
\xi=\col\{\xi_1,\xi_2,\cdots,\xi_m\},
$
and
\be\label{key2_semi}
\delta_{i,j,l}(x)=0,\quad \mbox{for } \; j<q_l, \; i=1,2,\cdots,m.
\ee

In this chapter, we would like to explore the application of
the normal form (\ref{form1-square_semi})-(\ref{key2_semi}) in solving the
problem of semi-global stabilization for nonlinear systems (\ref{nonsys_semi}).
The normal form (\ref{form1-square_semi})-(\ref{key2_semi}) does not
require a vector relative degree. The problem of semi-global
stabilization of system (\ref{nonsys_semi}) with a vector
relative degree has been well-studied in the literature.
For example, the work of \cite{byas91,suph91} solved the semi-global
stabilization problem for nonlinear systems with vector relative degrees,
{\em i.e.}, in the form of (\ref{vrd3}), but the globally asymptotically stable zero dynamics
\[
\dot\eta=f_0(\eta,\xi)
\]
is driven only
by
$
\xi_{i,1}, \;i=1,2,\cdots,m,
$
the states at the top of the $m$ chains of integrators.
The system with globally asymptotically stable zero dynamics is said to be of minimum phase.
The works of \cite{liss94,tess92} generalized this result of \cite{byas91,suph91}
by allowing $f_0$ to be dependent on any one state of each of the $m$
chains of integrators. More specifically, the system considered in
\cite{liss94,tess92} can be represented as follows,
\be\label{semi-global-sys-90}
\left\{\begin{array}{rcl}
\dot \eta &=& f_0(\eta,\xi_{1,\ell_1},\xi_{2,\ell_2},\cdots,\xi_{m,\ell_m}),\cr
\dot\xi_{i,j} &=& \xi_{i,j+1},\;\;
j=1,2,\cdots,r_i-1,\cr
\dot\xi_{i,q_i} &=& v_i,\cr
y_i &=& \xi_{i,1},\;\; i=1,2,\cdots,m,
\end{array}
\right.
\ee
where
$
1\leq \ell_i\leq r_i+1, \; i=1,2,\cdots,m,
$
and $\xi_{i,q_i+1}\equiv v_i$.
The peaking phenomenon, which was identified in \cite{suph91} as a
main obstacle to semi-global stabilization, in such systems is eliminated by stabilizing
part of linear subsystem with a high-gain linear control and the remaining part of the
linear subsystem with a small, bound nonlinear control \cite{tess92}.
The reference \cite{liss94}
shows that the same problem can be solved by linear state feedback laws, as those in \cite{tess92},
depend only on the linear states. The fundamental issue in the design of such linear
state feedback laws is to induce a specific time-scale structure in the linear part
of the closed-loop system.
This time-scale structure consists of a very slow and a very fast time scale, which are the
results of a linear state feedback of the high-and-low-gain nature.

Note that in \cite{byas91,liss94,suph91,tess92}, the system is considered to be minimum phase,
which means that its zero dynamic $\dot\eta=f(\eta,0,\cdots,0)$ have a globally asymptotically stable
equilibrium at the origin.

In this chapter, we consider the semi-global stabilization problem for
the following minimum phase nonlinear system,
\be\label{semi-global-sys}
\left\{\begin{array}{rcl}
\dot \eta &=& f_0(\eta,\xi_{1,\ell_1},\xi_{2,\ell_2},\cdots,\xi_{m,\ell_m}),\cr
\!\!\dot\xi_{i,j} &=& \xi_{i,j+1}\!+\!
\displaystyle\sum_{l=1}^{i-1}\delta_{i,j,l}(\eta,\xi)v_l,\;
j\!=\!1,2,\cdots,q_i\!-\!1,\!\!\!\!\cr
\!\!\!\!\!\dot\xi_{i,q_i} &=& v_i,\cr
y_{\rmd,i} &=& \xi_{i,1},\;\; i=1,2,\cdots,m,
\end{array}
\right.
\ee
where
$
q_1\leq q_2\leq\cdots\leq q_m,
$
$
\xi_i\!=\!\col\{\xi_{i,1},\xi_{i,2},\cdots,\xi_{i,q_i}\},
$
\newline
$
\xi\!=\!\col\{\xi_1,\xi_2,\cdots,\xi_m\},
$
%$
%v\!=\!\col\{v_1,v_2,\cdots,v_m\},
%$
and
\be\label{key33}
\delta_{i,j,l}(\eta,\xi)=0,\quad \mbox{for } \; j<q_l, \; i=1,2,\cdots,m.
\ee
\be\label{key4_semi}
\ell_1\leq q_1+1,\quad \ell_i\leq q_1, \quad i=2,3,\cdots,m,
\ee
with $\xi_{1,q_1+1}\equiv v_1$.
%Denote $n_0=\sum_{i=1}^m q_i$.
Note that the zero dynamics of (\ref{semi-global-sys})
is given by
\[
\dot\eta=f_0(\eta,0,\cdots,0).
\]

As explained earlier, no vector relative degree is required for systems to
be decomposed into the above normal form.

Note that in \cite{liss94,tess92}, $\delta_{i,j,l}=0$.
That is, the systems considered in \cite{liss94,tess92} are a cascade of
a linear subsystem with the zero dynamic, which is the only source of nonlinearity.

The remainder of this chapter is organized as follows.
Section~\ref{main_semi} presents our solution to the semi-global stabilization problem for
nonlinear systems without vector relative degrees. Design examples are presented to
illustrate how the proposed design approach works.
A brief conclusion to the chapter is drawn in Section~\ref{conclusion_semi}.

\section{Main Results}
\label{main_semi}

\begin{defi}\it
The system (\ref{semi-global-sys}) is semi-globally stabilizable by state feedback
if, for any compact set of initial conditions $\calX_0$ of the state space,
there exists a smooth state feedback
\be\label{controller1}
v=\alpha_{\calX_0}(\eta,\xi)
\ee
such that the equilibrium $(0,0)$ of the closed-loop system
(\ref{semi-global-sys}) and (\ref{controller1}) is locally asymptotically stable
and $\calX_0$ is contained in its domain of attraction.
\end{defi}

Theorem~\ref{stab_general} can be modified to deal with the semi-global stabilization problem.

\begin{theo}\label{lemma_semi_1}\it
Consider a system in the form
\be\label{mix_form_semi}
\!\!\!\!\left\{\begin{array}{rcl}
\dot \eta &=& f_0(\eta,\xi_{1,1},\xi_{2,1},\cdots,\xi_{m,1}),\cr
\!\!\dot\xi_{i,j} &=& \xi_{i,j+1}\!+\!
\displaystyle\sum_{l=1}^{i-1}\delta_{i,j,l}(\eta,\xi)v_l,\;\;
j=1,2,\cdots,q_i-1,\!\!\cr
\dot\xi_{i,q_i} &=& v_i,\cr
y_i &=& \xi_{i,1},\;\; i=1,2,\cdots,m,
\end{array}
\right.
\ee
where
$
q_1\leq q_2\leq\cdots\leq q_m,
$
$
\xi\!=\!\col\{\xi_1,\xi_2,\cdots,\xi_m\},
$
$
\xi_i\!=\!\col\{\xi_{i,1},\cdots,\xi_{i,q_i}\},
$
\be\label{key2_semi2}
\delta_{i,j,l}=0,\quad \mbox{for } \; j<q_l, \; i=1,2,\cdots,m.
\ee
and all functions
$\delta_{i,j,l}$ and $f_0$ are smooth.
Then the system is semi-globally stabilizable if
\ben
\item
for any compact set $\calZ\subset \RR^{n_0}$, there exist $c$,
$\phi_{i,1}(\eta)$,
$i=1,2,\cdots,m$,
and a smooth, positive definite Lyapunov function $W(\eta)$, such that
\[
\calZ\subset\{\eta: W(\eta)\leq c\},
\]
\[
\dot W=\frac{\partial W}{\partial \eta}f_0(\eta,\phi_{1,1}(\eta),\phi_{2,1}(\eta),\cdots,\phi_{m,1}(\eta))< 0,
\]
\[\hspace{40mm}
\quad \forall \eta\in \{\eta: W(\eta)\leq c\}\setminus\{0\};
\]
\item
there exists an ordered list $\kappa$ containing all variables of $\xi$
such that, for $j=1,2,\cdots,q_i-1$, $l=1,2,\cdots, i-1$, $i=1,2,\cdots,m$,
\ben
\item
$\xi_{i,j}$ appears earlier than $\xi_{i,j+1}$ in $\kappa$;
\item for $\delta_{i,j,l}\neq0$, the variables $\xi_l$ appear earlier than $\xi_{i,j}$
in $\kappa$;
\item $\delta_{i,j,l}$ depends only on $\eta$, $\xi_{\ell,1}$, $\ell=1,2,\cdots,m$, and
the variables that appear no later than $\xi_{i,j}$
in $\kappa$.
\een
\een
\end{theo}

In what follows, we will present an algorithm for constructing
a family of feedback laws that semi-globally stabilize the
system (\ref{semi-global-sys}). This algorithm consists of
two steps.

We first find positive constants $c_{i,k}$, such that the polynomials
\[
p_i(s)=s^{\ell_i-1}+c_{i,\ell_i-2}\;s^{\ell_i-2}+\cdots
+c_{i,1}\;s+c_{i,0},\quad
i=1,2,\cdots,m,
\]
have all roots with negative real parts. Define
\[
\xi_{i,\ell_i}^\star=-\varepsilon^{\ell_i-1}c_{i,0} \xi_{i,1}
-\varepsilon^{\ell_i-2}c_{i,1} \xi_{i,2}-\cdots
-\varepsilon c_{i,\ell_i-2} \xi_{i,\ell_i-1},\hspace{10mm}
\]
\be\label{low1}
\hspace{60mm} i=1,2,\cdots,m.
\ee
where $\varepsilon>0$.
Consider
\be\label{part1}
\left\{\begin{array}{rcl}
\dot \eta &=& f_0(\eta,\xi_{1,\ell_1}^\star,\xi_{2,\ell_2}^\star,\cdots,\xi_{m,\ell_m}^\star),\cr
\!\!\dot\xi_{i,j} &=& \xi_{i,j+1},\;\;
j=1,2,\cdots,\ell_i-2,\cr
\dot\xi_{i,\ell_i-1} &=& \xi_{i,\ell_i}^\star,\;\; i=1,2,\cdots,m.
\end{array}
\right.
\ee
Denote $z=\col\{\xi_{1,1},\xi_{1,2},\cdots,\xi_{\ell_1-1},\cdots,\cdots,\xi_{m,1},\xi_{m,2},\cdots,\xi_{\ell_m-1}\}$.

%all states $\xi_{i,j}$, $j=1,2,\cdots,\ell_i-1,\;\; i=1,2,\cdots,m$, as $z$.

Following \cite{liss94,tess92}, the dynamics of (\ref{part1}) with $\xi_{i,\ell_i}^\star$ given by
(\ref{low1}) has a locally
asymptotically stable equilibrium at the origin of
$
(\eta,z).
$
Moreover,
the domain of attraction of this equilibrium can be made arbitrarily
large by decreasing the value of the low gain parameter $\varepsilon$.

\begin{lemma}\label{low-gain-part}\it
Consider the system (\ref{part1}) with $\xi_{i,\ell_i}^\star$ given by (\ref{low1}).
Suppose that its zero dynamics $\dot\eta=f_0(\eta,0,\cdots,0)$ has a globally asymptotically
stable equilibrium at the origin.
%Then, for any $R>0$, there exists an $\varepsilon^\star>0$ such that,
%for any $\varepsilon\in(0, \varepsilon^\star]$, the system (\ref{part1}) is locally
%asymptotically stable and, moreover,
%\[
%\|(\eta(0),z(0))\|\leq R
%\Longrightarrow
%\left\{\begin{array}{l}
%\lim_{t\rightarrow\infty}\eta(t)=0,\cr
%\lim_{t\rightarrow\infty}z(t)=0.
%\end{array}\right.
%\]
Then for any compact set $\calY$,
there exist $\varepsilon$, $c$ and a smooth, positive definite, Lyapunov function $W(\eta,z)$, such that
\[
\calY \subset \{(\eta,z): W (\eta,z)\leq c\},
\]
and
\[
\dot W< 0 ,\quad \forall (\eta,z)\in \{(\eta,z): W (\eta,z)\leq c\}\setminus \{0\}.
\]
\end{lemma}

\medskip

Once the virtual inputs $\xi_{i,\ell_i}^\star,\;i=1,2,\cdots,m$, have been obtained,
both \cite{tess92} and \cite{liss94} design the overall controller by using linear high-gain
state feedback. This is possible because the systems considered there are
linear except the zero dynamics. In our situation, the system is in the form
of (\ref{semi-global-sys}). Because of the nonlinearities $\delta_{i,j,l}(\eta,\xi)v_l$,
we have to resort to backstepping procedure as described in Theorem~\ref{lemma_semi_1},
where a special case of (\ref{semi-global-sys}), {\em i.e.}, $\ell_i=1,\; i=1,2,\cdots,m$, is considered.

Consider the dynamics of (\ref{part1}) with $\xi_{i,\ell_i}^\star$ given by
(\ref{low1}) as the zero dynamics of the system (\ref{mix_form_semi}), by Theorem~\ref{lemma_semi_1} and Lemma~\ref{low-gain-part},
we have

\begin{theo}\label{semi-theo}\it
Consider the system (\ref{semi-global-sys}) with (\ref{key33}) and (\ref{key4_semi}).
Assume that all functions $\delta_{i,j,l}$ and $f_0$ are smooth.
If
\ben
\item
the zero dynamics
$
\dot\eta=f_0(\eta,0,0,\cdots,0)
$
is globally asymptotically stable
at the equilibrium $\eta=0$;
\item
there exists an ordered list $\kappa$ containing all variables of $\xi$
with $\xi_{s,p}$, $p=1,2,\cdots,\ell_s-1$, $s=1,2,\cdots,m$, being its first
$\sum_{s=1}^m (\ell_s-1)$ variables,
such that, for $j=1,2,\cdots,q_i-1$, $l=1,2,\cdots, i-1$, $i=1,2,\cdots,m$,
\ben
\item
$\xi_{i,j}$ appears earlier than $\xi_{i,j+1}$ in $\kappa$;
\item
the variables $\xi_l$ appear earlier than $\xi_{i,j}$
in $\kappa$ if $\delta_{i,j,l}\neq0$;
\item $\delta_{i,j,l}$ depends only on $\eta$, $\xi_{\ell,1}$, $\ell=1,2,\cdots,m$, and
the variables that appear no later than $\xi_{i,j}$
in $\kappa$.
\een
\een
Then the system is semi-globally stabilizable.
That is, for any compact set $\calX_0$ of the state space of $(\eta,\xi)$,
there exists a state feedback $v$ that locally asymptotically stabilizes
the system with $\calX_0$ contained in the domain of attraction.
\end{theo}

\begin{example}\rm
Consider a
three inputs three outputs
system in the form of
(\ref{semi-global-sys}) with three chains of integrators of lengths $\{3,4,4\}$,
and
$
\ell_1=2,\; \ell_2=1,\; \ell_3=3.
$
\be\label{example344_semi}
\left\{\begin{array}{rcl}
\dot\eta&=&f_0(\eta,\xi_{1,2},\xi_{2,1},\xi_{3,3}),\cr
\dot\xi_{1,j}&=&\xi_{1,j+1},\cr
\dot\xi_{1,3}&=&v_1,\cr
\dot\xi_{2,j}&=&\xi_{2,j+1},\cr
\dot\xi_{2,3}&=&\xi_{2,4}\!+\!\delta_{2,3,1}(
\eta;\xi_1;\xi_{2,1},\xi_{2,2},\xi_{2,3},\xi_{3,1},\xi_{3,2})v_1,\cr
\dot\xi_{2,4}&=&v_2,\cr
\dot\xi_{3,j}&=&\xi_{3,j+1},\quad j=1,2,\cr
\dot\xi_{3,3}&=&\xi_{3,4}+\delta_{3,3,1}(\eta;
\xi_1;\xi_{2,1},\xi_{2,2},\xi_{2,3},\xi_{3,1},\xi_{3,2},\xi_{3,3})v_1,\cr
\dot\xi_{3,4}&=&v_3.
\end{array}\!\!\right.
\ee
Suppose its zero dynamics
$
\dot\eta=f_0(\eta,0,0,0)
$
is globally asymptotically stable at the origin.
Clearly, this system satisfies the conditions in Theorem~\ref{semi-theo}
with
\[
\kappa=\{\xi_{1,1},\xi_{3,1},\xi_{3,2};\xi_{2,1},\xi_{1,2},\xi_{2,2},\xi_{1,3},\xi_{2,3},\xi_{3,3},\xi_{2,4},\xi_{3,4}\}.
\]

Consider the subsystem
\[
\left\{\begin{array}{rcl}
\dot\eta&=&f_0(\eta,\xi_{1,2}^\star,\xi_{2,1}^\star,\xi_{3,3}^\star),\cr
\dot\xi_{1,1}&=&\xi_{1,2}^\star,\cr
\dot\xi_{3,1}&=&\xi_{3,2},\cr
\dot\xi_{3,2}&=&\xi_{3,3}^\star,\cr
\end{array}\right.
\]
with
$
\xi_{1,2}^\star=-\varepsilon \xi_{1,1},\;
$
$
\xi_{2,1}^\star=0,\;
$
$
\xi_{3,3}^\star=-\varepsilon^2 \xi_{3,1}-\varepsilon \xi_{3,2}.
$
It is semi-globally asymptotically stable.

In what follows, we will illustrate how to implement
the level-by-level backstepping to find $v_1$, $v_2$ and $v_3$.
To carry out the backstepping on the first level variable $\xi_{2,1}$ to $\xi_{2,2}$,
we consider the following subsystem,
\[
\left\{\begin{array}{rcl}
\dot\eta&=&f_0(\eta,\xi_{1,2}^\star,\xi_{2,1},\xi_{3,3}^\star),\cr
\dot\xi_{1,1}&=&\xi_{1,2}^\star,\cr
\dot\xi_{2,1}&=&\xi_{2,2},\cr
\dot\xi_{3,1}&=&\xi_{3,2},\cr
\dot\xi_{3,2}&=&\xi_{3,3}^\star.\cr
\end{array}\right.
\]
with $\xi_{2,2}$ as the virtual input.
By backstepping, the desired input is given as
\[
\xi_{2,2}^\star=\phi_{2,2}(\eta,\xi_{1,1},\xi_{2,1},\xi_{3,1},\xi_{3,2}).
\]
Now consider backstepping from the second level variables.
To backstep from $\xi_{1,2}$ to $\xi_{1,3}$, we consider
\[
\left\{\begin{array}{rcl}
\dot\eta&=&f_0(\eta,\xi_{1,2},\xi_{2,1},\xi_{3,3}^\star),\cr
\dot\xi_{1,1}&=&\xi_{1,2},\cr
\dot\xi_{1,2}&=&\xi_{1,3},\cr
\dot\xi_{2,1}&=&\xi_{2,2}^\star,\cr
\dot\xi_{3,1}&=&\xi_{3,2},\cr
\dot\xi_{3,2}&=&\xi_{3,3}^\star,
\end{array}\right.
\]
with $\xi_{1,3}$ as the virtual input. The desired input is given as
\[
\xi_{1,3}^\star=\phi_{1,3}(\eta,\xi_{1,1},\xi_{1,2},\xi_{2,1},\xi_{3,1},\xi_{3,2}).
\]
To backstep from $\xi_{2,2}$ to $\xi_{2,3}$, we view $\xi_{2,3}$ as the virtual input of
\[
\left\{\begin{array}{rcl}
\dot\eta&=&f_0(\eta,\xi_{1,2},\xi_{2,1},\xi_{3,3}^\star),\cr
\dot\xi_{1,1}&=&\xi_{1,2},\cr
\dot\xi_{1,2}&=&\xi_{1,3}^\star,\cr
\dot\xi_{2,1}&=&\xi_{2,2},\cr
\dot\xi_{2,2}&=&\xi_{2,3},\cr
\dot\xi_{3,1}&=&\xi_{3,2},\cr
\dot\xi_{3,2}&=&\xi_{3,3}^\star.\cr
\end{array}\right.
\]
And the desired input is given as
\[
\xi_{2,3}^\star=\phi_{2,3}(\eta,\xi_{1,1},\xi_{1,2},\xi_{2,1},\xi_{2,2},\xi_{3,1},\xi_{3,2}).
\]
Next consider backstepping from the third level variables.
To backstep from $\xi_{1,3}$ to $v_1$ in the subsystem,
\[
\left\{\begin{array}{rcl}
\dot\eta&=&f_0(\eta,\xi_{1,2},\xi_{2,1},\xi_{3,3}^\star),\cr
\dot\xi_{1,1}&=&\xi_{1,2},\cr
\dot\xi_{1,2}&=&\xi_{1,3},\cr
\dot\xi_{1,3}&=&v_1,\cr
\dot\xi_{2,1}&=&\xi_{2,2},\cr
\dot\xi_{2,2}&=&\xi_{2,3}^\star,\cr
\dot\xi_{3,1}&=&\xi_{3,2},\cr
\dot\xi_{3,2}&=&\xi_{3,3}^\star,
\end{array}\right.
\]
we get
\[
v_1=v_1(\eta,\xi_1,\xi_{2,1},\xi_{2,2},\xi_{3,1},\xi_{3,2}).
\]
To backstep from $\xi_{2,3}$ to $\xi_{2,4}$, we consider
\[
\left\{\begin{array}{rcl}
\dot\eta&=&f_0(\eta,\xi_{1,2},\xi_{2,1},\xi_{3,3}^\star),\cr
\dot\xi_{1,1}&=&\xi_{1,2},\cr
\dot\xi_{1,2}&=&\xi_{1,3},\cr
\dot\xi_{1,3}&=&v_1,\cr
\dot\xi_{2,1}&=&\xi_{2,2},\cr
\dot\xi_{2,2}&=&\xi_{2,3},\cr
\dot\xi_{2,3}&=&\xi_{2,4}\!+\!\delta_{2,3,1}(
\eta;\xi_1;\xi_{2,1},\xi_{2,2},\xi_{2,3},\xi_{3,1},\xi_{3,2})v_1,\cr
\dot\xi_{3,1}&=&\xi_{3,2},\cr
\dot\xi_{3,2}&=&\xi_{3,3}^\star,
\end{array}\right.
\]
and obtain
\[
\xi_{2,4}^\star=\phi_{2,4}(\eta,\xi_1,\xi_{2,1},\xi_{2,2},\xi_{2,3},\xi_{3,1},\xi_{3,2}).
\]
Similarly, to backstep from $\xi_{3,3}$ to $\xi_{3,4}$, we obtain
\[
\xi_{3,4}^\star=\phi_{3,4}(\eta,\xi_1,\xi_{2,1},\xi_{2,2},\xi_{2,3},\xi_{3,1},\xi_{3,2},\xi_{3,3}).
\]
Finally, backstepping from the fourth level variables, we obtain
\[
\left.\begin{array}{rcl}
v_2&=&v_2(\eta,\xi_1,\xi_2,\xi_{3,1},\xi_{3,2},\xi_{3,3}),\cr
v_3&=&v_3(\eta,\xi_1,\xi_2,\xi_3.
\end{array}\right.
\]
The inputs $v_1$, $v_2$ and $v_3$
%as given by (\ref{v1}) and (\ref{v3})
semi-globally asymptotically
stabilize the origin of the system
(\ref{example344_semi}).
\end{example}

In what follows, we give an example
which requires the mixed chain-by-chain and level-by-level
backstepping design procedure.

\begin{example}\rm
Consider a %three input three output
system in the form of
(\ref{semi-global-sys}) with three chains of integrators of lengths $\{2,4,4\}$,
\be\label{equ_delta4}
\left\{\begin{array}{rcll}
\dot\eta &=& f_0(\eta,\xi_{1,2},\xi_{2,2},\xi_{3,2}),\cr
\dot\xi_{1,1} &=& \xi_{1,2},\cr
\dot\xi_{1,2} &=& v_1,\cr
\dot\xi_{2,1} &=& \xi_{2,2},\cr
\dot\xi_{2,2} &=& \xi_{2,3}+\delta_{2,2,1}(\eta,\xi_1,\xi_{2,1},
     \xi_{2,2},\xi_{3,1})v_1,\cr
\dot\xi_{2,3} &=& \xi_{2,4}+\delta_{2,3,1}(\eta,\xi_1,\xi_{2,1},
     \xi_{2,2},\xi_{2,3},\xi_{3,1},\xi_{3,2})v_1,\cr
\dot\xi_{2,4} &=& v_2,\cr
\dot\xi_{3,1} &=& \xi_{3,2},\cr
\dot\xi_{3,2} &=& \xi_{3,3}+\delta_{3,2,1}(\eta,\xi_1,\xi_{2,1},
     \xi_{2,2},\xi_{3,1},\xi_{3,2})v_1,\cr
\dot\xi_{3,3} &=& \xi_{3,4}+\delta_{3,3,1}(\eta,\xi_1,\xi_2,\xi_{3,1},
    \xi_{3,2},\xi_{3,3})v_1,\cr
\dot\xi_{3,4} &=& v_3,
\end{array}\right.
\ee
where
$
\ell_1=\ell_2=\ell_3=2.
$
Suppose its zero dynamics
$
\dot\eta=f_0(\eta,0,0,0)
$
is globally asymptotically stable at the origin.
It is obvious that the system satisfies the conditions in Theorem~\ref{semi-theo}
with
\[
\kappa=\{\xi_{1,1},\xi_{2,1},\xi_{3,1};\xi_{1,2},\xi_{2,2},
\xi_{3,2},\xi_{2,3},\xi_{2,4},\xi_{3,3},\xi_{3,4}\}.
\]
We first find the low-gain feedback,
\[
\xi_{1,2}^\star=-\varepsilon \xi_{1,1},\;\;
\xi_{2,2}^\star=-\varepsilon \xi_{2,1},\;\;
\xi_{3,2}^\star=-\varepsilon \xi_{3,1}.
\]
Then we carry out a mixed chain-by-chain and level-by-level backstepping
in the
order of
$\xi_{1,2},\xi_{2,2},
\xi_{3,2},\xi_{2,3},\xi_{2,4},$ $\xi_{3,3},\xi_{3,4}$
to obtain
\[
\left.\begin{array}{rcl}
v_1&=&v_1(\eta,\xi_1,\xi_{2,1},\xi_{3,1}),\cr
v_2&=&v_2(\eta,\xi_1,\xi_2,\xi_{3,1},\xi_{3,2}),\cr
v_3&=&v_3(\eta,\xi_1,\xi_2,\xi_3).
\end{array}\right.
\]
\end{example}

\medskip
\medskip

\begin{example}\rm
Consider
\[
\left\{\begin{array}{rcl}
\dot\eta &=& -\eta+(v_1+\xi_{2,2})\sin\eta,\cr
\dot\xi_{1,1} &=& \xi_{1,2},\cr
\dot\xi_{1,2} &=& v_1,\cr
\dot\xi_{2,1} &=& \xi_{2,2},\cr
\dot\xi_{2,2} &=& \xi_{2,3}+\eta v_1,\cr
\dot\xi_{2,3} &=& v_2.
\end{array}\right.
\]
Obviously the system satisfies the conditions in Theorem~\ref{semi-theo}
with
\[
q_1=2,\quad q_2=3,\quad \ell_1=3,\quad \ell_2=2.
\]
Choosing all poles of the linear slow subsystems to be $-\varepsilon$,
we obtain
\[
\left.\begin{array}{rcl}
v_1&=&- \varepsilon^2 \xi_{1,1}-2 \varepsilon \xi_{1,2},\cr
\xi_{2,2}^\star&=&-\varepsilon \xi_{2,1}.
\end{array}\right.
\]
Next, we view $\xi_{2,3}$ as a virtual input. By backstepping, the desired input is given as follows,
\[
\xi_{2,3}^\star=-(\varepsilon+1)\xi_{2,1}-(\varepsilon+1)\xi_{2,2}-\eta(\sin \eta - \varepsilon^2 \xi_{1,1}-2 \varepsilon \xi_{1,2}).
\]
Finally, backstepping one more time, we get
\[
\left.\begin{array}{l}
v_2=-(2\varepsilon+1)\xi_{2,1}-(2\varepsilon+3)\xi_{2,2}-(\varepsilon+2)\xi_{2,3}
-\eta(\sin\eta+2v_1-\varepsilon v_1-\varepsilon^2\xi_{1,2}) \vspace{5pt} \cr
\qquad-[-\eta+(v_1+\xi_{2,2})\sin\eta](\sin\eta+\eta\cos\eta+v_1).
\end{array}\right.
\]
The Lyapunov function is given as
\[
\left.\begin{array}{l}
V=\Bigl\{\eta^2+\varepsilon^2\xi_{1,1}^2+(\varepsilon\xi_{1,1}+\xi_{1,2})^2+\xi_{2,1}^2+(\xi_{2,2}+\varepsilon\xi_{2,1})^2\cr
\qquad+[\xi_{2,3}+(\varepsilon+1)\xi_{2,1}+(\varepsilon+1)\xi_{2,2}+\eta(\sin \eta - \varepsilon^2 \xi_{1,1}-2 \varepsilon \xi_{1,2})]^2\Bigl\}/2,
\end{array}\right.
\]
and
\[
\left.\begin{array}{l}
\dot V=\eta[-\eta+(- \varepsilon^2 \xi_{1,1}-2 \varepsilon \xi_{1,2} -\varepsilon \xi_{2,1})\sin\eta]
-\varepsilon^3\xi_{1,1}^2\cr
\hspace{15mm}-\varepsilon^2\xi_{1,1}\xi_{1,2} -\varepsilon\xi_{1,2}^2-\varepsilon\xi_{2,1}^2
  -(\xi_{2,2}+\varepsilon\xi_{2,1})^2\cr
\hspace{20mm}  -[\xi_{2,3}+(\varepsilon+1)\xi_{2,1}+(\varepsilon+1)\xi_{2,2}
+\eta(\sin \eta - \varepsilon^2 \xi_{1,1}-2 \varepsilon \xi_{1,2})]^2.
\end{array}\right.
\]
%Define
%\[
%{\cal D}=\bigl\{(\eta,\xi): \;|\xi_{1,1}|< r,\; |\xi_{1,2}|< r,\; |\xi_{2,1}|< r\bigl\}.
%\]
%For any compact set of initial conditons $\calX_0\subset\RR^6$, let $c$ and $r$ be such that
%\[
%\calX_0\subset\Bigl\{(\eta,\xi)\in\RR^6 : \; V(\eta,\xi)\leq c\Bigl\}\subset {\cal D}.
%\]
%By selecting
%\[
%0<\varepsilon<(\sqrt{9+4/r}-3)/2,
%\]
%we have
%\[
%\dot V< 0, \quad \forall(\eta,\xi)\in \Bigl\{(\eta,\xi)\in\RR^6 : \; V(\eta,\xi)\leq c\Bigl\}\setminus \{0\}.
%\]
%Therefore, the equilibrium $(\eta,\xi)=0$ of the closed-loop system is locally asymptotically stable, and $\calX_0$ is contained in its domain of attraction.
Shown in Figs.~\ref{f1a} and \ref{f1b} are
state trajectories of the closed-loop system with different initial conditions.

\begin{figure}[ht!]
\centerline{\includegraphics[width=100mm]{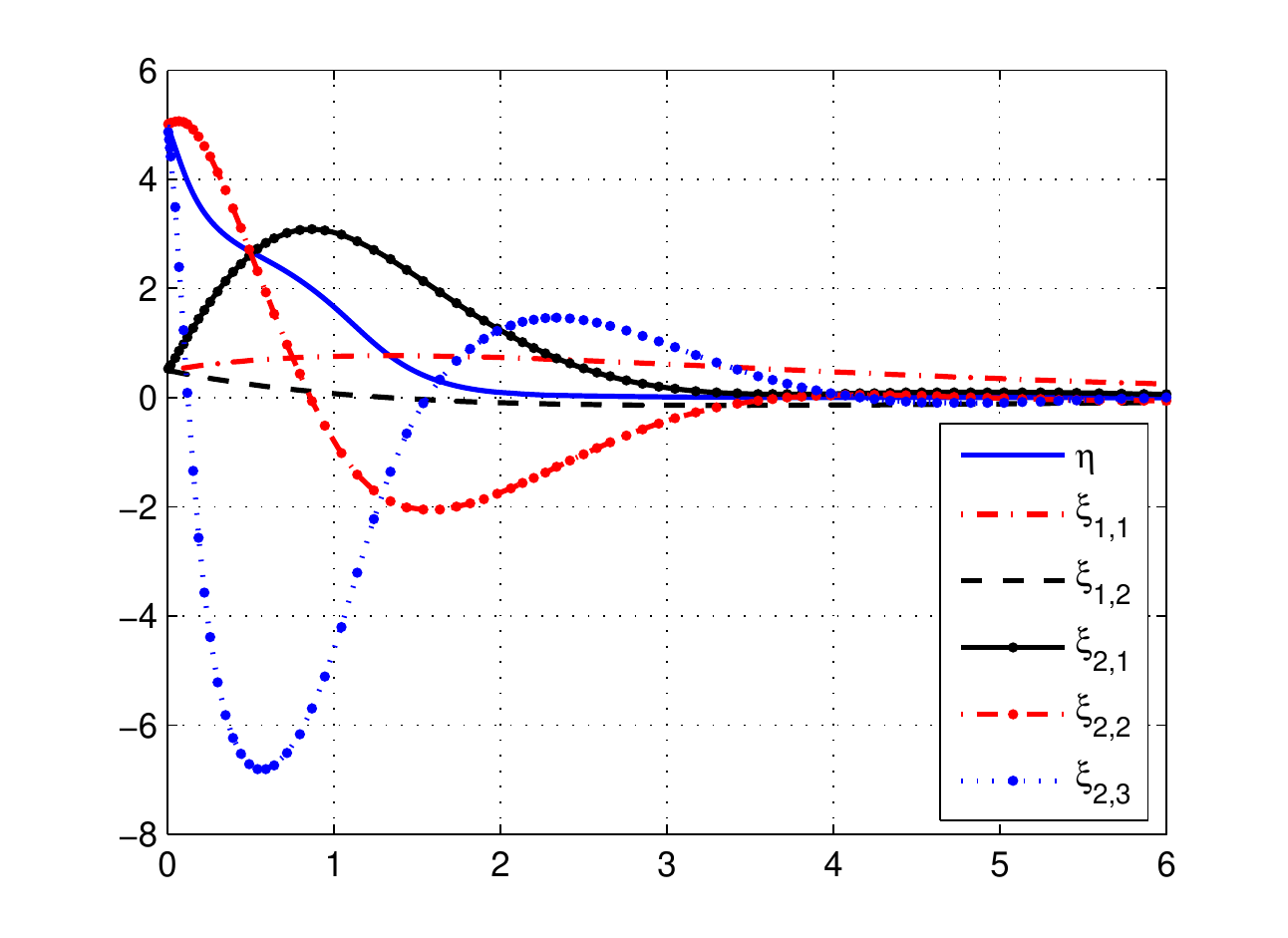}}
\vspace{-5mm}
\caption{State trajectories with the initial condition $(5, 0.5, 0.5, 0.5, 5, 5)$ and $\varepsilon=0.5$.}
\label{f1a}
\centerline{\includegraphics[width=100mm]{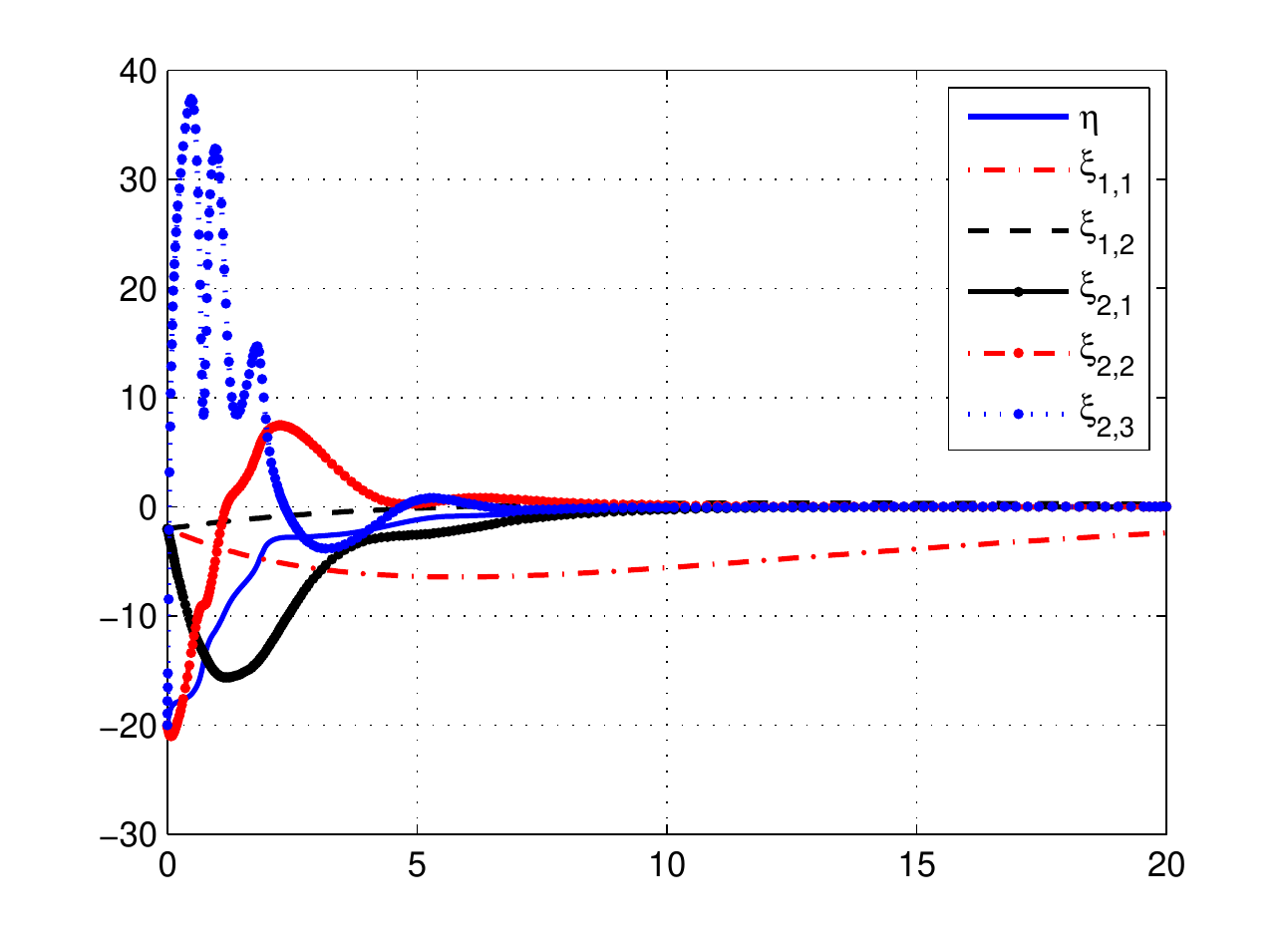}}
\vspace{-5mm}
\caption{State trajectories with the initial condition $(-20,-2,-2,-2, -20, $ $-20)$ and $\varepsilon=0.15$.}
\label{f1b}
\end{figure}

\end{example}

\section{Summary of the Chapter}
\label{conclusion_semi}

In this chapter, we showed how the
structural decomposition in Chapter 3 can be used to solve the
semi-global stabilization of a class of multiple input multiple output systems
without vector relative degrees. The design procedure
involved several existing design techniques in
nonlinear stabilization, including low gain
feedback and different forms of backstepping design procedures in Chapter 4.

\chapter{Disturbance Attenuation for Nonlinear Systems}

The problems of disturbance attenuation and almost disturbance decoupling
play a central role in control theory. In this chapter, by
employing the structural decomposition of multiple
input multiple output nonlinear systems in Chapter 3 and the backstepping
procedures
%that were developed based on the structure of the decomposed system
in Chapter 4, we show that these two problems can be
solved for a larger class of nonlinear systems.

\section{Introduction and Problem Statement}

Consider the problems of disturbance attenuation and almost disturbance
decoupling with internal stability for nonlinear systems affine in control,
\be
\label{nonsys_dist}
\left\{\begin{array}{rcl}\dot x &=& f(x)+g(x)u+p(x,w),\cr
                           y &=& h(x),\end{array}\right.
\ee
where $x\in\RR^n$, $u\in\RR^m$, $y\in\RR^m$ and $w\in\RR$ are the state, input,
output and disturbance, respectively, and the mappings $f$, $g$, $p$ and $h$ are
smooth with $f(0)=0$ and $h(0)=0$.
The problem of almost disturbance decoupling was originally formulated and solved
in \cite{willems1981ais} for linear systems and was later extended to single input single
output (SISO) minimum phase nonlinear systems in \cite{isidori1996nad,marino1989add,marino1994nha}.
It was further extended to SISO non-minimum phase systems in
\cite{isidori1996gad,lin1998almost}.

The problem of almost disturbance decoupling is, for any {\em a priori} given arbitrarily small
scalar $\gamma>0$, to find a feedback law such that the $\calL_2$ gain from
the disturbance to the output is less than or equal to $\gamma$. A practical
solution to the almost disturbance decoupling problem would require the
resulting closed-loop system to be globally or locally asymptotically stable as well.
Here in this chapter, we will focus on the requirement of global asymptotic
stability. The problem of disturbance attenuation is a less stringent one
in that it does not require the bound on the resulting $\calL_2$ to be arbitrarily
small. The problem of almost disturbance decoupling
is a special case of disturbance attenuation.

The problem of disturbance attenuation (or almost disturbance decoupling) with
stability can be solved by establishing the dissipativity of
the system \cite{isnc99}. That is, the problem of disturbance
attenuation with stability (or almost disturbance decoupling)
for a given system is, for a given (arbitrarily small)
scalar $\gamma>0$, to find a feedback law $u=u(x)$ such that the
resulting closed-loop system is strictly dissipative with respect
to the supply rate
$
q(w,y)=\gamma^2w^2-y^2,
$
which is equivalent to finding a feedback law $u=u(x)$ such that,
for some smooth, positive definite and proper function $V(x)$,
%{\em i.e.},
%\be
%\label{dissipation}
%\underline{\alpha}(\|x\|)\leq V(x)\leq \overline\alpha(\|x\|),\;x\in\RR^n
%\ee
%for some class $\calK_\infty$ functions $\underline{\alpha}$ and $\overline\alpha$,
the dissipation inequality
\begin{eqnarray}
\label{dissipation}
\frac{\partial V}{\partial x}\left(f(x)\!+\!g(x)u(x)\!+\!p(x,w)\right)
\leq -\alpha(\|x\|)\!+\!\gamma^2w^2\!-\!h^2(x),\quad \nonumber \\ \;x\in\RR^n,w\in\RR,
\end{eqnarray}
holds for some class $\calK_\infty$ function $\alpha$.

The inequality (\ref{dissipation}) guarantees that
the response of the closed-loop system in the absence of disturbance is
globally asymptotically stable and, with $x(0)=0$,
\[
\int_0^\infty y^2(t)dt\leq \gamma^2 \int_0^\infty w^2(t)dt,
\]
for every $\calL_2$ disturbance $w$.

The solution to the problem of disturbance attenuation and
almost disturbance decoupling usually resorts to transforming the
nonlinear systems into certain structural normal forms. For example,
in \cite{isidori1996nad,marino1989add,marino1994nha}, the problem of almost
disturbance decoupling problem with stability was solved for systems in
the following normal form
\be
\label{mmi}
\left\{\begin{array}{rcl}
\dot z &=& f_0(z,\xi_1)+p_0(z,\xi_1)w,\cr
\dot\xi_i &=& \xi_{i+1}+p_i(z,\xi_1,\xi_2,\cdots,\xi_i)w,\quad i=1,2,\cdots,r-1,\cr
\dot\xi_r &=& u+p_r(z,\xi_1,\xi_2,\cdots,\xi_r)w,\cr
y &=& \xi_1,
\end{array}\right.
\ee
where $f_0(0,0)=0$. A critical assumption made there is that
the system is of minimum phase, that is, the zero dynamics
$\dot z=f_0(z,0)$ is globally asymptotically stable.

The work \cite{isidori1996gad} relaxes the minimum phase assumption by
allowing part of the zero dynamics to be unstable as long as it is
unaffected by the disturbance and is stabilizable though the output of the system.
That is, it is assumed that the dynamic of $z$ in (\ref{mmi}) takes the following
form,
\be\label{two_z}
\left\{\begin{array}{rcl}
\dot z_1 &=& f_1(z_1,z_2,\xi_1)+p_0(z_1,z_2,\xi_1)w,\cr
\dot z_2 &=& f_2(z_2,\xi_1),
\end{array}\right.
\ee
where the $z_1$ subsystem is globally asymptotically stable at $z_1=0$ and
there exists some smooth $v_2(z)$ such that $\dot z=f_2(z_2,v(z_2))$ is globally
asymptotically stable at $z_2=0$. In a further note \cite{lin1998almost}, it was
pointed out that, under some further structural assumption on the $z_2$ subsystem,
the problem of almost disturbance decoupling with stability can
still be solved even if the $z_2$ subsystem is affected by the disturbance.

In an effort to solve the problem of disturbance attenuation for
multiple input multiple output nonlinear systems, a normal form
for square invertible systems was developed in \cite{isnc95,isnc99,scgn99}.

In Chapter 3, we studied the
structural properties of affine-in-control nonlinear systems beyond
the case of square invertible systems. We proposed an algorithm that
identifies a set of integers that are equivalent to the infinite
zero structure of linear systems and leads to
a normal form representation that corresponds to these integers as well
as to the system invertibility structure.

This new normal form facilitates the control
design. As shown in Chapter 4, it allows
the development of some new backstepping design procedures, the
level-by-level backstepping and the mixed chain-by-chain and
level-by-level backstepping. These new backstepping procedures
lead to the stabilization of a larger class of
systems that the conventional chain-by-chain backstepping design procedure
cannot stabilize.
The objective of this chapter is to show that the backstepping
design procedures of Chapter 4 can also be utilized to
solve the problems of disturbance attention and almost disturbance
decoupling for a larger class of multiple input multiple output
systems.

The remainder of this chapter is organized as follows. In Section \ref{preliminary_dist},
we recall some results on the problems of disturbance attenuation and almost disturbance decoupling for SISO systems.
We will also describe
the level-by-level  and the mixed chain-by-chain and
level-by-level backstepping design procedures.
Section~\ref{main_dist} presents our solutions to the problems of
disturbance attenuation and almost disturbance decoupling.
A brief conclusion to the chapter is drawn in Section~\ref{conclusion_dist}.

\section{Preliminary Results}
\label{preliminary_dist}

We first recall the follow result on disturbance attenuation with stability from \cite{isnc99}.
This result will serve as a building block in our design procedures.

\begin{lemma}\label{backsteping_w}\em
Consider a system described by
\be
\left\{\begin{array}{rcl}\dot z &=& f_0(z,\xi)+p_0(z,w),\cr
                         \dot \xi &=& u+f_1(z,\xi)+p_1(z,\xi,w),\cr
                           y &=& h(z,\xi),\end{array}\right.
\ee
where $(z,\xi)\in\RR^n \times \RR$, $f_0(0,0)=0$ and $f_1(0,0)=0$.
Assume that
\begin{eqnarray*}
\|p_0(z,w)\| &\leq & R_0(z) |w|,\;\forall z,w,\\
|p_1(z,\xi,w)| &\leq& R_1(z,\xi) |w|,\;\forall z,\xi,w,
\end{eqnarray*}
for some smooth real-valued functions $R_0(z)$ and $R_1(z,\xi)$.
Suppose that there exist a number $\gamma>0$, a smooth real-valued function $v(z)$ with $v(0)=0$,
a smooth positive definite and radially unbounded function $V(z)$, and a class $\calK_\infty$
function $\alpha_0(\cdot)$ such that
\begin{eqnarray}
\label{main_condition_1}
\frac{\partial V}{\partial z} [f_0(z,v(z))+p_0(z,w)] \leq -\alpha_0(\|z\|)+\gamma^2w^2-h^2(z,v(z)), \quad \nonumber \\
\forall z,\xi,w,\quad
\end{eqnarray}
that is, there exists a smooth $v(z)$ such that the subsystem
\[
\left\{\begin{array}{rcl}\dot z &=& f_0(z,v(z))+p_0(z,w),\cr
                           y &=& h(z,v(z)),\end{array}\right.
\]
is strictly dissipative with respect to the supply rate
$ q(w,y)=\gamma^2w^2-y^2$.

Then, for every $\epsilon>0$, there exist a smooth feedback law $u=u(z,\xi)$,
a smooth positive definite and radially unbounded function $W(z,\xi)$,
and a class $\calK_\infty$ function $\alpha(\cdot)$ such that
\begin{eqnarray*}
&& \frac{\partial W}{\partial z} \left(f_0(z,\xi)+p_0(z,w)\right)
 + \frac{\partial W}{\partial \xi} \left(u(z,\xi)+f_1(z,\xi)+p_1(z,\xi,w)\right)\qquad\qquad\quad\\
&&\qquad \leq -\alpha(\|\col\{z,\xi\}\|)\!+\!(\gamma+\epsilon)^2w^2 \!-\! h^2(z,\xi),\;
\forall z,\xi,w,
\end{eqnarray*}
or equivalently, there exist a smooth feedback law $u=u(z,\xi)$
such that the resulting closed-loop system is strictly dissipative with respect
to the supply rate $q(w,y)=(\gamma+\epsilon)^2w^2-y^2$.
\end{lemma}

In \cite{isnc99}, the possibility of fulfilling the main condition (\ref{main_condition_1}) in Lemma~\ref{backsteping_w} is discussed.
In the context of the almost disturbance decoupling problem,
suppose that the $z$-subsystem
\be
\label{almost_z_1a}
\dot z=f(z,\xi_1,w)
\ee
can be decomposed as
\be\label{almost_z_1}
\left\{\begin{array}{rcl}\dot z_1 &=& f_1(z_1,z_2,\xi_1,w),\cr
                           z_2 &=& f_2(z_2,\xi_1),\end{array}\right.
\ee
where $z_1$ represents ``stable component" and $z_2$ represents ``unstable but stabilizable component."

\begin{lemma}\em
\label{z1z2}
Consider system (\ref{almost_z_1a}) which can be decomposed as (\ref{almost_z_1}). Suppose that
\ben
\item
there exists a smooth positive definite and radially unbounded function $V_1(z_1)$ such that
\begin{eqnarray*}
\frac{\partial V_1}{\partial z_1} f_1(z_1,z_2,\xi_1,w)
\leq -\alpha_1(\|z_1\|)+\gamma_0^2 w^2+\gamma_0^2 \|z_2\|^2+\gamma_0^2 \xi_1^2,
\end{eqnarray*}
for some $\calK_\infty$ function $\alpha_1$ and some $\gamma_0>0$,
\item
there exist a smooth real-valued function $v_2(z_2)$ with $v_2(0)=0$,
and a smooth positive definite and radially unbounded function $V_2(z_2)$ such that
\[
\frac{\partial V_2}{\partial z_2} f_2(z_2,v_2(z_2))+v_2^2(z_2)
\leq -\alpha_2(\|z_2\|),
\]
for some $\calK_\infty$ function $\alpha_2$.
\een

Then, for every $\gamma>0$, there exist a smooth $v(z)$ with $v(0)=0$, and a smooth positive definite and radially unbounded function $V(z)$
such that
\[
\frac{\partial V}{\partial z} f(z,v(z),w)
\leq -\alpha(\|z\|)+\gamma^2 w^2-v^2(z),
\]
for some $\calK_\infty$ function $\alpha(\cdot)$.
\end{lemma}

\section{Disturbance Attenuation and Almost Disturbance Decoupling with Stability}
\label{main_dist}

Suppose that, by the algorithm in \cite{lisn08}, system (\ref{nonsys_dist}) is transferred into the following form,
\be\label{almost_dd2}
\left\{\begin{array}{rcl}
\dot z &=& f_0(z,\xi)+p_0(z,w),\cr
\dot\xi_{i,j} &=& \xi_{i,j+1}\!+\!
\displaystyle\sum_{l=1}^{i-1}\delta_{i,j,l}(z,\xi)v_l+p_{i,j}(z,\xi,w),\quad j=1,2,\cdots,q_i-1,\cr
\dot\xi_{i,q_i} &=& v_i+p_{i,q_i}(z,\xi)w,\cr
y_i &=& \xi_{i,1},\;\; i=1,2,\cdots,m,
\end{array}
\right.
\ee
where $\xi\!=\!\col\{\xi_1,\xi_2,\cdots,\xi_m\}$,
$\xi_i\!=\!\col\{\xi_{i,1},\xi_{i,2},\cdots,\xi_{i,q_i}\}$,
$q_1\leq q_2\leq\cdots\leq q_m$,
and functions
$f_0$, $p_0$ and $p_{i,j}$, $j=1,2,$ $\cdots,q_i, \;i=1,2,\cdots,m$ are smooth with
$f_0(0,0,\cdots,0)=0$.
Moreover,
\[
\delta_{i,j,l}=0,\quad \mbox{for } \; j<q_l, \; i=1,2,\cdots,m.
\]

As in the literature on the problems of disturbance attention and almost disturbance
decoupling for SISO systems, we assume that the zero dynamics is driven only
by the states on the top of the $m$ chains of integrators,
$\xi_{i,1}, \;\;i=1,2,\cdots,m$.
That is
$
\dot z = f_0(z,\xi_{1,1},\xi_{2,1},\cdots,\xi_{m,1})+p_0(z,w).
$

To apply the level-by-level backstepping,
we also assume that the coefficients $\delta_{i,j,l}$, $p_{i,j}(z,\xi,w)$, $j=1,2,\cdots,q_i$, $i=1,2,\cdots,m$, satisfy
the level-by-level triangular dependency on state variables.
We have following result on the problem disturbance attenuation with stability.

\begin{theo}\em
\label{backsteping_w_level}
Consider a system given by
\be
%\label{almost_dd2}
\left\{\begin{array}{rcl}
\dot z &=& f_0(z,\xi_{1,1},\xi_{2,1},\cdots,\xi_{m,1})+p_0(z,w),\cr
\dot\xi_{i,j} &=& \xi_{i,j+1}\!+\!
\displaystyle\sum_{l=1}^{i-1}\delta_{i,j,l}(z,\xi)v_l+p_{i,j}(z,\xi,w),\quad j=1,\cdots,q_i-1,\cr
\dot\xi_{i,q_i} &=& v_i+p_{i,q_i}(z,\xi)w,\cr
y_i &=& \xi_{i,1},\;\; i=1,\cdots,m,
\end{array}
\right.
\ee
where $\xi\!=\!\col\{\xi_1,\xi_2,\cdots,\xi_m\}$,
$\xi_i\!=\!\col\{\xi_{i,1},\xi_{i,2},\cdots,\xi_{i,q_i}\}$,
$q_1\leq q_2\leq\cdots\leq q_m$,
\[
\delta_{i,j,l}=0,\; \mbox{for } \; j<q_l, \; i=1,2,\cdots,m,
\]
and functions
$f_0$, $p_0$, $\delta_{i,j,l}$ and $p_{i,j}$,  $j=1,2,\cdots,q_i, \;i=1,2,$ $\cdots,m$, are smooth with
$f_0(0,0,\cdots,0)=0$.
Assume that
\begin{eqnarray*}
\|p_0(z,w)\| &\leq &  R_0(z) \|w\|,\qquad \forall z,w,\\
|p_{i,j}(z,\xi,w)| &\leq & R_{i,j}(z,\xi) \|w\|, \quad \forall z,\xi,w,\; j=1,2,\cdots,q_i,\; i=1,2,\cdots,m,
\end{eqnarray*}
for some smooth functions $R_0(z)$ and $R_{i,j}(z,\xi)$, $j=1,2,$ $\cdots,q_i,$
$i=1,2,\cdots,m$.
Suppose that
\ben
\item[1)]
there exist a number $\gamma>0$, smooth function $\phi_{i,1}(z)$, with $\phi_{i,1}(0)=0$, $i=1,2,\cdots,m$,
a smooth positive definite and radially unbounded function $V(z)$, and a class $\calK_\infty$
function $\alpha_0(\cdot)$ such that
\begin{eqnarray*}
\frac{\partial V}{\partial z} [f_0(z,\phi_{1,1}(z),\phi_{2,1}(z),\cdots,\phi_{m,1}(z))+p_0(z,w)] \hspace{40mm}\\
\leq -\alpha_0(\|z\|)+\gamma^2\|w\|^2
-\|\col\{\phi_{1,1}(z),\phi_{2,1}(z),\cdots,\phi_{m,1}(z)\}\|^2,\qquad \qquad\;\;\,
\end{eqnarray*}
for all $z$ and $w$.
\item[2)]
the functions $\delta_{i,j,l}(z,\xi)$ and $p_{i,j}(z,\xi,\cdot)$ depend only on variables
$z$ and $\xi_{\ell_\rmp,\ell_\rmb}$, with
\ben
\item
$1\leq\ell_\rmp\leq m$ and $\ell_\rmb=1$; or,
\item
$\ell_\rmb\leq j-1$; or
\item
$\ell_\rmb=j$ and $\ell_\rmp\leq i$.
\een
\een
Then, for every $\epsilon>0$, there exist smooth feedback laws $v_i=v_i(z,\xi)$, $i=1,2,\cdots,m$,
such that the resulting closed-loop system is strictly dissipative with respect to the supply rate
$
q(w,y)=(\gamma+\epsilon)^2\|w\|^2-\|y\|^2,
$
where $y=\col\{y_1,y_2,\cdots,y_m\}$.
\end{theo}

{\em Proof:}
The theorem can be proven by using the level-by-level backstepping design
procedure \cite{lisn08}. In each step of the procedure,
we use Lemma~\ref{backsteping_w}. Let $n_\rmd=\sum_{l=1}^m q_l$.

We start the backstepping with
\be\label{step0}
\left\{\begin{array}{rcl}
\dot z &=& f_0(z,\xi_{1,1},\phi_{2,1}(z),\phi_{3,1}(z),\cdots,\phi_{m,1}(z))+p_0(z,w),\cr
\dot\xi_{1,1} &=& \xi_{1,2}+p_{1,1}(z,\xi_{1,1},\phi_{2,1}(z),\phi_{3,1}(z),\cdots,\phi_{m,1}(z),w),\cr
y_1 &=& \xi_{1,1},\cr
y_i &=& \phi_{i,1}(z),\;\; i=2,3,\cdots,m.\cr
\end{array}
\right.
\ee
Here, $\xi_{1,2}$ is viewed as a virtual input. By Lemma~\ref{backsteping_w},
for every $\epsilon>0$, there exist a smooth feedback law $\xi_{1,2}=\phi_{1,2}(z,\xi_{1,1})$,
a smooth positive definite and radially unbounded function $W_{1,1}(z,\xi_{1,1})$,
and a class $\calK_\infty$ function $\alpha_{1,1}(\cdot)$ such that
\begin{eqnarray}
\frac{\partial W_{1,1}}{\partial z} [f_0(z,\xi_{1,1},\phi_{2,1}(z),\cdots,\phi_{m,1}(z))+p_0(z,w)]\nonumber\\
\hspace{10mm}+\frac{\partial W_{1,1}}{\partial \xi_{1,1}} [\phi_{1,2}(z,\xi_{1,1})+p_{1,1}(z,\xi_{1,1},w)]\nonumber\\
\label{step1}
\hspace{5mm}\leq -\alpha_{1,1}(\|\col\{z,\xi_{1,1}\}\|)+ (\gamma+\epsilon/n_\rmd)^2 \|w\|^2 \nonumber\\
\hspace{25mm} -\|\col\{\xi_{1,1},\phi_{2,1}(z),\cdots,\phi_{m,1}(z)\}\|^2,\qquad\,
\end{eqnarray}
for all $z, \xi_{1,1}$ and $w$.
That is, subsystem (\ref{step0}) with the feedback $\xi_{1,2}=\phi_{1,2}(z,\xi_{1,1})$
is strictly dissipative with respect to the supply rate
$
q(w,y)=(\gamma+\epsilon/n_\rmd)^2\|w\|^2-\|y\|^2.
$

Next, consider the subsystem
\be
\left\{\begin{array}{rcl}
\dot z &=& f_0(z,\xi_{1,1},\xi_{2,1},\phi_{3,1}(z),\cdots,\phi_{m,1}(z))+p_0(z,w),\cr
\dot\xi_{1,1} &=& \phi_{1,2}(z,\xi_{1,1})+p_{1,1}(z,\xi_{1,1},\xi_{2,1},\phi_{3,1}(z),\cdots,\phi_{m,1}(z),w),\cr
\dot\xi_{2,1} &=& \xi_{2,2}+p_{2,1}(z,\xi_{1,1},\xi_{2,1},\phi_{3,1}(z),\cdots,\phi_{m,1}(z),w),\cr
y_1 &=& \xi_{1,1},\cr
y_2 &=& \xi_{2,1},\cr
y_i &=& \phi_{i,1}(z),\;\; i=3,\cdots,m,\cr
\end{array}
\right.
\ee
where $\xi_{2,2}$ is viewed as a virtual input. By Lemma~\ref{backsteping_w}, and in view
of (\ref{step1}), there exist a smooth feedback law $\xi_{2,2}=\phi_{2,2}(z,\xi_{1,1},\xi_{2,1})$
such that the resulting closed-loop system is strictly dissipative with respect to the supply rate
$q(w,y)=(\gamma+2\epsilon/n_\rmd)^2\|w\|^2-\|y\|^2$.

Similarly, we step back from the remain states in the first-level, and obtain
$v_i=v_i( z;\xi_{1,1},\xi_{2,1},\cdots,\xi_{i,1}),$ $i=1,2,\cdots,b_1$,
where $b_1$ is the number of chains that contain exactly one integrator,
{\em i.e.}, $q_1=q_2=\cdots=q_{b_1}=1$.

For chains that contain more than one integrator, we have
\[
\xi_{i,2}=\phi_{i,2}( z;\xi_{1,1},\xi_{2,1},\cdots,\xi_{i,1}),\quad
i=b_1+1,b_1+2,\cdots,m.
\]
Thus, the following subsystem
\be\label{step3}
\left\{\begin{array}{rcl}
\dot z &=& f_0(z,\xi_{1,1},\xi_{2,1},\cdots,\xi_{m,1})+p_0(z,w),\cr
\dot\xi_{i,1} &=& v_i+p_{i,1}(z,\xi_{1,1},\xi_{2,1},\cdots,\xi_{m,1},w), \quad i=1,2,\cdots,b_1,\cr
\dot\xi_{i,1} &\!

\!

=& \phi_{i,2}(z,\xi_{1,1},\xi_{2,1},\cdots,\xi_{i,1})+p_{i,1}(z,\xi_{1,1},\xi_{2,1},\cdots,\xi_{m,1},w), \cr
&& \qquad\qquad\qquad i=b_1+1,b_1+2, \cdots,m,\cr
y_i &=& \xi_{i,1},\;\; i=1,2,\cdots,m,\cr
\end{array}
\right.
\ee
is strictly dissipative with respect to the supply rate
$q(w,y)=(\gamma+m\epsilon/n_\rmd)^2\|w\|^2-\|y\|^2$.

To proceed backstepping on the first state in the second-level, we view $\xi_{b_1+1,3}$ as a virtual input of the following subsystem,
which consists of (\ref{step3}) and the dynamics
\[
\dot\xi_{b_1+1,2} = \xi_{b_1+1,3}+p_{b_1+1,2}(z,\xi_{1,1},\xi_{2,1},\cdots,\xi_{m,1},\xi_{b_1+1,2},w).
\]
Again, by Lemma~\ref{backsteping_w}, there exists a smooth feedback law
\[
\xi_{b_1+1,3} = \phi_{b_1+1,3}(z,\xi_{1,1},\xi_{2,1},\cdots,\xi_{m,1},\xi_{b_1+1,2}),
\]
such that the resulting closed-loop subsystem is strictly dissipative with respect to the supply rate
$q(w,y)=[\gamma+(m+1)\epsilon/n_\rmd]^2\|w\|^2-\|y\|^2$.

Continuing in this way, we finally obtain
\begin{eqnarray*}
v_i=v_i( z;\xi_{1,1},\xi_{2,1},\cdots,\xi_{m,1};\xi_{1,2},\xi_{2,2},\cdots,\xi_{m,2};
\cdots; \\
\xi_{1,q_m-1}, \xi_{2,q_m-1},\cdots,\xi_{m,q_m-1};\xi_{i,q_m}),
\end{eqnarray*}
for chains that contain $q_m$ integrators, such that
such that the closed-loop system is strictly dissipative with respect to the supply rate
$q(w,y)=(\gamma+\epsilon)^2\|w\|^2-\|y\|^2$.
\squ

\medskip

The level-by-level backstepping procedure enlarges
the class of systems for which the disturbance attenuation problem can be solved.
The triangular dependency requirement in Theorem~\ref{backsteping_w_level}
can be further weakened if we mix the chain-by-chain backstepping and the
level-by-level backstepping and implement it
on a same system.
The following result includes the chain-by-chain backstepping and level-by-level as special cases.

\begin{theo}\label{da_mimo}\em
Consider a system in the form
\be\label{general}
\left\{\begin{array}{rcl}
\dot  z &=& f_0( z,\xi_{1,1},\xi_{2,1},\cdots,\xi_{m,1})+p_0(z,w),\cr
\dot\xi_{i,j} &=& \xi_{i,j+1}+\displaystyle \sum_{l=1}^{i-1}\delta_{i,j,l}( z,\xi)v_l,+p_{i,j}(z,\xi,w),\quad j=1,2,\cdots,q_i-1,\cr
\dot\xi_{i,q_i} &=& v_i+p_{i,q_i}(z,\xi,w),\cr
y_i &=& \xi_{i,1},\;\; i=1,2,\cdots,m,
\end{array}
\right.
\ee
where $\xi=\col\{\xi_1,\xi_2,\cdots,\xi_m\}$, $\xi_i=\col\{\xi_{i,1},\xi_{i,2},\cdots,$ $\xi_{i,q_i}\}$,
$q_1\leq q_2\leq \cdots \leq q_m$, and functions
$f_0$, $p_0$, $\delta_{i,j,l}$ and $p_{i,j}$, $i=1,2,\cdots,m$ are smooth.
Assume that
\begin{eqnarray*}
\|p_0(z,w)\| &\leq & R_0(z) \|w\|,\;\forall z,w,\\
|p_{i,j}(z,\xi,w)| &\leq & R_{i,j}(z,\xi) \|w\|,\quad \forall z,\xi,w,\; j=1,2,\cdots,q_i,\; i=1,2,\cdots,m,
\end{eqnarray*}
for some smooth functions $R_0(z)$ and $R_{i,j}(z,\xi)$,
$j=1,2,$ $\cdots,q_i,$ $i=1,2,\cdots,m$,
Suppose that
\ben
\item[1)]
Condition 1) in Theorem~\ref{backsteping_w_level} holds,
\item[2)]
there exists an ordered list $\kappa$ containing all variables of $\xi$
such that, for $j=1,2,\cdots,q_i-1$, $l=1,2,\cdots, i-1$, $i=1,2,\cdots,m$,
\ben
\item
$\xi_{i,j}$ appears earlier than $\xi_{i,j+1}$ in $\kappa$;
\item for $\delta_{i,j,l}\neq0$, the variables $\xi_l$ appear earlier than $\xi_{i,j}$
in $\kappa$;
\item
$\delta_{i,j,l}(z,\xi)$ and $p_{i,j}(z,\xi,\cdot)$ depend only on $z$, $\xi_{\ell,1}$, $\ell=1,2,\cdots,m$, and
the variables that appear no later than $\xi_{i,j}$
in $\kappa$.
\een
\een
Then, for every $\epsilon>0$, there exist smooth feedback laws $v_i=v_i(z,\xi)$, $i=1,2,\cdots,m$,
such that the resulting closed-loop system is strictly dissipative with respect to the supply rate
$
q(w,y)=(\gamma+\epsilon)^2\|w\|^2-\|y\|^2,
$
where $y=\col\{y_1,y_2,\cdots,y_m\}$.
\end{theo}

{\em Proof:} The backstepping can be carried out one state by one state in the order
of the list $\kappa$. Suppose that, after backstepping $\ell$ state variables,
we want to backstep from $\xi_{i,j}$, the $\ell+1$-th element in the list $\kappa$,
to $\xi_{i,j+1}$.
Let $n_\rmd=\sum_{l=1}^m q_l$.
Denote all the state variables that come before $\xi_{i,j}$ in the list $\kappa$ as $Z$.
By Condition 2), we can describe the subsystem of $Z$ and $\xi_{i,j}$ as
\be\label{sub2}
\left\{\begin{array}{rcl}
\dot  Z &=& F_0( Z,\xi_{i,j})+P_0(Z,w),\cr
\dot\xi_{i,j} &=& \xi_{i,j+1}+\displaystyle \sum_{l=1}^{i-1}\delta_{i,j,l}(Z,\xi_{i,j})v_l(Z)+P_{i,j}(Z,\xi_{i,j},w),\cr
y_\iota &=& \xi_{\iota,1},\;\; \iota=1,2,\cdots,m,
\end{array}
\right.
\ee
where $\xi_{i,j+1}$ is viewed as a virtual input. By Lemma~\ref{backsteping_w}, there exists a smooth function
\[
\xi_{i,j+1}=\phi_{i,j+1}(Z,\xi_{i,j})
\]
such that the resulting closed-loop subsystem is strictly dissipative with respect to the supply rate
$
q(w,y)=[\gamma+ (\ell+1)\epsilon/n_\rmd]^2\|w\|^2-\|y\|^2.
$
By backstepping through all the state variables in the list $\kappa$, we can find the desired
feedback laws $v_i=v_i(z,\xi)$, $i=1,2,\cdots,m$.
\squ

\medskip

We illustrate the backstepping procedure of Theorem~\ref{da_mimo} by the following example.

\begin{example}\rm
Consider a system in the form of
(\ref{general}), with $q_1=2$ and $q_2=3$,
\be\label{example23_dist}
\left\{\begin{array}{rcl}
\dot z &=& f_0(z,\xi_{1,1},\xi_{2,1})+ R_0(z)w,\cr
\dot\xi_{1,1} &=& \xi_{1,2}+\xi_{2,1}w,\cr
\dot\xi_{1,2} &=& v_1+\xi_{1,2}w,\cr
\dot\xi_{2,1} &=& \xi_{2,2}+z\sin w,\cr
\dot\xi_{2,2} &=& \xi_{2,3}+\xi_{2,1}v_1+zw,\cr
\dot\xi_{2,3} &=& v_2 ,\cr
y_1&=&\xi_{1,1},\cr
y_2&=&\xi_{2,1}.
\end{array}\right.
\ee

We first note that the triangular dependency condition (\ref{tri1}), needed for the
conventional backstepping, does not hold for this system. We will thus resort to Theorem \ref{da_mimo}.
Obviously, Condition 2) in Theorem~\ref{da_mimo} holds with the ordered list
$\kappa=\{\xi_{2,1},\xi_{1,1},\xi_{1,2},\xi_{2,2},\xi_{2,3}\}$.
Suppose that
there exist a number $\gamma>0$, smooth functions $\phi_{i,1}(z)$, with $\phi_{i,1}(0)=0$, $i=1,2$,
a smooth positive definite and radially unbounded function $V(z)$, and a class $\calK_\infty$
function $\alpha_0(\cdot)$ such that
\begin{eqnarray*}
\frac{\partial V}{\partial z} [f_0(z,\phi_{1,1}(z),\phi_{2,1}(z))+R_0(z)w] \qquad\qquad\qquad\qquad\cr
\qquad\qquad\leq -\alpha_0(\|z\|)+\gamma^2\|w\|^2-\|\col\{\phi_{1,1}(z),\phi_{2,1}(z)\}\|^2,
\end{eqnarray*}
for all $z$ and $w$.

Consider the following subsystem,
\be\label{equ_delta6}
\left\{\begin{array}{rcl}
\dot z &=& f_0(z,\phi_{1,1}(z),\xi_{2,1})+ R_0(z)w,\cr
\dot\xi_{2,1} &=& \xi_{2,2}+z\sin w,\cr
y_1&=&\phi_{1,1}(z),\cr
y_2&=&\xi_{2,1}.
\end{array}\right.
\ee
View $\xi_{2,2}$ as a virtual input.
By Lemma~\ref{backsteping_w}, for every $\epsilon>0$, there exists a smooth feedback $\xi_{2,2}=\phi_{2,2}(z,\xi_{2,1})$ such that
subsystem (\ref{equ_delta6}) is strictly dissipative with respect to the supply rate $q(w,y)=(\gamma+\epsilon/5)^2w^2-\|y\|^2$.

Next consider
\be\label{equ_delta7}
\left\{\begin{array}{rcl}
\dot z &=& f_0(z,\xi_{1,1},\xi_{2,1})+ R_0(z)w,\cr
\dot\xi_{1,1} &=& \xi_{1,2}+\xi_{2,1}w,\cr
\dot\xi_{2,1} &=& \phi_{2,2}(z,\xi_{2,1})+z\sin w,\cr
y_1&=&\xi_{1,1},\cr
y_2&=&\xi_{2,1},
\end{array}\right.
\ee
and view $\xi_{1,2}$ as a virtual input.
Again, by Lemma~\ref{backsteping_w}, there exists a smooth feedback $\xi_{1,2}=\phi_{1,2}(z,\xi_{2,1},$ $\xi_{1,1})$ such that
subsystem (\ref{equ_delta7}) is strictly dissipative with respect to the supply rate $q(w,y)=(\gamma+2\epsilon/5)^2w^2-\|y\|^2$.

By backstepping in a similar way through $\xi_{1,2},\xi_{2,2},\xi_{2,3}$,
we obtain the smooth feedback laws
\[
\left.\begin{array}{rcl}
v_1 &=& v_1(z,\xi_{2,1},\xi_{1,1},\xi_{1,2}),\cr
v_2 &=& v_2(z,\xi_{2,1},\xi_{1,1},\xi_{2,2},\xi_{2,3}),\cr
\end{array}\right.
\]
such that
system (\ref{example23_dist}) is strictly dissipative with respect to the supply rate $q(w,y)=(\gamma+\epsilon)^2w^2-\|y\|^2$.
\end{example}
\medskip

As the problem of almost disturbance decoupling is a special case of the problem of
disturbance attenuation, the following result on almost disturbance
decoupling with stability is a corollary to Theorem \ref{da_mimo}.

\newpage
\begin{coro}\label{add_mimo}\em
Consider a system in the form
\be\label{general2}
\left\{\begin{array}{rcl}
\dot  z &=& f_0( z,\xi_{1,1},\xi_{2,1},\cdots,\xi_{m,1},w),\cr
\dot\xi_{i,j} &=& \xi_{i,j+1}+\displaystyle \sum_{l=1}^{i-1}\delta_{i,j,l}( z,\xi)v_l,+p_{i,j}(z,\xi,w),\quad j=1,2,\cdots,q_i-1,\cr
\dot\xi_{i,q_i} &=& v_i+p_{i,q_i}(z,\xi,w),\cr
y_i &=& \xi_{i,1},\;\; i=1,2,\cdots,m,
\end{array}
\right.
\ee
where $\xi=\col\{\xi_1,\xi_2,\cdots,\xi_m\}$, $\xi_i=\col\{\xi_{i,1},\xi_{i,2},\cdots,$ $\xi_{i,q_i}\}$,
$q_1\leq q_2\leq \cdots \leq q_m$, and functions
$f_0$, $\delta_{i,j,l}$ and $p_{i,j}$, $i=1,2,\cdots,m$, are smooth.
Assume that
\begin{eqnarray*}
|p_{i,j}(z,\xi,w)|\leq R_{i,j}(z,\xi) \|w\|, \;\forall z,\xi,w,\quad j=1,2,\cdots,q_i,\; i=1,2,\cdots,m,
\end{eqnarray*}
for some smooth functions $R_{i,j}(z,\xi)$, $j=1,2,\cdots,q_i,$ $i=1,2,\cdots,m$.
Suppose that
\ben
\item[1)]
for every $\gamma_0>0$, there exist smooth $\phi_{i,1}(z)$ with $\phi_{i,1}(0)=0$, $i=1,2,\cdots,m$,
and a smooth positive definite and radially unbounded function $V(z)$
such that
\begin{eqnarray*}
\frac{\partial V}{\partial z} f_0(z,\phi_{1,1}(z),\phi_{2,1}(z),\cdots,\phi_{m,1}(z),w)\qquad\qquad\\
\leq -\alpha(\|z\|)+\gamma_0^2 \|w\|^2
-\|\col\{\phi_{1,1}(z),\phi_{2,1}(z),\cdots,\phi_{m,1}(z)\}\|^2,\;\forall z,w,
\end{eqnarray*}
for some $\calK_\infty$ function $\alpha(\cdot)$.

\item[2)]
Condition 2) in Theorem~\ref{da_mimo} holds.

\een
Then, for every $\gamma>0$, there exist smooth feedback laws $v_i=v_i(z,\xi)$, $i=1,2,\cdots,m$,
such that the resulting closed-loop system is strictly dissipative with respect to the supply rate
$
q(w,y)=\gamma^2\|w\|^2-\|y\|^2,
$
where $y=\col\{y_1,y_2,\cdots,y_m\}$.
\end{coro}

\medskip

We next consider further the fulfillment of Condition 1) in Corollary~\ref{add_mimo}.
It is a generalization of Lemma~\ref{z1z2} to multiple input multiple output systems.
Suppose that the $z$-subsystem
\be\label{almost_z_2a}
\dot z=f_0(z,\xi_{1,1},\xi_{2,1},\cdots,\xi_{m,1},w)
\ee
can be decomposed as
\be\label{almost_z_2}
\left\{\begin{array}{rcl}\dot z_1 &=& f_1(z_1,z_2,\xi_{1,1},\xi_{2,1},\cdots,\xi_{m,1},w),\cr
                           z_2 &=& f_2(z_2,\xi_{1,1},\xi_{2,1},\cdots,\xi_{m,1}),\end{array}\right.
\ee
where $z_1$ represents ``stable component" and $z_2$ represents ``unstable but stabilizable component."
We have the following result.

\newpage
\begin{coro}\label{non_mini_mimo}\em
Consider system (\ref{almost_z_2a}) which can be decomposed as (\ref{almost_z_2}). Suppose that
\ben
\item[1)]
there exists a smooth positive definite and radially unbounded function $V_1(z_1)$ such that
\begin{eqnarray*}
&&\frac{\partial V_1}{\partial z_1} f_1( z_1,z_2,\xi_{1,1},\xi_{2,1},\cdots,\xi_{m,1},w)\\
&&\qquad \leq -\alpha_1(\|z_1\|)+\gamma_0^2\|z_2\|^2+\gamma_0^2\|w\|^2
+\gamma_0^2\|\col\{\xi_{1,1},\xi_{2,1},\cdots,\xi_{m,1}\}\|^2,
\end{eqnarray*}
for some $\calK_\infty$ function $\alpha_1$ and some $\gamma_0>0$, and
\item[2)]
there exist smooth functions $\bar v_i(z_2)$ with $\bar v_i(0)=0$ for $i=1,2,\cdots,m$,
and a smooth positive definite and radially unbounded function $V_2(z_2)$ such that
\begin{eqnarray*}
\frac{\partial V_2}{\partial z_2} f_2(z_2,\bar v_1(z_2),\bar v_2(z_2),\cdots,\bar v_m(z_2))\hspace{30mm}\\
+\|\col\{\bar v_1(z_2),\bar v_2(z_2),\cdots,\bar v_m(z_2)\}\|^2 \\
\hspace{70mm}\leq -\alpha_2(\|z_2\|)
\end{eqnarray*}
for some $\calK_\infty$ function $\alpha_2$.
\een
Then, for every $\gamma>0$, there exist smooth $v_i(z)$ with $v_i(0)=0$, $i=1,2,\cdots,m$, and a smooth positive definite and radially unbounded function $V(z)$
such that
\begin{eqnarray*}
\frac{\partial V}{\partial z} f_0(z,v_1(z),v_2(z),\cdots,v_m(z),w)\leq -\alpha(\|z\|)\qquad\qquad\\
+\gamma^2 \|w\|^2-\|\col\{v_1(z),v_2(z),\cdots,v_m(z)\}\|^2,\;
 \forall z,w,
\end{eqnarray*}
for some $\calK_\infty$ function $\alpha(\cdot)$.
\end{coro}

Corollary~\ref{non_mini_mimo} provides, for the system in Corollary~\ref{add_mimo}, a starting point
from which the backstepping can be carried out. We next use a numerical example with unstable
zero dynamics to illustrate Corollary~\ref{non_mini_mimo}.

\begin{example}\rm
Consider a system in the form of
(\ref{general2}) with $q_1=1$ and $q_2=2$,
\be\label{add_exam}
\left\{\begin{array}{rcl}
\dot z &=& z+\xi_{1,1}+\xi_{2,1},\cr
\dot\xi_{1,1} &=& v_1+\xi_{2,1}w,\cr
\dot\xi_{2,1} &=& \xi_{2,2}+z w,\cr
\dot\xi_{2,2} &=& v_2+(\cos\xi_{1,1})\sin w,\cr
y_1 &=& \xi_{1,1},\cr
y_2 &=& \xi_{2,1}.
\end{array}\right.
\ee

Note that the dependency requirement in Corollary \ref{add_mimo} holds and
zero dynamics satisfies the conditions in Corollary~\ref{non_mini_mimo}.
The zero dynamics $\dot z=z$ is unstable.
View $\xi_{1,1}$ and $\xi_{2,1}$ as virtual input of
\be\label{zzz}
\dot z = z+\xi_{1,1}+\xi_{2,1}.
\ee
Then $\xi_{1,1}=\phi_{1,1}(z)=-2z$ and $\xi_{2,1}=0$ stabilizes (\ref{zzz})
with Lyapunov function $V_0=z^2/2$.

Condition 2) of Corollary~\ref{add_mimo} holds with $\kappa=\{\xi_{2,1},\xi_{1,1},$ $\xi_{2,2}\}$.

To begin the mixed chain-by-chain and level-by-level backstepping procedure, we consider the subsystem
\be\label{add_exam1}
\left\{\begin{array}{rcl}
\dot z &=& -z+\xi_{2,1},\cr
\dot\xi_{2,1} &=& \xi_{2,2}+z w,\cr
\end{array}\right.
\ee
and view $\xi_{2,2}$ as a virtual input.
Consider the Lyapunov function
$
V_1=V_0+\xi_{2,1}^2/2=z^2/2+\xi_{2,1}^2/2.
$
Its time derivative is given by
\[
\dot V_1=\dot V_0+z\xi_{2,1}+\xi_{2,1}\xi_{2,2}+\xi_{2,1}zw.
\]
Let
\[
\xi_{2,2}=\phi_{2,2}(z,\xi_{2,1})=-z-\xi_{2,1}-\frac{3}{4\gamma^2}\xi_{2,1}(1+z^2),
\]
which renders
\[
\dot V_1\leq -\xi_{2,1}^2+ \frac{\gamma^2}{3}\|w\|^2-z^2.
\]

We next consider
\be\label{add_exam2}
\left\{\begin{array}{rcl}
\dot z &=& z+\xi_{1,1}+\xi_{2,1},\cr
\dot\xi_{2,1} &=& -z -\xi_{2,1}-\frac{3}{4\gamma^2}\xi_{2,1}(1+z^2)+z w,\cr
\dot\xi_{1,1} &=& v_1+\xi_{2,1}w.
\end{array}\right.
\ee
Letting
$
V_2=V_1+(\xi_{1,1}+2z)^2/2,
$
we have
\[
\dot V_2\leq \dot V_1 +(\xi_{1,1}+2z)(v_1+3z+\xi_{1,1}+\xi_{2,1}+\xi_{2,1}w).
\]
Let
\[
v_1=-\frac{11}{3}z-\frac{4}{3}\xi_{1,1}-\xi_{2,1}-\frac{3}{4\gamma^2} (\xi_{1,1}+2z) (1+\xi_{2,1}^2).
\]
We have
\begin{eqnarray*}
\dot V_2&\leq& -\xi_{2,1}^2+ \frac{2\gamma^2}{3} w^2-z^2-(\xi_{1,1}+2z)^2/3\\
&\leq& -\xi_{2,1}^2+ \frac{2\gamma^2}{3} w^2-4z^2 -\xi_{1,1}^2.
\end{eqnarray*}
Finally, consider
\be\label{add_exam3}
\left\{\begin{array}{rcl}
\dot z &=& z+\xi_{1,1}+\xi_{2,1},\cr
\dot\xi_{1,1} &=& v_1+\xi_{2,1}w,\cr
\dot\xi_{2,1} &=& \xi_{2,2}+z w,\cr
\dot\xi_{2,2} &=& v_2+(\cos\xi_{1,1}) \sin w,\cr
\end{array}\right.
\ee
for which we let
$
V_3=V_2+(\xi_{2,2}-\phi_{2,2})^2/2.
$
Thus,
\begin{eqnarray*}
\dot V_3 &=& \dot V_2+(\xi_{2,2}-\phi_{2,2})\left(\xi_{2,1}+v_2+(\cos\xi_{1,1}) \sin w-\dot \phi_{2,2}\right)\\
 &=& \dot V_2+(\xi_{2,2}-\phi_{2,2}) \left( v_2 + \Psi + \Phi w+ (\cos\xi_{1,1}) \sin w \right),
\end{eqnarray*}
where
\begin{eqnarray*}
\Phi &=& z +\frac{3}{4\gamma^2}(z+z^3),\\
\Psi &=& z+\xi_{1,1}+2\xi_{2,1}+\xi_{2,2} +\frac{3}{4\gamma^2}\xi_{2,2}(1+z^2)
+\frac{3}{2\gamma^2}\xi_{2,1}z(z+\xi_{1,1}+\xi_{2,1}).
\end{eqnarray*}
Let
\begin{eqnarray*}
v_2&=&-\xi_{2,2}+\phi_{2,2}-\Psi
-\frac{3}{4\gamma^2}(\xi_{2,2}-\phi_{2,2}) \left( 1+(|\Phi|+|\cos\xi_{1,1}|)^2\right).
\end{eqnarray*}
We have
\begin{eqnarray*}
\dot V_3&\leq& -(\xi_{2,2}-\phi_{2,2})^2  -\xi_{2,1}^2+ \gamma^2 w^2-4z^2 -\xi_{1,1}^2\\
 &\leq&  -\xi_{2,1}^2-\xi_{1,1}^2 + \gamma^2 w^2,
\end{eqnarray*}
from which we have
\[
\int_0^t\|y(\tau)\|^2 d\tau\leq \gamma^2 \int_0^t w(\tau)^2 d\tau,
\]
in the absence of initial condition. In the presence of
initial condition $x(0)$, we have
\[
\int_0^t\|y(\tau)\|^2 d\tau\leq \gamma^2 \int_0^t w(\tau)^2 d\tau+V_3(x(0)).
\]
\begin{figure}[ht!]
\centerline{
\includegraphics[width=100mm]{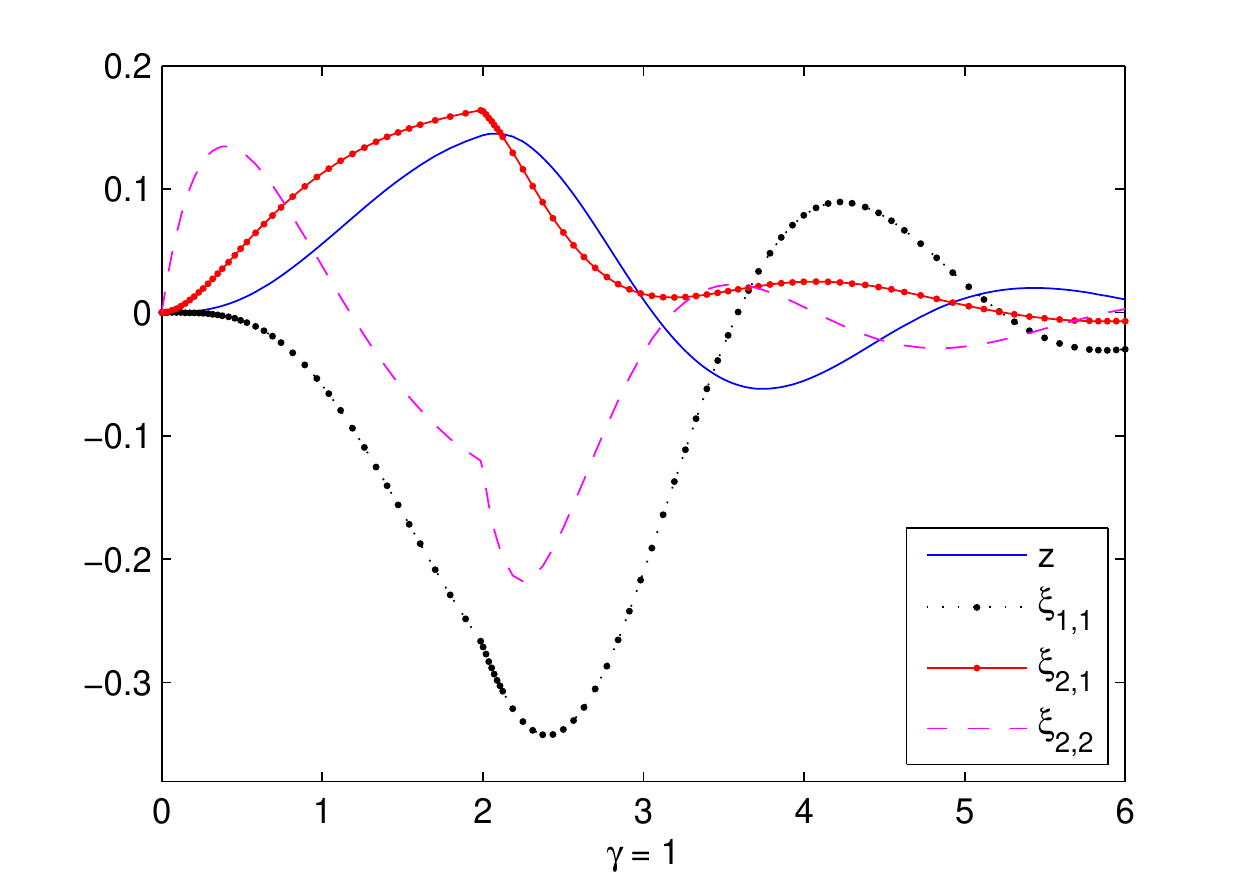}
}
\centerline{
\includegraphics[width=100mm]{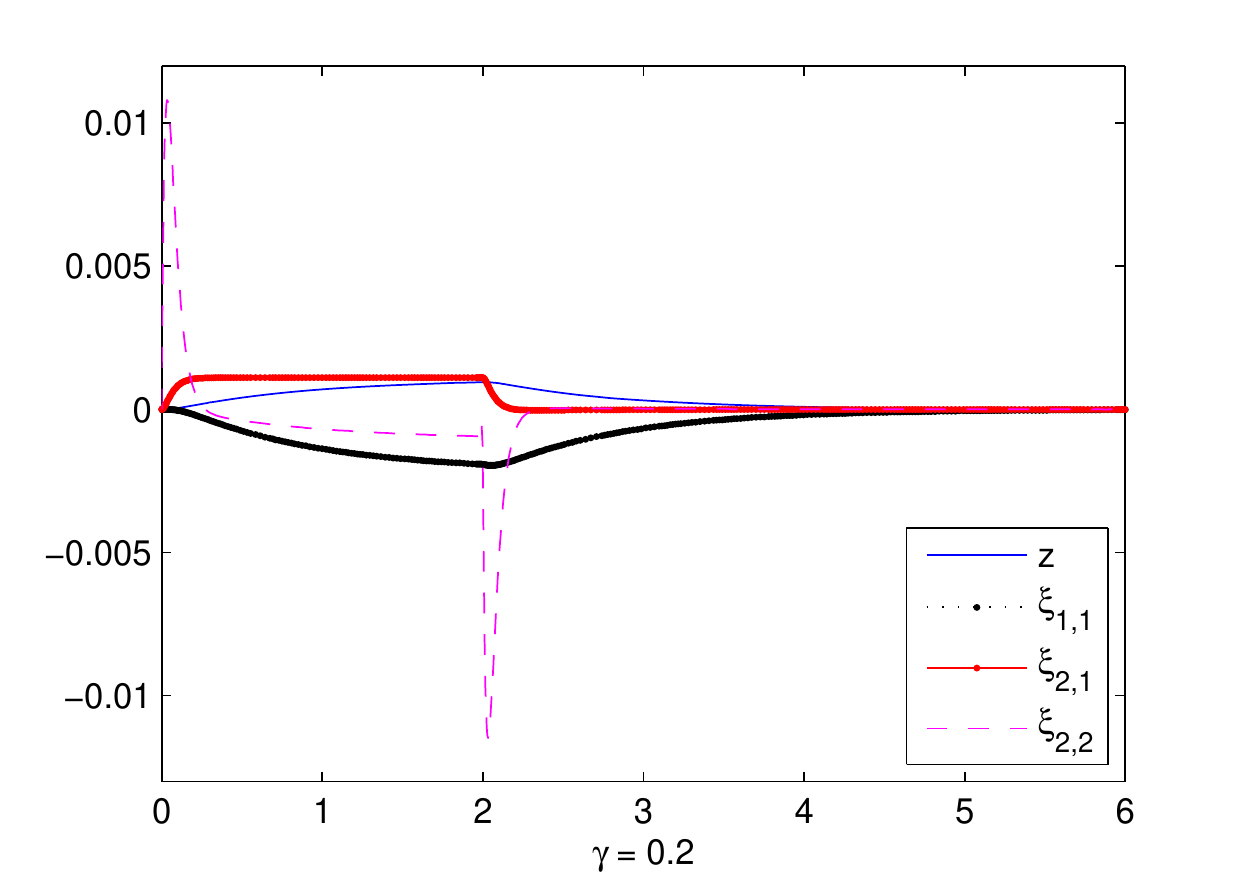}
}
\caption{
State trajectories with $x(0)=0$ and $w(t)=1(t)-1(t-2)$.}
\label{zero}
\end{figure}
\begin{figure}[ht!]
\centerline{
\includegraphics[width=100mm]{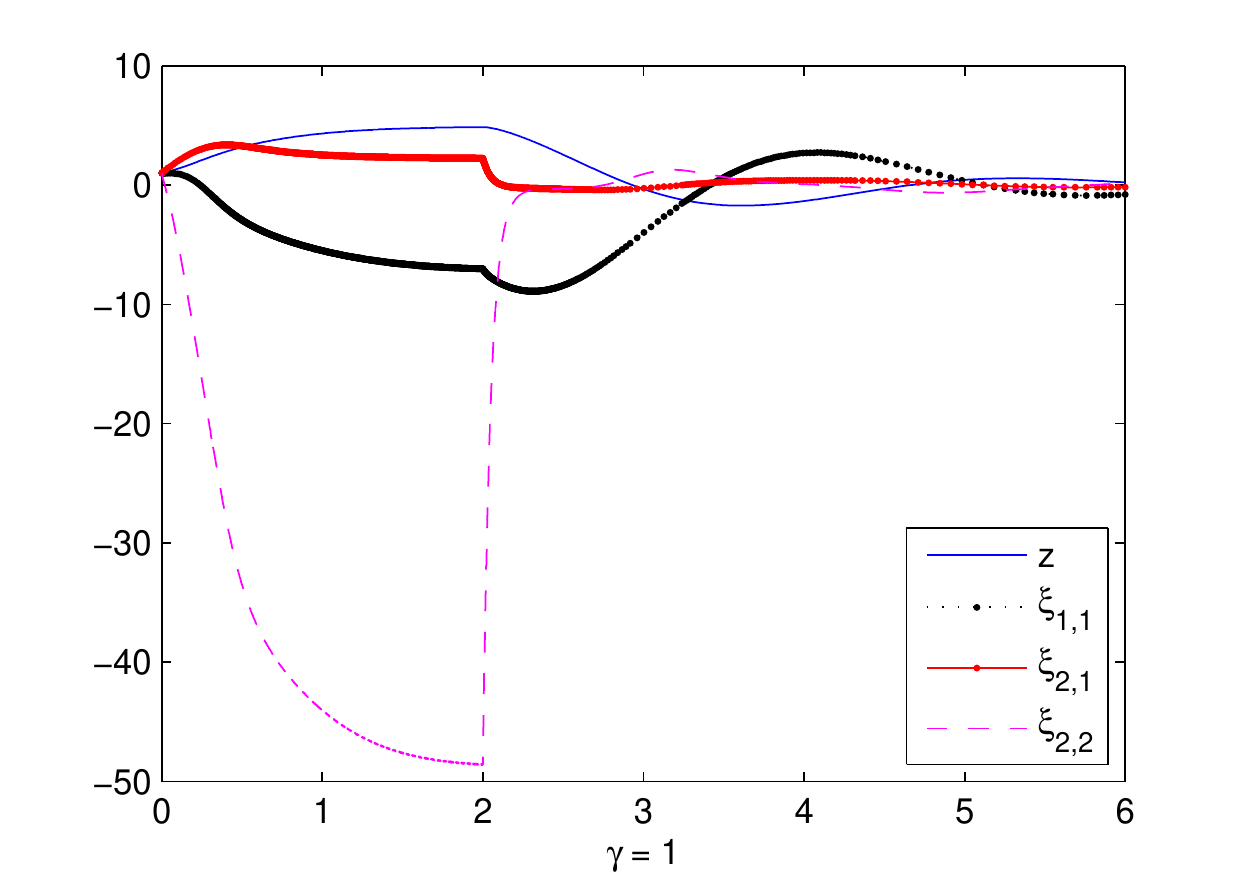}
}
\centerline{
\includegraphics[width=100mm]{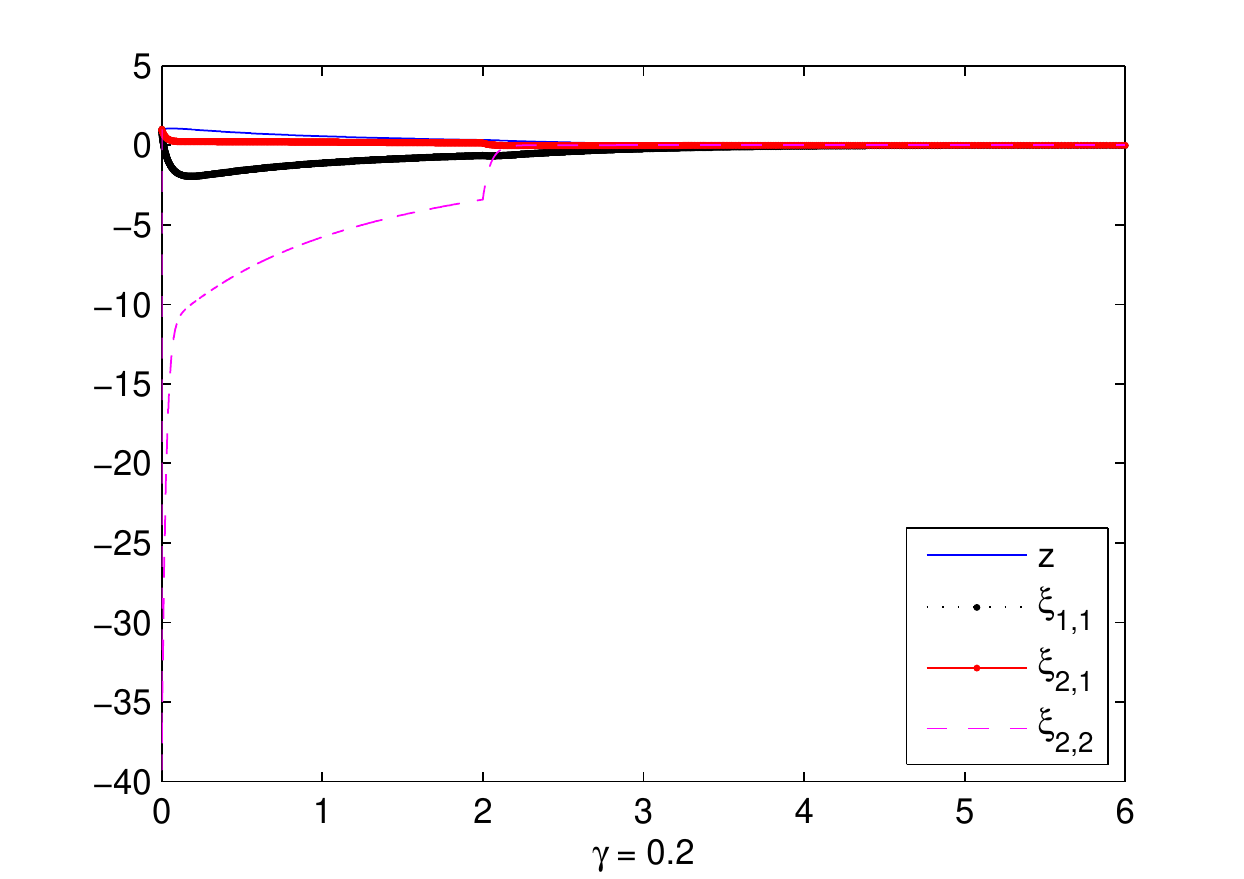}
}
\caption{
State trajectories with $x(0)=[1\;1\;1\;1]^T$ and $w(t)=\left(1(t)-1(t-2)\right)\times 10$.}
\label{nonzero}
\end{figure}
Shown in Fig.~\ref{zero} are some state trajectories of the closed-loop system
with $x(0)=0$ and $w(t)=1(t)-1(t-2)$.
Shown in Fig.~\ref{nonzero} are some state trajectories with $x(0)=[1\;1\;1\;1]^T$
and $w(t)=\left(1(t)-1(t-2)\right)\times 10$.

\end{example}

\section{Summary of the Chapter}
\label{conclusion_dist}

In this chapter, we have revisited the problems disturbance attenuation and almost
disturbance decoupling for nonlinear systems and showed how a recently developed
structural decomposition of multiple input multiple output systems and the
new backstepping design procedure it motivates can lead to the solution
of these two problems for a larger class of systems.

\chapter{Summary}

In this note,
we have obtained a few further results in differential geometric nonlinear control theory.
We first developed normal forms for nonlinear system affine in control.
Then, based on these normal forms, we revisited stabilization, semi-global stabilization and disturbance attenuation.

We presented constructive algorithms for decomposing an
affine nonlinear system into its normal form representations.
Such algorithms generalize the existing results in several ways.
They require fewer restrictive assumptions on the system
and apply to general multiple input multiple output noninear systems that do not necessarily
have the same number of inputs and outputs.
The resulting normal forms reveal various nonlinear extensions of
linear system structural properties.
These algorithms and the resulting normal forms
are thus expected to facilitate the solution of several nonlinear
control problems.

We exploited the properties of the
structural decomposition for the stabilization of multiple input and multiple output
systems, and showed that this
decomposition simplifies the conventional chain-by-chain backstepping design procedure and
motivates a new level-by-level backstepping design procedure that is able to
stabilize some systems for which the conventional backstepping procedure is
not applicable.
The chain-by-chain and level-by-level backstepping procedures
can be combined to form a mixed backstepping design
technique. The enlarged class of systems that can be stabilized by this mixed backstepping
design procedure is characterized in the form of a theorem.

We then showed how the
structural decomposition can be used to solve the problem of
semi-global stabilization for a class of multiple input multiple output systems
without vector relative degrees. The design procedure
involved several existing design techniques in
nonlinear stabilization, including low gain
feedback and different forms of backstepping design procedures.

We also revisited the problems of disturbance attenuation and almost
disturbance decoupling for nonlinear systems and showed how the
structural decomposition of nonlinear systems and the
new backstepping design procedures it motivates can lead to the solution
of these two problems for a larger class of systems.

For the future research, we are interested in the problems of
non-interacting control, tracking and regulation of nonlinear systems.
These control problems can be dealt with based on the normal forms in Chapter 3.

Output feedback control is a more challenging problem.
The normal forms proposed in the note, which reveals system structure at infinity,
will also facilate the construction of high gain observers, which will result in output feedback laws.

The structural algorithms can be applied to general nonlinear systems that are not necessarily
square invertible.
We have only considered their application to square invertible nonlinear systems in the note.
We will to utilize these normal forms to study control problems for non-invertible nonlinear systems,
in particular, underactuated nonlinear systems.

%\include{preliminaries}
%\include{linear}
%\include{nonlinear}
%\include{chapter}
%\include{chapter}
%\appendix
%\include{appendix}

\backmatter%%%%%%%%%%%%%%%%%%%%%%%%%%%%%%%%%%%%%%%%%%%%%%%%%%%%%%%
%\include{solutions}
%\include{referenc}


%\bibliographystyle{IEEEtran}
%\bibliography{C:/Xinmin/Dropbox/XMPaper/string,C:/Xinmin/Dropbox/XMPaper/nonlinear,C:/Xinmin/Dropbox/XMPaper/linear,C:/Xinmin/Dropbox/XMPaper/medical,C:/Xinmin/Dropbox/XMPaper/xm_papers}

\begin{thebibliography}{10}
\providecommand{\url}[1]{#1}
\csname url@samestyle\endcsname
\providecommand{\newblock}{\relax}
\providecommand{\bibinfo}[2]{#2}
\providecommand{\BIBentrySTDinterwordspacing}{\spaceskip=0pt\relax}
\providecommand{\BIBentryALTinterwordstretchfactor}{4}
\providecommand{\BIBentryALTinterwordspacing}{\spaceskip=\fontdimen2\font plus
\BIBentryALTinterwordstretchfactor\fontdimen3\font minus
  \fontdimen4\font\relax}
\providecommand{\BIBforeignlanguage}[2]{{%
\expandafter\ifx\csname l@#1\endcsname\relax
\typeout{** WARNING: IEEEtran.bst: No hyphenation pattern has been}%
\typeout{** loaded for the language `#1'. Using the pattern for}%
\typeout{** the default language instead.}%
\else
\language=\csname l@#1\endcsname
\fi
#2}}
\providecommand{\BIBdecl}{\relax}
\BIBdecl

\bibitem{elliott1997bookreview}
D.~L. Elliott, ``{Book Reviews of Nonlinear Control Systems, Alberto Isidori,
  1995},'' \emph{IEEE Transactions on Automatic Control}, vol.~42, no.~7, pp.
  1043--1044, 1997.

\bibitem{nind90}
H.~Nijmeijer and A.~J. van~der Schaft, \emph{Nonlinear Dynamical Control
  Systems}.\hskip 1em plus 0.5em minus 0.4em\relax Springer, 1990.

\bibitem{hermann1963accessibility}
R.~Hermann, ``{On the Accessibility Problem in Control Theory},'' in
  \emph{International Symposium on Nonlinear Differential Equations and
  Nonlinear Mechanics: proceedings}.\hskip 1em plus 0.5em minus 0.4em\relax
  Academic Press, 1963, p. 325.

\bibitem{hermann1968differential}
------, \emph{{Differential geometry and the calculus of variations}}.\hskip
  1em plus 0.5em minus 0.4em\relax Academic Press, 1968.

\bibitem{hermann1977nonlinear}
R.~Hermann and A.~J. Krener, ``{Nonlinear controllability and observability},''
  \emph{IEEE Transactions on Automatic Control}, vol.~22, no.~5, pp. 728--740,
  1977.

\bibitem{brockett1979feedback}
R.~W. Brockett, ``{Feedback invariants for nonlinear systems},'' \emph{A link
  between science and applications of automatic control}, pp. 1115--1120, 1979.

\bibitem{brockett1983asymptotic}
------, ``Asymptotic stability and feedback stabilization,'' in
  \emph{Differential Geometric Control Theory}.\hskip 1em plus 0.5em minus
  0.4em\relax Birkh{\"a}user, Dec. 1983, pp. 181--191.

\bibitem{isidori1981nonlinear}
A.~Isidori, A.~Krener, C.~Gori-Giorgi, and S.~Monaco, ``{Nonlinear decoupling
  via feedback: a differential geometric approach},'' \emph{IEEE Transactions
  on Automatic Control}, vol.~26, no.~2, pp. 331--345, 1981.

\bibitem{hirschorn1981mathcal}
R.~Hirschorn, ``{$(A, B)$-invariant distributions and disturbance decoupling of
  nonlinear systems},'' \emph{SIAM Journal on Control and Optimization},
  vol.~19, p.~1, 1981.

\bibitem{mosi73}
A.~S. Morse, ``Structural invariants of linear multivariable systems,''
  \emph{SIAM Journal on Control}, vol.~11, pp. 446--465, 1973.

\bibitem{wolm79}
W.~M. Wonham, \emph{{Linear Multivariable Control: A Geometric
  Approach}}.\hskip 1em plus 0.5em minus 0.4em\relax Springer-Verlag, 1979.

\bibitem{wonham1970decoupling}
W.~M. Wonham and A.~S. Morse, ``{Decoupling and pole assignment in linear
  multivariable systems: A geometric approach},'' \emph{SIAM Journal on
  Control}, vol.~8, no.~1, pp. 1--18, 1970.

\bibitem{wonham1972feedback}
------, ``{Feedback invariants of linear multivariable systems},''
  \emph{Automatica}, vol.~8, no.~1, pp. 93--100, 1972.

\bibitem{basile1969controlled}
G.~Basile and G.~Marro, ``{Controlled and conditioned invariant subspaces in
  linear system theory},'' \emph{Journal of Optimization Theory and
  Applications}, vol.~3, no.~5, pp. 306--315, 1969.

\bibitem{basile1992controlled}
------, \emph{{Controlled and conditioned invariants in linear system
  theory}}.\hskip 1em plus 0.5em minus 0.4em\relax Prentice Hall Englewood
  Cliffs, New Jersey, 1992.

\bibitem{isnc95}
A.~Isidori, \emph{Nonlinear Control Systems}, 3rd~ed.\hskip 1em plus 0.5em
  minus 0.4em\relax Springer, 1995.

\bibitem{isnc99}
------, \emph{{Nonlinear Control Systems II}}.\hskip 1em plus 0.5em minus
  0.4em\relax Springer, 1999.

\bibitem{byas91}
C.~I. Byrnes and A.~Isidori, ``{Asymptotic stabilization of minimum phase
  nonlinear systems},'' \emph{IEEE Transactions on Automatic Control}, vol.~36,
  pp. 1122--1137, 1991.

\bibitem{saac89}
S.~S. Sastry and A.~Isidori, ``{Adaptive control of linearizable systems},''
  \emph{IEEE Transactions on Automatic Control}, vol.~34, no.~11, pp.
  1123--1131, 1989.

\bibitem{marino1994equivalence}
R.~Marino, W.~Respondek, and A.~J. van~der Schaft, ``{Equivalence of nonlinear
  systems to input-output prime forms},'' \emph{SIAM Journal of Control and
  Optimization}, vol.~32, p. 387, 1994.

\bibitem{byls88}
C.~I. Byrnes and A.~Isidori, ``{Local stabilization of minimum-phase nonlinear
  systems},'' \emph{Systems \& Control Letters}, vol.~11, no.~1, pp. 9--17,
  1988.

\bibitem{scgn99}
B.~Schwartz, A.~Isidori, and T.~J. Tarn, ``Global normal forms for mimo
  nonlinear systems, with application to stabilization and disturbance
  attenuation,'' \emph{Mathematics of Control, Signals and Systems}, vol.~12,
  pp. 121--142, 1999.

\bibitem{chu2002numerical}
D.~Chu, X.~Liu, and R.~C. Tan, ``On the numerical computation of a structural
  decomposition in systems and control,'' \emph{IEEE Transactions on Automatic
  Control}, vol.~47, no.~11, pp. 1786--1799, 2002.

\bibitem{liu2003problem}
X.~Liu, B.~M. Chen, and Z.~Lin, ``On the problem of general structural
  assignments of linear systems through sensor/actuator selection,''
  \emph{Automatica}, vol.~39, no.~2, pp. 233--241, 2003.

\bibitem{liu2005linear}
------, ``Linear systems toolkit in matlab: structural decompositions and their
  applications,'' \emph{Journal of Control Theory and Applications}, vol.~3,
  no.~3, pp. 287--294, 2005.

\bibitem{liu2006symbolic}
X.~Liu, Z.~Lin, and B.~Chen, ``Symbolic realization of asymptotic time-scale
  and eigenstructure assignment design method in multivariable control,''
  \emph{International Journal of Control}, vol.~79, no.~11, pp. 1471--1484,
  2006.

\bibitem{liu2007further}
X.~Liu, Z.~Lin, and B.~M. Chen, ``Further results on structural assignment of
  linear systems via sensor selection,'' \emph{Automatica}, vol.~43, no.~9, pp.
  1631--1639, 2007.

\bibitem{chu2007numerical}
D.~Chu, X.~Liu, and V.~Mehrmann, ``A numerical method for computing the
  hamiltonian schur form,'' \emph{Numerische Mathematik}, vol. 105, no.~3, pp.
  375--412, 2007.

\bibitem{chen2008interconnection}
B.~M. Chen, X.~Liu, and Z.~Lin, ``Interconnection of kronecker canonical form
  and special coordinate basis of multivariable linear systems,'' \emph{Systems
  \& Control Letters}, vol.~57, no.~1, pp. 28--33, 2008.

\bibitem{lisn08}
X.~Liu and Z.~Lin, ``On stabilization of nonlinear systems affine in control,''
  in \emph{Proc. 2008 American Control Conference}, Seattle, WA, Jun. 2008, pp.
  4123--4128.

\bibitem{liu2009assignment}
X.~Liu, Z.~Lin, and B.~M. Chen, ``Assignment of complete structural properties
  of linear systems via sensor selection,'' \emph{IEEE Transactions on
  Automatic Control}, vol.~54, no.~9, pp. 2072--2086, 2009.

\bibitem{liu2009semi3}
X.~Liu and Z.~Lin, ``On semi-global stabilization of minimum phase nonlinear
  systems without vector relative degrees,'' \emph{Science in China. Series F,
  Information sciences}, vol.~52, no.~11, pp. 2153--2162, 2009.

\bibitem{liu2010disturbance}
------, ``Further results on disturbance attenuation for multiple input
  multiple output nonlinear systems,'' in \emph{Proc. 2010 American Control
  Conference}, Baltimore, MD, Jun. 2010.

\bibitem{liu2011normal}
------, ``On normal forms of nonlinear systems affine in control,'' \emph{IEEE
  Transactions on Automatic Control}, vol.~56, no.~2, pp. 239--253, 2011.

\bibitem{liu2012backstepping}
------, ``On the backstepping design procedure for multiple input nonlinear
  systems,'' \emph{International Journal of Robust and Nonlinear Control},
  vol.~22, no.~8, pp. 918--932, 2012.

\bibitem{liu2012grid}
X.~Liu, Z.~Lin, and S.~Acton, ``A grid-based bayesian approach to robust visual
  tracking,'' \emph{Digital Signal Processing}, vol.~22, no.~1, pp. 54--65,
  2012.

\bibitem{liu2013dynamic}
X.~Liu, T.~Iwasaki, and F.~Fish, ``Dynamic modeling and gait analysis of batoid
  swimming,'' in \emph{American Control Conference, 2013. ACC'13.}\hskip 1em
  plus 0.5em minus 0.4em\relax IEEE, 2013.

\bibitem{belcher2014development}
A.~H. Belcher, X.~Liu, Z.~Grelewicz, E.~Pearson, and R.~D. Wiersma,
  ``Development of a 6dof robotic motion phantom for radiation therapy,''
  \emph{Medical physics}, vol.~41, no.~12, p. 121704, 2014.

\bibitem{liu2015robotic}
X.~Liu, A.~H. Belcher, Z.~Grelewicz, and R.~D. Wiersma, ``Robotic real-time
  translational and rotational head motion correction during frameless
  stereotactic radiosurgery,'' \emph{Medical Physics}, vol.~42, no.~6, pp.
  2757--2763, 2015.

\bibitem{belcher2016spatial}
A.~H. Belcher, X.~Liu, Z.~Grelewicz, and R.~D. Wiersma, ``Spatial and
  rotational quality assurance of 6dof patient tracking systems,''
  \emph{Medical Physics}, vol.~43, no.~6, pp. 2785--2793, 2016.

\bibitem{liu2016modeling}
X.~Liu, F.~Fish, R.~S. Russo, S.~S. Blemker, and T.~Iwasaki, ``Modeling and
  optimality analysis of pectoral fin locomotion,'' in \emph{Neuromechanical
  Modeling of Posture and Locomotion}.\hskip 1em plus 0.5em minus 0.4em\relax
  Springer, 2016, pp. 309--332.

\bibitem{liu2017design}
X.~Liu and T.~Iwasaki, ``Design of coupled harmonic oscillators for
  synchronization and coordination,'' \emph{IEEE Transactions on Automatic
  Control}, vol.~62, no.~8, 2017.

\bibitem{liu2017use}
X.~Liu, C.~Pelizzari, A.~H. Belcher, Z.~Grelewicz, and R.~D. Wiersma, ``Use of
  proximal operator graph solver for radiation therapy inverse treatment
  planning,'' \emph{Medical Physics}, vol.~44, no.~4, pp. 1246--1256, 2017.

\bibitem{agrachev2004control}
A.~A. Agrachev and Y.~L. Sachkov, \emph{{Control Theory from the Geometric
  Viewpoint}}.\hskip 1em plus 0.5em minus 0.4em\relax Springer Heidelberg,
  2004.

\bibitem{bullo2005geometric}
F.~Bullo and A.~Lewis, \emph{{Geometric control of mechanical systems}}.\hskip
  1em plus 0.5em minus 0.4em\relax Springer Berlin, 2005.

\bibitem{doolin1990introduction}
B.~F. Doolin and C.~Martin, \emph{{Introduction to differential geometry for
  engineers}}.\hskip 1em plus 0.5em minus 0.4em\relax Marcel Dekker Inc, 1990.

\bibitem{waner2005introduction}
S.~Waner and G.~C. Levine, \emph{{Introduction to Differential Geometry and
  General Relativity }}.\hskip 1em plus 0.5em minus 0.4em\relax Hofstra
  University, 2005.

\bibitem{byaf84}
C.~I. Byrnes and A.~Isidori, ``A frequency domain philosophy for nonlinear
  systems, with applications to stabilization and to adaptive control,'' in
  \emph{Proc. 23rdIEEE Conference on Decision and Control}, Dec. 1984, pp.
  1569--1573.

\bibitem{chgn03}
D.~Cheng and L.~Zhang, ``{Generalized normal form and stabilization of
  non-linear systems},'' \emph{International Journal of Control}, vol.~76,
  no.~2, pp. 116--128, 2003.

\bibitem{conc99}
G.~Conte, C.~H. Moog, and A.~M. Perdon, \emph{Nonlinear Control Systems: an
  Algebraic Setting}.\hskip 1em plus 0.5em minus 0.4em\relax Springer, 1999.

\bibitem{diar07}
Z.~Ding, ``{Asymptotic rejection of asymmetric periodic disturbances in
  output-feedback nonlinear systems},'' \emph{Automatica}, vol.~43, no.~3, pp.
  555--561, 2007.

\bibitem{foga94}
T.~I. Fossen, \emph{{Guidance and Control of Ocean Vehicles}}.\hskip 1em plus
  0.5em minus 0.4em\relax Chichester, 1994.

\bibitem{hiio79m}
R.~M. Hirschorn, ``Invertibility of multivariable nonlinear control systems,''
  \emph{IEEE Transactions on Automatic Control}, vol.~24, pp. 855--865, 1979.

\bibitem{isnd81}
A.~Isidori, A.~J. Krener, C.~Gori-Giorgi, and S.~Monaco, ``Nonlinear decoupling
  via feedback: a differential geometric approach,'' \emph{IEEE Transactions on
  Automatic Control}, vol.~26, pp. 331--345, 1981.

\bibitem{jiau04}
Z.~P. Jiang, I.~Mareels, D.~J. Hill, and J.~Huang, ``{A unifying framework for
  global regulation via nonlinear output feedback: from ISS to iISS},''
  \emph{IEEE Transactions on Automatic Control}, vol.~49, no.~4, pp. 549--562,
  2004.

\bibitem{kane06}
G.~Kaliora, A.~Astolfi, and L.~Praly, ``{Norm estimators and global output
  feedback stabilization of nonlinear systems With ISS inverse dynamics},''
  \emph{IEEE Transactions on Automatic Control}, vol.~51, no.~3, pp. 493--498,
  2006.

\bibitem{kaof05}
D.~Karagiannis, Z.~P. Jiang, R.~Ortega, and A.~Astolfi, ``{Output-feedback
  stabilization of a class of uncertain non-minimum-phase nonlinear systems},''
  \emph{Automatica}, vol.~41, no.~9, pp. 1609--1615, 2005.

\bibitem{khns02}
H.~K. Khalil, \emph{Nonlinear Systems}, 2nd~ed.\hskip 1em plus 0.5em minus
  0.4em\relax Prentice Hall, 2002.

\bibitem{liberzon2004output}
D.~Liberzon, ``{Output--input stability implies feedback stabilization},''
  \emph{Systems \& Control Letters}, vol.~53, no. 3-4, pp. 237--248, 2004.

\bibitem{manc96}
R.~Marino and P.~Tomei, \emph{Nonlinear Control Design: Geometric, Adaptive and
  Robust}.\hskip 1em plus 0.5em minus 0.4em\relax London: Prentice-Hall, 1996.

\bibitem{pees96}
K.~Y. Pettersen and O.~Egeland, ``{Exponential stabilization of an
  underactuated surface vessel},'' in \emph{Proc. 35th IEEE Conference on
  Decision and Control}, vol.~1, 1996, pp. 967--972.

\bibitem{sett02}
R.~Sepulchre, M.~Arcak, and A.~R. Teel, ``{Trading the stability of finite
  zeros for global stabilization ofnonlinear cascade systems},'' \emph{IEEE
  Transactions on Automatic Control}, vol.~47, no.~3, pp. 521--525, 2002.

\bibitem{siam81}
S.~N. Singh, ``A modified algorithm for invertibility in nonlinear systems,''
  \emph{IEEE Transactions on Automatic Control}, vol.~26, pp. 595--598, 1981.

\bibitem{tegs94}
A.~R. Teel and L.~Praly, ``{Global stabilizability and observability imply
  semi-global stabilizability by output feedback},'' \emph{Systems \& Control
  Letters}, vol.~22, no.~5, pp. 313--325, 1994.

\bibitem{esof92}
F.~Esfandiari and H.~K. Khalil, ``{Output feedback stabilization of fully
  linearizable systems},'' \emph{International Journal of Control}, vol.~56,
  no.~5, pp. 1007--1037, 1992.

\bibitem{krna95}
M.~Krsti{\'c}, I.~Kanellakopoulos, and P.~V. Kokotovi{\'c}, \emph{{Nonlinear
  and Adaptive Control Design}}.\hskip 1em plus 0.5em minus 0.4em\relax John
  Wiley \& Sons, New York, 1995.

\bibitem{sags90}
A.~Saberi, P.~V. Kokotovi{\'c}, and H.~J. Sussmann, ``{Global stabilization of
  partially linear composite systems},'' \emph{SIAM Journal of Control and
  Optimization}, vol.~28, no.~6, pp. 1491--1503, 1990.

\bibitem{tetf95}
A.~R. Teel and L.~Praly, ``{Tools for semiglobal stabilization by partial state
  and output feedback},'' \emph{SIAM Journal of Control and Optimization},
  vol.~33, pp. 1443--1488, 1995.

\bibitem{atas99}
A.~N. Atassi and H.~K. Khalil, ``{A separation principle for the stabilization
  of a class of nonlinear systems},'' \emph{IEEE Transactions on Automatic
  Control}, vol.~44, no.~9, 1999.

\bibitem{huot03}
J.~Huang, ``{On the solvability of the regulator equations for a class of
  nonlinear systems},'' \emph{IEEE Transactions on Automatic Control}, vol.~48,
  no.~5, pp. 880--885, 2003.

\bibitem{husi97}
L.~R. Hunt and G.~Meyer, ``{Stable inversion for nonlinear systems},''
  \emph{Automatica}, vol.~33, no.~8, pp. 1549--1554, 1997.

\bibitem{jido98}
Z.~P. Jiang and L.~Praly, ``{Design of robust adaptive controllers for
  nonlinear systems with dynamic uncertainties},'' \emph{Automatica}, vol.~34,
  no.~7, pp. 825--840, 1998.

\bibitem{lios02}
D.~Liberzon, A.~S. Morse, and E.~D. Sontag, ``{Output-input stability and
  minimum-phase nonlinear systems},'' \emph{IEEE Transactions on Automatic
  Control}, vol.~47, no.~3, pp. 422--436, 2002.

\bibitem{oria03}
R.~Ortega, L.~Hsu, and A.~Astolfi, ``{Immersion and invariance adaptive control
  of linear multivariable systems},'' \emph{Systems \& Control Letters},
  vol.~49, no.~1, pp. 37--47, 2003.

\bibitem{saas87}
P.~Sannuti and A.~Saberi, ``A special coordinate basis of multivariable linear
  systems -- finite and infinite zero structure, squaring down and
  decoupling,'' \emph{International Journal of Control}, vol.~45, pp.
  1655--1704, 1987.

\bibitem{ross70}
H.~H. Rosenbrock, \emph{State Space and Multivariable Theory}.\hskip 1em plus
  0.5em minus 0.4em\relax New York: John-Wiley, 1970.

\bibitem{hatf76}
M.~L.~J. Hautus, ``The formal \mbox{L}aplace transform for smooth linear
  systems,'' in \emph{Mathematical Systems Theory, Lect. Notes Econ. Math.
  Syst.}, vol. 131, 1976, pp. 29--47.

\bibitem{isnf83}
A.~Isidori, ``{Nonlinear feedback, structure at infinity and the input-output
  linearization problem},'' in \emph{Mathematical theory of networks and
  systems}.\hskip 1em plus 0.5em minus 0.4em\relax Springer, 1983.

\bibitem{niza85}
H.~Nijmeijer and J.~Schumacher, ``{Zeros at infinity for affine nonlinear
  control systems},'' \emph{IEEE Transactions on Automatic Control}, vol.~30,
  no.~6, pp. 566--573, 1985.

\bibitem{flan86}
M.~Fliess, ``{A new approach to the structure at infinity of nonlinear
  systems},'' \emph{Systems \& Control Letters}, vol.~7, no.~5, pp. 419--421,
  1986.

\bibitem{mond88}
C.~H. Moog, ``{Nonlinear decoupling and structure at infinity},''
  \emph{Mathematics of Control, Signals and Systems}, vol.~1, no.~3, pp.
  257--268, 1988.

\bibitem{beri89}
M.~D.~D. Benedetto, J.~W. Grizzle, and C.~H. Moog, ``Rank invariants of
  nonlinear systems,'' \emph{SIAM Journal of Control and Optimization},
  vol.~27, pp. 658--672, 1989.

\bibitem{respondek1990right}
W.~Respondek, ``Right and left invertibility of nonlinear control systems,'' in
  \emph{Nonlinear Controllability and Optimal Control}, New York and Basel,
  1990, pp. 133--176.

\bibitem{chls04}
B.~M. Chen, Z.~Lin, and Y.~Shamash, \emph{Linear Systems Theory: A Structural
  Decomposition Approach}.\hskip 1em plus 0.5em minus 0.4em\relax Boston:
  Birkh{\"a}user, 2004.

\bibitem{lils04}
\BIBentryALTinterwordspacing
Z.~Lin, B.~M. Chen, and X.~Liu, \emph{Linear Systems Toolkit}, 2004. [Online].
  Available: \url{http://linearsystemskit.net.}
\BIBentrySTDinterwordspacing

\bibitem{brtr65}
R.~W. Brockett and M.~R. Mesarovic, ``The reproducibility of multivariable
  systems,'' \emph{Journal of Mathematics Analysis Application}, vol.~11, pp.
  548--563, 1965.

\bibitem{saio69}
M.~K. Sain and J.~L. Massey, ``Invertibility of linear time-invariant dynamical
  systems,'' \emph{IEEE Transactions on Automatic Control}, vol.~14, pp.
  141--149, 1969.

\bibitem{siio69}
L.~M. Silverman, ``Inversion of multivariable linear systems,'' \emph{IEEE
  Transactions on Automatic Control}, vol.~14, pp. 270--276, 1969.

\bibitem{nijmeuer1988dynamic}
H.~Nijmeijer and W.~Respondek, ``{Dynamic input-output decoupling of nonlinear
  control systems},'' \emph{IEEE Transactions on Automatic Control}, vol.~33,
  no.~11, pp. 1065--1070, 1988.

\bibitem{respondek1988local}
W.~Respondek and H.~Nijmeijer, ``{On local right-invertibility of nonlinear
  control systems},'' \emph{Control Theory and Advanced Technology}, vol.~4,
  pp. 325--348, 1988.

\bibitem{suph91}
H.~J. Sussmann and P.~Kokotovi{\'c}, ``{The peaking phenomenon and the global
  stabilization of nonlinearsystems},'' \emph{IEEE Transactions on Automatic
  Control}, vol.~36, no.~4, pp. 424--440, 1991.

\bibitem{liss94}
Z.~Lin and A.~Saberi, ``{Semi-global stabilization of minimum phase nonlinear
  systems in special normal form via linear high-and-low-gain state
  feedback},'' \emph{International Journal of Robust and Nonlinear Control},
  vol.~4, pp. 353--362, 1994.

\bibitem{tess92}
A.~R. Teel, ``{Semi-global stabilization of minimum phase nonlinear systems in
  special normal forms},'' \emph{Systems \& Control Letters}, vol.~19, no.~3,
  pp. 187--192, 1992.

\bibitem{willems1981ais}
J.~Willems, ``{Almost invariant subspaces: An approach to high gain feedback
  design--Part I: Almost controlled invariant subspaces},'' \emph{IEEE
  Transactions on Automatic Control}, vol.~26, no.~1, pp. 235--252, 1981.

\bibitem{isidori1996nad}
A.~Isidori, ``{A note on almost disturbance decoupling for nonlinear minimum
  phase systems},'' \emph{Systems \& Control Letters}, vol.~27, no.~3, pp.
  191--194, 1996.

\bibitem{marino1989add}
R.~Marino, W.~Respondek, and A.~J. Van~der Schaft, ``{Almost disturbance
  decoupling for single-input single-output nonlinear systems},'' \emph{IEEE
  Transactions on Automatic Control}, vol.~34, no.~9, pp. 1013--1017, 1989.

\bibitem{marino1994nha}
R.~Marino, W.~Respondek, A.~J. Van~der Schaft, and P.~Tomei, ``{Nonlinear
  $H_\infty$ almost disturbance decoupling},'' \emph{Systems \& Control
  Letters}, vol.~23, no.~3, pp. 159--68, 1994.

\bibitem{isidori1996gad}
A.~Isidori, ``{Global almost disturbance decoupling with stability for non
  minimum-phase single-input single-output nonlinear systems},'' \emph{Systems
  \& Control Letters}, vol.~28, no.~2, pp. 115--122, 1996.

\bibitem{lin1998almost}
Z.~Lin, ``{Almost disturbance decoupling with global asymptotic stability for
  nonlinear systems with disturbance-affected unstable zero dynamics},''
  \emph{Systems \& Control Letters}, vol.~33, no.~3, pp. 163--169, 1998.

\end{thebibliography}

%\printindex

%%%%%%%%%%%%%%%%%%%%%%%%%%%%%%%%%%%%%%%%%%%%%%%%%%%%%%%%%%%%%%%%%%%%%%

\end{document}